\documentclass[a4paper]{article}
\title{$P$-adic $L$-functions for GL$_2$}
\date{}

\usepackage[utf8]{inputenc}
\usepackage[all]{xy}

\author{Daniel Barrera Salazar and Chris Williams}

\newcommand{\Addresses}{{
  \bigskip
  \footnotesize

 Daniel Barrera Salazar; Universitat Polit\`{e}cnica de Catalunya\\Campus Nord, Calle Jordi Girona, 1-3, 08034 Barcelona, Spain\par\nopagebreak
\texttt{daniel.barrera.salazar@upc.edu}

  \medskip
Chris Williams; Imperial College London\\ South Kensington Campus\\London SW7 2AZ, United Kingdom\par\nopagebreak
\texttt{christopher.williams@imperial.ac.uk}

}}

\usepackage{etoolbox}
\newtoggle{long}
\togglefalse{long}

\usepackage{comment}
\usepackage{verbatim}
\iftoggle{long}{
\newenvironment{longversion}{}{} 
\newenvironment{shortversion}{\comment}{\endcomment} 

}{

\newenvironment{shortversion}{}{} 
\newenvironment{longversion}{\comment}{\endcomment} 

}

\usepackage[margin=1.35in]{geometry}
\usepackage{preamble}
\usepackage{bbm}

\newcommand{\lb}{\\ \\}

\newcommand{\coker}{\mathrm{coker}}

\newcommand{\inc}{\mathrm{inc}_p}

\newcommand{\clgp}[1]{\mathrm{Cl}_F^+(#1)}

\newcommand{\nclgp}{\mathrm{Cl}_F^+}

\newcommand{\rayclgp}{\mathcal{D}(\clgp{p^\infty},L)}

\newcommand{\dist}{\mathcal{D}_\lambda(L)}
\newcommand{\distplus}{\mathcal{D}_\lambda^+(L)}
\newcommand{\distplusf}{\mathcal{D}_\lambda^{\ff,+}(L)}
\newcommand{\loc}{\mathcal{A}_\lambda(L)}
\newcommand{\locplus}{\mathcal{A}_{\lambda}^+(L)}
\newcommand{\locplusf}{\mathcal{A}_{\lambda}^{\ff,+}(L)}

\newcommand{\poly}{\mathcal{L}_2(V_{\lambda}(L)^*)}
\newcommand{\polyt}{V_{\lambda}(L)^*}

\newcommand{\polyC}{\mathcal{L}_1(V_{\lambda}(\C)^*)}

\newcommand{\oms}{\mathrm{H}_{\mathrm{c}}^q(Y_1(\n),\mathcal{L}_2(\dist))}

\newcommand{\hc}{\mathrm{H}_{\mathrm{c}}^{q}}
\newcommand{\hct}{\mathrm{H}_{\mathrm{c}}}
\newcommand{\hcusp}{\mathrm{H}_{\mathrm{c}}^q}

\newcommand{\epsphi}{\varepsilon_\varphi}

\newcommand{\roiplus}{\roi_{F,+}^\times}

\newcommand{\y}{\mathbf{y}}

\newcommand{\phip}{\varphi_{p-\mathrm{fin}}}

\newcommand{\sho}[1]{\mathcal{L}_1(#1)}
\newcommand{\sht}[1]{\mathcal{L}_2(#1)}
\newcommand{\LL}{\mathcal{L}}
\newcommand{\polyb}{(V_\lambda(\C)^*)}

\newcommand{\kk}{\mathbf{k}}
\newcommand{\vv}{\mathbf{v}}
\newcommand{\jj}{\mathbf{j}}
\newcommand{\XX}{\mathbf{X}}
\newcommand{\YY}{\mathbf{Y}}
\newcommand{\xx}{\mathbf{x}}
\newcommand{\yy}{\mathbf{y}}
\newcommand{\ttt}{\mathbf{t}}
\newcommand{\nn}{\mathbf{n}}
\newcommand{\zz}{\mathbf{z}}

\newcommand{\pmr}{\{\pm 1\}^{\Sigma(\R)}}

\newcommand{\ev}{\mathrm{Ev}}

\newcommand{\evclo}{\ev^{a_{ij}}_{\ff,\jj,1}}

\newcommand{\evy}{\ev^{a_{\y}}_{\ff,\dagger}}

\newcommand{\evyclo}{\ev^{a_{\y}}_{\ff,\jj,1}}
\newcommand{\evyclt}{\ev^{a_{\y}}_{\ff,\jj,2}}

\newcommand{\roifp}{\roi_F\otimes_{\Z}\Zp}

\newcommand{\rr}{\mathbf{r}}

\newcommand{\autoc}{E(\ff)F_\infty^1\backslash F_\infty^+}

\newcommand{\Diff}{\mathfrak{D}}

\begin{document}

\maketitle

%
%

\begin{abstract}Since Rob Pollack and Glenn Stevens used overconvergent modular symbols to construct $p$-adic $L$-functions for non-critical slope rational modular forms, the theory has been extended to construct $p$-adic $L$-functions for non-critical slope automorphic forms over totally real and imaginary quadratic fields by the first and second authors respectively. In this paper, we give an analogous construction over a general number field. In particular, we start by proving a control theorem stating that the specialisation map from overconvergent to classical modular symbols is an isomorphism on the small slope subspace. We then show that if one takes the modular symbol attached to a small slope cuspidal eigenform, then one can construct a ray class distribution from the corresponding overconvergent symbol, that moreover interpolates critical values of the $L$-function of the eigenform. We prove that this distribution is independent of the choices made in its construction. We define the $p$-adic $L$-function of the eigenform to be this distribution.
\end{abstract}

%
%

\section*{Introduction}
The study of $L$-functions has proved extremely fruitful in number theory for almost two centuries, and there are a wealth of research papers relating their critical values to important arithmetic information. A much more recent branch of the theory is the construction and study of $p$-adic $L$-functions, which are natural $p$-adic analogues of classical (complex) $L$-functions. These $p$-adic $L$-functions are naturally distributions on certain ray class groups that interpolate the algebraic parts of critical classical $L$-values. Such $p$-adic $L$-functions have been constructed in a number of cases; for example, one can attach $p$-adic $L$-functions to Dirichlet characters, number fields and rational elliptic curves. Where they exist, these objects have had a number of interesting applications. For example, the \emph{Iwasawa main conjectures} are a wide-ranging series of conjectures predicting deep links between $p$-adic $L$-functions and Selmer groups attached to Galois representations. The Iwasawa main conjecture has been proved by Skinner and Urban for a large class of elliptic curves (see \cite{SU14}). If the main conjecture holds for an elliptic curve $E$, then the order of vanishing of the $p$-adic $L$-function of $E$ is directly related to the rank of the $p$-Selmer group of $E$. Under finiteness of $\Sha(E/\Q)$, this is enough to deduce a $p$-adic analogue of the Birch and Swinnerton-Dyer conjecture (see \cite{MTT86, Dis16} for details of the conjecture). Moreover, the Iwasawa main conjecture has been used to prove the $p$-part of the leading term formula in the (classical) Birch and Swinnerton-Dyer conjecture in analytic ranks $0$ and $1$ (see \cite{JSW15,Cas17, CCSS17}).
$\lb$
The first constructions of $p$-adic $L$-functions for classical modular forms were given by Mazur and Swinnerton-Dyer in \cite{MSD74}, followed by a variety of other constructions. In particular, in 2011, Pollack and Stevens gave an alternative construction using overconvergent modular symbols in \cite{PS11}. Until recently, $p$-adic $L$-functions of automorphic forms for $\GLt$ over more general number fields had been constructed only in isolated cases. For the most general results previously known, see \cite{Haran87} or \cite{Dep16}, where such $p$-adic $L$-functions are constructed for weight 2 (also known as parallel weight $0$) forms that are ordinary at $p$.
$\lb$
Pollack and Stevens' construction of $p$-adic $L$-functions for \emph{small slope} classical modular forms is both beautiful and computationally effective. The first author generalised their approach to the case of Hilbert modular forms in \cite{Bar13}, whilst the second author generalised their approach to Bianchi modular forms (that is, modular forms for $\GLt$ over imaginary quadratic fields) in \cite{Wil17}. These two generalisations use very different methods, owing to the different difficulties that arise in the respective cases. In this paper, we generalise these results further to construct $p$-adic $L$-functions for small-slope automorphic forms for $\GLt$ over any number field.

\subsection*{Summary of the results}
The construction of these $p$-adic $L$-functions is essentially completed  via a blend of the methods used previously by the authors in their respective PhD theses. We now give a quick summary of the argument. Throughout the paper, we take $\Phi$ to be a cohomological cuspidal automorphic eigenform of weight $\lambda$ and level $\Omega_1(\n)$ over a number field $F$, where $\lambda$ and $\Omega_1(\n)$ are defined as in Section \ref{notation}. We write $d = r_1+2r_2$ for the degree of $F$, where $r_1$ (resp.\ $r_2$) denotes the number of real (resp.\ complex) places of $F$.
$\lb$
Let $q = r_1 + r_2$. The space of \emph{modular symbols} of level $\Omega_1(\n)$ and weight $\lambda$ is the compactly supported cohomology space $\hc(Y_1(\n),\mathcal{V}_\lambda)$, where $Y_1(\n)$ is the locally symmetric space associated to $\Omega_1(\n)$ and $\mathcal{V}_\lambda$ is a suitable sheaf of polynomials on $Y_1(\n)$ depending on the weight. The \emph{Eichler-Shimura isomorphism} gives a Hecke-equivariant isomorphism between this cohomology group and the direct sum of certain spaces of automorphic forms. In particular, to each automorphic form $\Phi$ as above -- an inherently analytic object -- one can associate a canonical modular symbol (up to scaling) in a way that preserves Hecke data. In passing from an analytic to an algebraic object, we obtain something that is in some ways easier to study.$\lb$
Using \emph{evaluation maps}, which were described initially by Dimitrov for totally real fields in \cite{Dim13} and which we have generalised to the case of arbitrary number fields, we relate this modular symbol to critical values of the $L$-function of the automorphic form. We show that these results have an algebraic analogue; that is, we can pass to a cohomology class with coefficients in a sufficiently large number field, and then relate this to the algebraic part of the critical $L$-values of $\Phi$.
\begin{longversion} In particular, in the most technical part of the paper, we prove\end{longversion}
\begin{shortversion}In particular, we give a sketch proof of\end{shortversion} the following result (see Theorem \ref{evphithm} in the paper for a more precise formulation):
\begin{mthmnum}
For each Hecke character $\varphi$ of $F$, there is a map 
\[\mathrm{Ev}_\varphi : \hc(Y_1(\n),\mathcal{V}_\lambda(A)) \longrightarrow A\]
such that if $\Phi$ is a cuspidal automorphic form of weight $\lambda$, with associated $A$-valued modular symbol $\phi_A$ (for $A$ either $\C$ or a sufficiently large number field), then
\[\ev_\varphi(\phi_A) = (*)L(\Phi,\varphi),\]
where $L(\Phi,\cdot)$ is the $L$-function attached to $\Phi$ and $(*)$ is an explicit factor.
\end{mthmnum}
All of this is rather classical in nature, and makes explicit results that are, in theory, `well-known' (although the authors could not find the results in the generality they require in the existing literature). At this point, we start using new $p$-adic methods. Henceforth, assume that $(p)|\n$, and take $L$ to be a (sufficiently large) finite extension of $\Qp$. We define the space of \emph{overconvergent modular symbols of level $\Omega_1(\n)$ and weight $\lambda$} to be the compactly supported cohomology of $Y_1(\n)$ with coefficients in an (infinite-dimensional) space of $p$-adic distributions equipped with an action of $\Omega_1(\n)$ that depends on $\lambda$. 
$\lb$ 
For each prime $\pri|p$ in $F$, we have the Hecke operator $U_{\pri}$ at $\pri$ on both automorphic forms and (classical and overconvergent) modular symbols, induced from the action of the matrix $\smallmatrd{1}{0}{0}{\pi_{\pri}}$, where $\pi_{\pri}\in L$ is a fixed uniformiser at $\pri$. There is a natural \emph{specialisation map} from overconvergent to classical modular symbols that is equivariant with respect to these operators. 
$\lb$
In Section \ref{slopedecomp}, we prove that for any $h_{\pri}\in\Q$, the space of overconvergent modular symbols admits a slope $\leq h_{\pri}$ decomposition (as defined in Definition \ref{slope decomposition}) with respect to the $U_{\pri}$ operator. 
\begin{mdef*}
Let $M$ be an $L$-vector space with an action of a collection of operators $\{U_{\pri}:\pri|p\}$. Where it exists, we denote the slope $\leq h_\pri$ subspace with respect to the $U_{\pri}$ operator by $M^{\leq h_{\pri},U_{\pri}}$. If $\mathrm{h} \defeq (h_{\pri})_{\pri|p}$ is a collection of rationals indexed by the primes above $p$, we define
\[M^{\leq \mathbf{h}} \defeq \bigcap_{\pri|p}M^{\leq h_{\pri},U_{\pri}}\]
to be the slope $\leq\mathbf{h}$-subspace at $p$.
\end{mdef*}
\begin{mdef*}
Let $p\roi_F = \prod \pri^{e_{\pri}}$ be the decomposition of $p$ in $F$, and for each $\pri|p$ let $h_{\pri} \in \Q$. Let $\Sigma$ denote the set of complex embeddings of $F$, and write the weight $\lambda$ as $\lambda = ((k_\sigma),(v_\sigma)) \in \Z[\Sigma]^2$. For each $\sigma \in \Sigma$, there is a unique prime $\pri(\sigma)|p$ corresponding to $\sigma$, and to denote this we write $\sigma \sim \pri$. Define $k_{\pri}^0 \defeq \min\{k_\sigma: \sigma\sim\pri\}$ and $v_{\pri}(\lambda) \defeq \sum_{\sigma\sim\pri}v_{\sigma}.$
$\lb$
We say that the slope $\mathbf{h}\defeq (h_{\pri})_{\pri|p}$ is \emph{small} if $h_{\pri} < (k_{\pri}^0 + v_{\pri}(\lambda) + 1)/e_{\pri}$ for each $\pri|p$.
\end{mdef*}
There is a surjective Hecke-equivariant specialisation map $\rho$ from the space of overconvergent modular symbols to the space of classical modular symbols of fixed weight. In Section \ref{controltheorem}, we prove the following \emph{control theorem}:
\begin{mthmnum}
Let $\mathbf{h} \in \Q^{\{\pri|p\}}$ be a small slope. Then the restriction of the specialisation map $\rho$ to the slope $\leq \mathbf{h}$ subspaces of the spaces of modular symbols is an isomorphism.
\end{mthmnum}
In particular, to a small slope cuspidal eigenform -- that is, an eigenform whose associated modular symbol lives in some small-slope subspace of the space of classical modular symbols -- one can attach a \emph{unique} small-slope overconvergent eigenlift of its associated modular symbol.%
$\lb$
Let $\Psi$ be an overconvergent eigensymbol. We can use a slightly different version of the evaluation maps from previously to construct a distribution $\mu_{\Psi}$ on the narrow ray class group $\clgp{p^\infty}$ attached to $\Psi$, closely following the work of the first author in \cite{Bar13}. We prove that the distribution we define is independent of the choice of class group representatives made. Via compatibility between classical and overconvergent evaluation maps, this distribution then interpolates the critical values of the $L$-function of $\Phi$, and we hence define the $p$-adic $L$-function to be this distribution. To summarise, the main result of this paper is:
\begin{mthmnum}
Let $\Phi$ be a small slope cuspidal eigenform over $F$. Let $\phi_{\Phi}$ be the ($p$-adic) classical modular symbol attached to $\Phi$, and let $\Psi_{\Phi}$ be its (unique) small-slope overconvergent eigenlift. Let $\mu_{\Phi}$ be the distribution on $\clgp{p^\infty}$ attached to $\Psi_{\Phi}$.
$\lb$ 
If $\varphi$ is a critical Hecke character, then we can define a canonical locally algebraic character $\varphi_{p-\mathrm{fin}}$ on $\clgp{p^\infty}$ associated to $\varphi$. Then
\[\mu_{\Phi}(\varphi_{p-\mathrm{fin}}) = (*)L(\Phi,\varphi),\]
where $(*)$ is an explicit factor.
\end{mthmnum}
\begin{mdef*}We define the \emph{$p$-adic $L$-function of $\Phi$} to be the distribution $\mu_{\Phi}$ on $\clgp{p^\infty}$.
\end{mdef*}
For a precise notion of which characters are critical, and the factor $(*)$, see Theorem \ref{summary}.
$\lb$
In the case that $F$ is totally real or imaginary quadratic, given slightly tighter conditions on the slope one can prove that the distribution we obtain is \emph{admissible}, that is, it satisfies a growth property that then determines the distribution uniquely. In the general situation, it is rather more difficult to define the correct notion of admissibility; we discuss this further in Section \ref{admissibility}. We instead settle for proving that our construction is independent of choice, so that it is indeed reasonable to \emph{define} the $p$-adic $L$-function in this manner.

\subsection*{Structure of the paper}
Sections \ref{setup} to \ref{algebraic} of the paper focus on the classical side of the theory. The main results of this part of the paper come in Sections \ref{autocyclesec} and \ref{algebraic}, where we relate modular symbols to $L$-values using evaluation maps. Sections \ref{distributions} to \ref{controltheorem} focus on proving the control theorem, allowing us to lift small slope classical eigensymbols to overconvergent symbols. Section \ref{constructdist} then uses evaluation maps to define a distribution attached to an overconvergent eigensymbol. In Section \ref{interpolation}, we prove compatibility results between overconvergent and classical evaluation maps that allow us to prove interpolation properties of this distribution. Our results are summarised fully in Section \ref{padiclfns}.\\
\\
\textbf{Acknowledgements:} We would like to thank David Loeffler for encouraging us to work on this project and for his invaluable comments on the final draft of this paper, as well as for many helpful conversations on the subject. Whilst we were working on this paper, we also met with David Hansen -- who was independently working on a similar project -- to discuss the work we had completed at the time; since submission, he and John Bergdall have released their preprint on the Hilbert case (see \cite{BH17}). The second author would also like to thank Adrian Iovita and the Centre de Recherches Math\'{e}matiques in Montreal for supporting his vist in Spring 2015, during which a large portion of the work in this paper was carried out. Finally, we would like to thank the anonymous referee for their comments and corrections, which greatly improved the paper. \\
\\
The first author was funded by the Centre de Recherches Math\'{e}matiques in Montreal and the European Research Council (ERC) under the European Union's Horizon 2020 research and innovation programme (grant agreement No. 682152), whilst the second author was supported by an EPSRC DTG doctoral grant at the University of Warwick.

%
%

\section{Notation, Hecke characters and automorphic forms}\label{setup}

\subsection{Notation}\label{notation}
This section will serve as an index for the notation that we will use during this paper. Let $p$ be a prime, and fix -- once and for all -- embeddings $\mathrm{inc} : \overline{\Q} \hookrightarrow \C$ and $\inc :\overline{\Q}\hookrightarrow \overline{\Q}_p.$ Let $F$ be a number field of degree $d = r_1 + 2r_2$, where $r_1$ is the number of real embeddings and $r_2$ the number of pairs of complex embeddings of $F$. Write $q = r_1 + r_2$. We write $\Sigma$ for the set of all infinite embeddings of $F$. Let $\Sigma(\R)$ denote the set of real places of $F$ and let $\Sigma(\C)$ be the set containing a (henceforth fixed) choice of embedding from each pair of complex places, so that 
\[\Sigma = \Sigma(\R) \cup \Sigma(\C)\cup c\Sigma(\C),\]
where $c$ denotes complex conjugation. We write $\Diff$ for the different of $F$ and $D$ for the discriminant of $F$. For each finite place $v$ in $F$, fix (once and for all) a uniformiser $\pi_v$ in the completion $F_v$.
$\lb$
Let $\A_F = F_\infty \times \A_F^f$ denote the adele ring of $F$, with infinite adeles $F_\infty \cong F \otimes_{\Q}\R$ and finite adeles $\A_F^f$. Let $\widehat{\roi}_F \cong \widehat{\Z}\otimes_{\Z}\roi_F$ denote the integral (finite) adeles. Let $F_\infty^+ \cong \R_{>0}^{r_1}\times(\C^\times)^{r_2}$ be the connected component of the identity in $F_\infty^\times$.
$\lb$
Let $\n \subset \roi_F$ be an ideal with $(p)|\n$. This will be our level; write
\[\Omega_1(\n) \defeq \left\{\matr \in \GLt\left(\widehat{\roi}_F\right): c \equiv 0 \newmod{\n}, d \equiv 1 \newmod{\n}\right\}.\]
This is an open compact subgroup of $\GLt(\A_F^f).$ Let $K_\infty^+ \defeq \SOt(\R)^{r_1}\times\SUt(\C)^{r_2}$, a subgroup of the standard maximal compact subgroup $K_\infty$ of $\GLt(F_\infty)$, and let $Z_\infty \defeq Z(\GLt(F_\infty)) \cong (F\otimes_{\Q}\R)^\times$ (with $Z_\infty^+$ the connected component of $Z_\infty$ including the identity). Then the locally symmetric space associated to $\Omega_1(\n)$ is 
\[Y_1(\n) \defeq \GLt(F)\backslash \GLt(\A_F) /\Omega_1(\n)K_\infty^+ Z_\infty.\]
For an ideal $\ff \subset \roi_F$, we define $U(\ff)$ to be the set of elements of $\widehat{\roi_F}$ that are congruent to $1 \newmod{\ff}$, and denote the narrow ray class group modulo $\ff$ by 
\[\mathrm{Cl}_F^+(\ff) \defeq F^\times\backslash \A_F^\times/U(\ff)F_\infty^+.\]
When $\ff =\roi_F$, we write simply $\mathrm{Cl}_F^+$ (the narrow class group of $F$). Write $h$ for the narrow class number of $F$ and choose fixed representatives $I_1, ..., I_h$ of the narrow class group, coprime to $\n$, represented by ideles $a_1,...,a_h$, with $(a_i)_v = 1$ for all $v|\n\infty$.
$\lb$
Throughout, $\lambda = (\kk,\vv) \in \Z[\Sigma]^2$, with $\kk \geq 0$, will denote an admissible weight (to be defined in Definition \ref{admissiblew}). If $\rr \in \Q[\Sigma]$ is parallel, then we write $[\rr]$ for the unique rational such that $\rr = [\rr]\ttt$, where $\ttt = (1,...,1) \in \Z[\Sigma]$.
$\lb$
For a ring $A$ and an integer $k$, we define $V_k(A)$ to be the ring of homogeneous polynomials in two variables of degree $k$ over $A$. This has a natural left $\GLt(A)$-action given by
\[\matr\cdot f(X,Y) = f(bY+dX,aY+cX).\]
For $\kk \in \Z[\Sigma]$, we write $V_{\kk}(A) \defeq \bigotimes_{v}V_{k_v}(A).$ This has a natural $\GLt(A)^d$-action induced from that on each component. For $\lambda = (\kk,\vv)$ as above, we also write $V_\lambda(A)$ for the module $V_{\kk}(A)$ with $\GLt(A)^d$ action twisted by $\det^{\vv}$, that is, given by 
\[\gamma\cdot_{\lambda}f(\XX,\YY) = \left(\prod_{v\in\Sigma}\det(\gamma_v)^{v_v}\right) \gamma\cdot f(\XX,\YY), \hspace{12pt} \gamma = (\gamma_v)_{v\in\Sigma} \in \GLt(A)^d.\]

	\begin{shortversion}
	%
	%
	\subsection{Hecke characters}
	A \emph{Hecke character for $F$} is a continuous homomorphism $\varphi: F^\times\backslash \A_F^\times \rightarrow \C^\times$. For a place $v$ of $F$, we write $\varphi_v$ for the restriction of $\varphi$ to $F_v^\times$, where $F_v$ denotes the completion of $F$ at $v$. We will typically write $\ff$ for the conductor of $\varphi$. For an ideal $I \subset \roi_F$, write $\varphi_I\defeq\prod_{v|I}\varphi_v$. We write $\varphi_f \defeq \prod_{v\nmid\infty}\varphi_v$ and $\varphi_\infty \defeq \prod_{v|\infty}\varphi_v$.
	$\lb$
	We can identify a Hecke character $\varphi$ with a function on ideals of $F$ that has support on those that are coprime to the conductor in a natural way. Concretely, if $\mathfrak{q}$ is a prime ideal coprime to the conductor, define $\varphi(\mathfrak{q}) = \varphi(\pi_{\mathfrak{q}})$ (which is independent of the choice of uniformiser $\pi_{\mathfrak{q}}$), and if $\mathfrak{q}$ is not coprime to the conductor, define $\varphi(\mathfrak{q}) = 0$. In an abuse of notation, we also write $\varphi$ for this function.
	
	%
	%
	\subsubsection{Admissible Infinity Types}\label{heckechar}
	Let $\varphi$ be a Hecke character. There is a canonical decomposition $F_\infty^\times = \pmr \times F_\infty^+,$ and we write $\varphi_\infty^+ \defeq \varphi|_{F_\infty^+}$. We say $\varphi$ is \emph{arithmetic} if $\varphi_\infty^+$ takes the form
	\[\zz = (z_v)_{v|\infty} \longmapsto \zz^{\rr} = \prod_{v|\infty}z_v^{r_v}\]
	for some $\rr \in \Z[\Sigma]$, and we say $\rr$ is the \emph{infinity-type} of $\varphi$. Henceforth, all Hecke characters will be assumed to be arithmetic.
	$\lb$
	Define a character $\epsphi$ of the \emph{Weyl group} $\pmr$ attached to $\varphi$ by
	\[\epsphi(\iota) \defeq \varphi_\infty(\iota)\iota^{\rr},\]
	where we consider $\iota\in\pmr$ as an infinite idele by setting its entries at complex places to be 1. In the sequel, we will (in an abuse of notation) write $\epsphi$ for both this character of $\pmr$ and for the character of the ideles given by $\epsphi(x) = \epsphi((\text{sign}(x_v))_{v\in\Sigma(\R)}).$ Note then that $\varphi_\infty \epsphi$ is the unique algebraic character of $F_\infty^\times$ that restricts to $\varphi_\infty^+$ on $F_\infty^+$; namely, it is the character of $F_\infty^\times$ given by $\zz \mapsto \zz^{\rr}$. Note that if $F = \Q$ and $\varphi = |\cdot|$ is the norm character on $\A_{\Q}^\times$, then $\epsphi(-1) = -1$, even though $\varphi$ itself takes only positive values.
\begin{longversion} 
Not all elements of $\Z[\Sigma]$ can be realised as the infinity type of a Hecke character. In \cite{Hid94}, Chapter 3, a description of the `admissible' types are given. In particular, let $F_{\mathrm{CM}}$ be the maximal CM subfield of $F$ (or the maximal totally real subfield if no such CM field exists), and denote the set of infinite places by $\Sigma_{\mathrm{CM}}.$ There is a natural inflation map
	\begin{align*}\mathrm{Inf}: \Z[\Sigma_{\mathrm{CM}}] &\longrightarrow \Z[\Sigma],\\
	\sum_{\tau\in \Sigma_{\mathrm{CM}}}n_\tau \tau &\longmapsto \sum_{\substack{\sigma \in \Sigma\\\sigma|_{F_{\mathrm{CM}}} = \tau}} n_\tau \sigma.
	\end{align*}
	\begin{mdef}
	Let $\Xi_{\mathrm{CM}} \defeq \{\jj \in \Z[\Sigma_{\mathrm{CM}}] : \jj +c\jj \in \ttt_{\mathrm{CM}} \Z\},$ where $\ttt_{\mathrm{CM}} = (1,1,...,1).$ We define the set of \emph{admissible infinity types} to be 
	\[\Xi = \mathrm{Inf}(\Xi_{\CM}).\]
	\end{mdef}
	In more concrete terms, a necessary (but not sufficient) condition for $\rr \in \Xi$ is that $\rr + c\rr$ is parallel. This motivates the following piece of notation, which we will require in the sequel:
\end{longversion}
\begin{shortversion}
Not all elements of $\Z[\Sigma]$ can be realised as the infinity type of a Hecke character. In \cite{Hid94}, Chapter 3, a description of the set $\Xi \subset \Z[\Sigma]$ of `admissible' types are given. A necessary (but not sufficient) condition for $\rr \in \Xi$ is that $\rr + c\rr$ is parallel. This motivates the following piece of notation, which we will require in the sequel:
\end{shortversion}
	\begin{mdef}\label{bracket}
	Let $\rr\in\Z[\Sigma]$ be admissible, that is, let $\rr \in \Xi$. Then define $[\rr] \in \R$ to be the unique number such that
	\[\rr + c\rr = 2[\rr]\ttt.\]
	Note that, in particular, for any $\zeta \in F^\times$, we have $N((\zeta))^{[\rr]} = |\zeta|^{\rr},$ which we will use later.
	\end{mdef}
	
	In \cite{Wei56}, Weil then shows that:
	\begin{mprop}
	An element $\rr \in \Z[\Sigma]$ can be realised as the infinity type of a Hecke character of $F$ if and only if $\rr \in \Xi$, that is, $\rr$ is admissible.
	\end{mprop}
	For example, if $F$ is totally real (or more generally has any real embedding), then the only admissible infinity types are parallel. If $F$ is imaginary quadratic, then \emph{any} pair $(r,s)\in\Z[\Sigma]$ is admissible.
	
	%
	%
	
	\subsubsection{Hecke characters on ray class groups}\label{rcg}
	To a Hecke character $\varphi$ of conductor $\ff|p^\infty$, we can associate a locally analytic function $\phip$ on the $p$-adic analytic group 
	\[\clgp{p^\infty} \defeq F^\times \backslash \A_F^\times/U(p^\infty)F_\infty^+,\]
	where $U(p^\infty)$ is the group of elements of $\widehat{\roi}_F^\times$ that are congruent to 1 $\newmod{p^n}$ for all integers $n$ (that is, elements of $\widehat{\roi}_F^\times$ such that their components at primes above $p$ are all equal to $1$). By class field theory, $\clgp{p^\infty}$ is isomorphic to the Galois group of the maximal abelian extension of $F$ unramified outside $p$ and $\infty$. The $p$-adic $L$-function of an automorphic form over $F$ should be a distribution on this space, and to this end we discuss the structure of this space in the sequel.$\lb$
	Let $\varphi$ be a Hecke character with infinity type $\rr$ and associated character $\epsphi$ on $\pmr$, as above. Then there is a unique algebraic homomorphism $w^{\rr}: F^\times \longrightarrow \overline{\Q}^\times$ given by 
	\[w^{\rr}(\gamma) = \prod_{v \in \Sigma}\sigma_v(\gamma)^{r_v},\]
	where $\sigma_v$ is the complex embedding corresponding to the infinite place $v$. This then induces maps $w^{\rr}_\infty:(F\otimes_{\Q}\R)^\times \rightarrow \C^\times$ and $w^{\rr}_p: (F\otimes_{\Q}\Qp)^\times \rightarrow \overline{\Q}_p^\times \subset \Cp^\times$. Note that $w_\infty^{\rr}$ is equal to $\epsphi\varphi_\infty$, the unique algebraic character of $F_\infty$ that agrees with $\varphi_\infty$ on $F_\infty^+$.
	$\lb$
The finite part of any Hecke character takes algebraic values (see \cite{Wei56}). In particular, under our fixed embedding $\overline{\Q} \hookrightarrow \Cp$, we can see $\varphi_f$ as taking values in $\Cp^\times$. In particular, the following function is well-defined.
	\begin{mdef}We define $\phip$ to be the function
	\begin{align*}\phip : \A_F^\times &\longrightarrow \Cp^\times\\
	x &\longmapsto = \epsphi\varphi_f(x)w_p^{\rr}(x_p).
	\end{align*}
	\end{mdef}
	\begin{mprop}Let $\varphi$ be a Hecke character of conductor $\ff|(p^\infty).$ Then the function $\phip$ gives a well-defined function on the narrow ray class group $\clgp{p^\infty}$.
	\end{mprop}
	\begin{proof}
	By definition, $\phip$ is trivial on $F_\infty^+$. Since $w_\infty^{\rr}$ and $w_p^{\rr}$ are both induced from the same function on $F^\times$, we see that $\phip = \varphi = 1$ on $F^\times$. As $\varphi$ has conductor $\ff$, it is trivial on $U(\ff)$, and hence on $U(p^\infty)$. Finally, if $x \in U(p^\infty)$, then $x_p = x_\infty = 1$, so that $w_p^{\rr}(x_p) = w_\infty^{\rr}(x_\infty) = 1$. This completes the proof.
	\end{proof}

	%
	%
	\subsubsection{Gauss sums}\label{gsum}
	Let $\varphi$ be a Hecke character of conductor $\ff$. We can attach a \emph{Gauss sum} to $\varphi$ that has many of the desirable properties that Gauss sums of Dirichlet characters enjoy. We first introduce a more general exponential map on the adeles of $F$.
	\begin{mdef}\label{ef}
	Let $e_F$ be the unique continuous homomorphism
	\begin{align*}e_F: \A_F/F \longrightarrow \C^\times
	\end{align*}
	that satisfies
	\[
	x_\infty \longmapsto e^{2\pi i\mathrm{Tr}_{F/\Q}(x_\infty)},
	\]
	where $x_\infty$ is an infinite adele. We can describe $e_F$ explicitly as
	\[e_F(\xx) = \prod_{v\in \Sigma(\C)}e^{2\pi i\mathrm{Tr}_{\C/\R}(x_v)}\prod_{v\in\Sigma(\R)}e^{2\pi ix_v}\prod_{\lambda|\ell\text{ finite}}e_\ell(-\mathrm{Tr}_{F_\lambda/\Q_\ell}(x_\lambda)),\]
	where
	\[e_\ell\bigg(\sum_j c_j\ell^j\bigg) = e^{2\pi i \sum_{j<0}c_j\ell^j}.\]
	\end{mdef}
	Let $d$ be a (finite) idele representing the different $\Diff$.
	
	\begin{mdef}\label{gausssum}
	Define the \emph{Gauss sum attached to $\varphi$} to be 
	\[\tau(\varphi) \defeq \varphi(d^{-1})\sum_{b \in (\roi_F/\ff)^\times}\varphi_{\ff}(b)e_F(bd^{-1}(\pi_{\ff}^{-1})_{v|\ff}),\]
	where $(\pi_{\ff}^{-1})_{v|\ff}$ is the adele given by
	\[((\pi_{\ff}^{-1})_{v|\ff})_w \defeq \left\{\begin{array}{ll}\pi_w^{-v_w(\ff)} &: w|\ff\\
	0 &: \text{otherwise}.\end{array}\right.\]
	\end{mdef}
\begin{mrem}\label{gausssumprop}
 This definition, which is independent of the choice of $d$, is a natural one; in fact, it is the product of the $\epsilon$-factors over $v|\ff$, as defined by Deligne in \cite{Del72}. For this particular iteration of the definition, we have followed \cite{Hid94}, page 480 (though we have phrased the definition slightly differently by choosing more explicit representatives).
	\end{mrem}
	\begin{mprop}For $\zeta \in \roi_F$ non-zero, we have
	\[\varphi(d^{-1})\sum_{b \in (\roi_F/\ff)^\times}\varphi_{\ff}(b)e_F(\zeta bd^{-1}(\pi_{\ff}^{-1})_{v|\ff}) =\left\{\begin{array}{ll} \varphi_{\ff}(\zeta)^{-1}\tau(\varphi) &: ((\zeta),\ff) = 1\\
	0 &: \text{otherwise}\end{array}\right.,\]
	where the notation $((\zeta),\ff) = 1$ means that the two ideals are coprime.
	\end{mprop}
	\begin{proof}
	See \cite{Del72}, or, for an English translation, \cite{Tat79}. There is also an account of Gauss sums and their properties in \cite{Nar04}.
	\end{proof}
	
	\end{shortversion}

\begin{longversion}

%
%
\subsection{Hecke characters}

We will extensively refer to Hecke characters in the sequel, and here we recap some of the basic theory, in the process fixing the notation we shall use throughout the paper.
\subsubsection{Definitions}
\begin{mdef}A \emph{Hecke character for $F$} is a continuous homomorphism 
\[\varphi : F^\times \backslash \A_F^\times \longrightarrow \C^\times.\]
By restriction, for each place $v$ of $F$, we obtain a character $\varphi_v : F_v^\times \rightarrow \C^\times$, where $F_v$ denotes the completion of $F$ at $v$.
\end{mdef}
\begin{mdef}
\begin{itemize}
\item[(i)]For each Hecke character $\varphi$, there is a minimal ideal $\ff \subset \roi_F$ such that
\[\varphi|_{1+\ff\widehat{\roi_F}} = 1.\]
We call this ideal the \emph{conductor} of $\varphi$.
\item[(ii)]For an ideal $\ff \subset \roi_F$, write
\[\varphi_{\ff} \defeq \prod_{v|\ff}\varphi_v.\]
\item[(iii)]The \emph{infinite part} of $\varphi$ is
\[\varphi_\infty \defeq \prod_{v|\infty}\varphi_v = \varphi|_{(F\otimes_{\Q}\R)^\times}.\]
The \emph{finite part} of $\varphi$ is
\[\varphi_f \defeq \varphi/\varphi_\infty = \prod_{v\nmid\infty}\varphi_v.\]
\end{itemize}
\end{mdef}
We can identify a Hecke character $\varphi$ with a function on ideals of $F$ that has support on those that are coprime to the conductor in a natural way. In an abuse of notation, we also write $\varphi$ for this function.

%
%
\subsubsection{Admissible Infinity Types}\label{heckechar}
Let $\varphi$ be a Hecke character. There is a canonical decomposition $F_\infty^\times = \pmr \times F_\infty^+,$ and thus any Hecke character $\varphi$ gives a character $\varphi_\infty^+$ on $F_\infty^+$ by restriction. We say such a character of $F_\infty^+$ is \emph{arithmetic} if it takes the form 
\[\zz = (z_v)_{v|\infty} \longmapsto \zz^{\rr} = \prod_{v|\infty}z_v^{r_v}\]
for some $\rr \in \Z[\Sigma]$, and we say $\rr$ is the \emph{infinity-type} of $\varphi$. Henceforth, all Hecke characters will be assumed to be arithmetic.
$\lb$
We also define a character $\epsphi$ of the \emph{Weyl group} $\pmr$ attached to $\varphi$. Rather than defining $\epsphi$ simply by restriction, we do this more subtly. In particular, we can consider $\iota \in \pmr$ as an infinite idele by setting its entries at complex places to be $1$; then we define 
\[\epsphi(\iota) \defeq \varphi_\infty(\iota)\iota^{\rr},\]
where we consider $\iota\in\pmr$ as an infinite idele by setting its entries at complex places to be 1. In the sequel, we will (in an abuse of notation) write $\epsphi$ for both this character of $\pmr$ and for the character of the ideles given by $\epsphi(x) = \epsphi((\text{sign}(x_v))_{v\in\Sigma(\R)}).$ Note then that $\varphi_\infty \epsphi$ is the unique algebraic character of $F_\infty^\times$ that restricts to $\varphi_\infty^+$ on $F_\infty^+$; namely, it is the character of $F_\infty^\times$ given by $\zz \mapsto \zz^{\rr}$. 
\begin{mrem}
Note that if $F = \Q$ and $\varphi = |\cdot|$ is the norm character on $\A_{\Q}^\times$, then $\epsphi(-1) = -1$, even though $\varphi$ itself takes only positive values.
\end{mrem}
Not all elements of $\Z[\Sigma]$ can be realised as the infinity type of a Hecke character. In \cite{Hid94}, Chapter 3, a description of the `admissible' types are given. In particular, let $F_{\mathrm{CM}}$ be the maximal CM subfield of $F$ (or the maximal totally real subfield if no such CM field exists), and denote the set of infinite places by $\Sigma_{\mathrm{CM}}.$ There is a natural inflation map
\begin{align*}\mathrm{Inf}: \Z[\Sigma_{\mathrm{CM}}] &\longrightarrow \Z[\Sigma],\\
\sum_{\tau\in \Sigma_{\mathrm{CM}}}n_\tau \tau &\longmapsto \sum_{\substack{\sigma \in \Sigma\\\sigma|_{F_{\mathrm{CM}}} = \tau}} n_\tau \sigma.
\end{align*}
\begin{mdef}
Let $\Xi_{\mathrm{CM}} \defeq \{\jj \in \Z[\Sigma_{\mathrm{CM}}] : \jj +c\jj \in \ttt_{\mathrm{CM}} \Z\},$ where $\ttt_{\mathrm{CM}} = (1,1,...,1).$ We define the set of \emph{admissible infinity types} to be 
\[\Xi = \mathrm{Inf}(\Xi_{\CM}).\]
\end{mdef}
In more concrete terms, a necessary (but not sufficient) condition for $\rr \in \Xi$ is that $\rr + c\rr$ is parallel. This motivates the following piece of notation, which we will require in the sequel:
\begin{mdef}\label{bracket}
Let $\rr\in\Z[\Sigma]$ be admissible, that is, let $\rr \in \Xi$. Then define $[\rr] \in \R$ to be the unique number such that
\[\rr + c\rr = 2[\rr]\ttt.\]
Note that, in particular, for any $\zeta \in F^\times$, we have $N((\zeta))^{[\rr]} = |\zeta|^{\rr},$ which we will use later.
\end{mdef}

In \cite{Wei56}, Weil then shows that:
\begin{mprop}
An element $\rr \in \Z[\Sigma]$ can be realised as the infinity type of a Hecke character of $F$ if and only if $\rr \in \Xi$, that is, $\rr$ is admissible.
\end{mprop}
This restriction on the possible infinity types comes from the condition that any Hecke character $\varphi$ is trivial on $F^\times$. In particular, an element $\rr \in \Z[\Sigma]$ is admissible if and only if there exists an integer $n$ such that $\epsilon^{n\rr} = \prod_{v\in\Sigma}\epsilon^{nr_v} = 1$ for all $\epsilon \in \roi_F^\times.$
\begin{mexs}\begin{itemize}
\item[(i)] Suppose $F$ is totally real. As the unit group is as `big' as it can be relative to the degree of $F$, this condition is very restrictive, and indeed the only admissible infinity types are parallel.
\item[(ii)] Suppose $F$ is imaginary quadratic. Then the unit group is finite, and if we take $n$ to be its order, we see that \emph{any} element of $\Z[\Sigma]$ can be an infinity type.
\end{itemize}
\end{mexs}

%
%

\subsubsection{Hecke characters on ray class groups}\label{rcg}
To a Hecke character $\varphi$ of conductor $\ff|p^\infty$, we can associate a locally analytic function $\phip$ on the $p$-adic analytic group 
\[\clgp{p^\infty} \defeq F^\times \backslash \A_F^\times/U(p^\infty)F_\infty^+,\]
where $U(p^\infty)$ is the group of elements of $\widehat{\roi_F}$ that are congruent to 1 $\newmod{p^n}$ for all integers $n$ (that is, elements of $\widehat{\roi_F}$ such that their components at primes above $p$ are all equal to $1$). By class field theory, $\clgp{p^\infty}$ is isomorphic to the Galois group of the maximal abelian extension of $F$ unramified outside $p$. The $p$-adic $L$-function of an automorphic form over $F$ should be a distribution on this space, and to this end we discuss the structure of this space in the sequel.$\lb$
Let $\varphi$ be a Hecke character with infinity type $\rr$ and associated character $\epsphi$ on $\pmr$, as above. Then there is a unique algebraic homomorphism $w^{\rr}: F^\times \longrightarrow \overline{\Q}$ given by 
\[w^{\rr}(\gamma) = \prod_{v \in \Sigma}\sigma_v(\gamma)^{r_v},\]
where $\sigma_v$ is the complex embedding corresponding to the infinite place $v$. This then induces maps $w^{\rr}_\infty:(F\otimes_{\Q}\R)^\times \rightarrow \C^\times$ and $w^{\rr}_p: (F\otimes_{\Q}\Qp)^\times \rightarrow \overline{\Qp}^\times \subset \Cp^\times$. Note that $w_\infty^{\rr}$ is equal to $\epsphi\varphi_\infty$, the unique algebraic character of $F_\infty$ that agrees with $\varphi_\infty$ on $F_\infty^+$.
$\lb$
Fix an isomorphism $\C \cong \Cp$ that is compatible with our earlier embedding $\inc: \overline{\Q} \hookrightarrow \overline{\Q}_p.$ Then under this isomorphism, we can see $w_p^{\rr}$ as having values in $\C$. 
\begin{mdef}We define $\phip$ to be the function
\begin{align*}\phip : \A_F^\times &\longrightarrow \C^\times\\
x &\longmapsto w_\infty^{\rr}(x_\infty)^{-1}w_p^{\rr}(x_p)\varphi(x) = \epsphi\varphi_f(x)w_p^{\rr}(x_p).
\end{align*}
\end{mdef}
\begin{mprop}Let $\varphi$ be a Hecke character of conductor $\ff|(p^\infty).$ Then the function $\phip$ gives a well-defined function on the narrow ray class group $\clgp{p^\infty}$.
\end{mprop}
\begin{proof}
By definition, $\phip$ is trivial on $F_\infty^+$. As $w_\infty^{\rr}$ and $w_p^{\rr}$ are both induced from the same function on $F$, we see that $\phip$ is also trivial on $F^\times$. As $\varphi$ has conductor $\ff$, it is trivial on $U(\ff)$, and hence on $U(p^\infty)$. Finally, if $x \in U(p^\infty)$, then $x_p = x_\infty = 1$, so that $w_p^{\rr}(x_p) = w_\infty^{\rr}(x_\infty) = 1$. This completes the proof.
\end{proof}

%
%
\subsubsection{Gauss sums}\label{gsum}
Let $\varphi$ be a Hecke character of conductor $\ff$. We can attach a \emph{Gauss sum} to $\varphi$ that has many of the desirable properties that Gauss sums of Dirichlet characters enjoy. We first introduce a more general exponential map on the adeles of $F$.
\begin{mdef}\label{ef}
Let $e_F$ be the unique function
\begin{align*}e_F: \A_F/F \longrightarrow \C^\times
\end{align*}
that satisfies
\[
x_\infty \longmapsto e^{2\pi i\mathrm{Tr}_{F/\Q}(x_\infty)},
\]
where $x_\infty$ is an infinite adele. We can describe $e_F$ explicitly as
\[e_F(\xx) = \prod_{v\in \Sigma(\C)}e^{2\pi i\mathrm{Tr}_{\C/\R}(x_v)}\prod_{v\in\Sigma(\R)}e^{2\pi ix_v}\prod_{\lambda|\ell\text{ finite}}e_\ell(-\mathrm{Tr}_{F_\lambda/\Q_\ell}(x_\lambda)),\]
where
\[e_\ell\bigg(\sum_j c_j\ell^j\bigg) = e^{2\pi i \sum_{j<0}c_j\ell^j}.\]
\end{mdef}
Let $d$ be a (finite) idele representing the different $\Diff$.

\begin{mdef}\label{gausssum}
Define the \emph{Gauss sum attached to $\varphi$} to be 
\[\tau(\varphi) \defeq \varphi(d^{-1})\sum_{b \in (\roi_F/\ff)^\times}\varphi_{\ff}(b)e_F(bd^{-1}(\pi_{\ff}^{-1})_{v|\ff}),\]
where $(\pi_{\ff}^{-1})_{v|\ff}$ is the adele given by
\[((\pi_{\ff}^{-1})_{v|\ff})_w \defeq \left\{\begin{array}{ll}\pi_w^{-v_w(\ff)} &: w|\ff\\
0 &: \text{otherwise}.\end{array}\right.\]
\end{mdef}
\begin{mrems}\label{gausssumprop}
\begin{itemize}
\item[(i)] This definition is independent of the choice of $d$. 
\item[(ii)] This definition is a natural one; in fact, it is the product of the $\epsilon$-factors over $v|\ff$, as defined by Deligne in \cite{Del72}. For this particular iteration of the definition, we have followed \cite{Hid94}, page 480 (though we have phrased the definition slightly differently by choosing more explicit representatives).
\end{itemize}
\end{mrems}
\begin{mprop}For $\zeta \in \roi_F$ non-zero, we have
\[\varphi(d^{-1})\sum_{b \in (\roi_F/\ff)^\times}\varphi_{\ff}(b)e_F(\zeta bd^{-1}(\pi_{\ff}^{-1})_{v|\ff}) =\left\{\begin{array}{ll} \varphi_{\ff}(\zeta)^{-1}\tau(\varphi) &: ((\zeta),\ff) = 1\\
0 &: \text{otherwise}\end{array}\right.,\]
where the notation $((\zeta),\ff) = 1$ means that the two ideals are coprime.
\end{mprop}
\begin{proof}
See \cite{Del72}, or, for an English translation, \cite{Tat79}. There is also an account of Gauss sums and their properties in \cite{Nar04}.
\end{proof}

\end{longversion}

%
%

\subsection{Automorphic forms}\label{autoforms}
We now give a brief summary of the theory of automorphic forms for $\GLt$, fixing as we do so the notation and conventions we will use during the rest of the paper. For a more comprehensive survey, see \cite{Hid94}, Chapters 2 and 3, or for a more detailed account of the general theory, see \cite{Wei71}.
\begin{mdef}\label{admissiblew}
An element $\lambda = (\kk,\vv) \in \Z[\Sigma]\times\Z[\Sigma]$ is an \emph{admissible weight} if
\begin{itemize}
\item[(i)] we have $\kk =c\kk \geq 0$, and 
\item[(ii)] $\kk + 2\vv$ is parallel.
\end{itemize}
\end{mdef}
\begin{longversion}
We could use a more general notion of admissible by replacing condition (ii) with the condition that $\kk+2\vv \in \Xi$ (as defined above), but this does not affect our results. Note also that condition (ii) forces $k_v \equiv k_w \newmod{2}$ for all $v,w \in \Sigma(\R)$. 
$\lb$
\end{longversion}
Let $\lambda = (\kk,\vv)$ be an admissible weight as above. Recall the definition of $Y_1(\n)$ from Section \ref{notation}; we now define a representation $\rho$ of $K_\infty^+\times Z_\infty$ that will give us the appropriate `weight $\lambda$ automorphy condition'. We do this individually at each place.
\begin{itemize}
\item Suppose $v \in \Sigma(\C).$ Note that for any non-negative integer $n$, the space $V_n(\C)$ (as defined in Section \ref{notation}) is an irreducible right $\SUt(\C)$-module; write 
\[\widetilde{\rho}(n) : \SUt(\C)\longrightarrow \mathrm{GL}(V_n(\C))\]
for the corresponding antihomomorphism. Then define
\begin{align*}
\rho_v : \SUt(\C)\times\C^\times &\longrightarrow \mathrm{GL}(V_{2k_v+2}(\C))\\
(u,z)&\longmapsto \widetilde{\rho}(2k_v+2)(u)|z|^{-k_v-2v_v}.
\end{align*}
\item Suppose $v \in \Sigma(\R)$. Define 
\begin{align*}
\rho_v : \SOt(\R)\times\R^\times &\longrightarrow \C^\times\\
\left(r(\theta),x\right) &\longmapsto e^{ik\theta}x^{-k_v-2v_v},
\end{align*}
where $r(\theta) \defeq \smallmatrd{\cos(\theta)}{-\sin(\theta)}{\sin(\theta)}{\cos(\theta)}.$
\end{itemize}
Define $\kk^* \in \Z[\Sigma]$ by
\[k^*_v \defeq \left\{\begin{array}{ll}2k_v+2 &: v \in \Sigma(\C)\cup c\Sigma(\C),\\
0 &: v \in \Sigma(\R).\end{array}\right.\]
Now define 
\[\rho: K_\infty^+\times Z_\infty \longrightarrow \mathrm{GL}(V_{\kk^*}(\C))\]
by 
\[\rho \defeq \bigotimes_{v\in\Sigma(\C)\cup\Sigma(\R)} \rho_v.\]
\begin{mdef}
We say that a function 
\[\Phi: \GLt(\A_F) \longrightarrow V_{\kk^*}(\C)\]
is a \emph{cusp form of weight $\lambda$ and level $\Omega_1(\n)$} if it satisfies:
\begin{itemize}
\item[(i)](Automorphy condition) $\Phi(zgu) = \Phi(g)\rho(u,z)$ for $u \in K_\infty^+$ and $z \in Z_\infty \cong (F\otimes_{\Q}\R)^\times,$
\item[(ii)] (Level condition) $\Phi$ is right invariant under $\Omega_1(\n)$,
\item[(iii)] $\Phi$ is left invariant under $\GLt(F)$,
\item[(iv)] (Harmonicity/holomorphy condition) If we write $\Phi_\infty$ for the restriction of $\Phi$ to $\GLt(F_\infty^+)$, where $F_\infty^+$ is the connected component of the identity in $F_\infty$, then $\Phi_\infty$ is an eigenfunction of the operators $\delta_v$ for all places $v$, with
\[\delta_v(\Phi_\infty)= \left(\frac{k_v^2}{2}+ k_v\right)\Phi_\infty,\]
where $\delta_v$ is a component of the Casimir operator in the Lie algebra $\mathfrak{s}\mathfrak{l}_2(\C)\otimes_{\R}\F_v$ (see \cite{Hid93}, Section 1.3),
\item[(v)] (Growth condition) Let $B = \left\{\smallmatrd{t}{z}{0}{1} \in \GLt(\A_F)\right\}.$ Then $\Phi$ is $B$\emph{-moderate} in the sense that there exists $N\geq 0$ such that for every compact subset $S$ of $B$, we have
\[\left|\left|\Phi\left[\matrd{t}{z}{0}{1}\right]\right|\right| = O(|t|^N + |t|^{-N})\]
(for any fixed norm $||\cdot||$) uniformly over $\smallmatrd{t}{z}{0}{1} \in S$,
\item[(vi)] (Cuspidal condition) We have
\[\int_{F\backslash\A_F}\Phi(ug) du = 0,\]
where $\A_F \hookrightarrow \GLt(\A_F)$ via $u \mapsto \smallmatrd{1}{u}{0}{1}$, and $du$ is the Lebesgue measure on $\A_F$.
\end{itemize}
We write $S_{\lambda}(\Omega_1(\n))$ for the space of cusp forms of weight $\lambda$ and level $\Omega_1(\n)$.
\end{mdef}

\begin{longversion}
\begin{mrem}
This is a natural generalisation of classical modular forms to the setting of general number fields; an explanation of this is contained in \cite{Wei71}. When $F$ is totally real, the harmonicity condition reduces to a holomorphicity statement; the operators $\delta_v$ will play no further part in this paper.
\end{mrem}

A cusp form $\Phi$ of weight $\lambda$ and level $\Omega_1(\n)$ descends to give a collection of $h$ functions $F^i:\GLt^+(F_\infty) \rightarrow V_{\kk^*}(\C)$, for $i = 1,...,h$, where $h$ is the narrow class number and the superscript `$+$' denotes the connected component of the identity. Such a collection is non-canonical, depending on choices of representatives for the narrow class group of $F$. To obtain this decomposition, we note that $\GLt(\A_F)$ decomposes as the disjoint union of $h$ sets. To describe this, recall that we took $I_1, ... ,I_h$ to be a complete set of representatives for the class group, with idelic representatives $a_i$. Then define
\[g_i = \matrd{a_i}{0}{0}{1} \in \GLt(\A_F).\]
Then we have:
\begin{mthm}[Strong Approximation]\label{strongapprox} There is a decomposition
\[\GLt(\A_F) = \coprod\limits_{i=1}^h \GLt(F)\cdot g_i\cdot \left[\GLt^+(F_\infty)\times \Omega_1(\n)\right].\]
\end{mthm}
\begin{proof}
See \cite{Hid94}, equation (3.4b).
\end{proof}
It is now clear that $\Phi$ descends in the way claimed above via $F^i(g) \defeq \Phi(g_i g)$. Such functions are \emph{automorphic forms on $\GLt^+(F_\infty)$ of weight $\lambda$ and level }
\begin{equation}\label{gammai}\Gamma_1^i(\n) \defeq \SLt(F) \cap g_i\Omega_1(\n)\GLt^+(F_\infty)g_i^{-1}.
\end{equation}
This process, when $F$ is an imaginary quadratic field, is described in \cite{Byg98}. Note here that if we take $a_1 = 1$, then $\Gamma_1^1(\n)$ is nothing other than the usual $\Gamma_1(\n) \leq \SLt(\roi_F)$. We can go even further; define
\begin{equation}\label{uhp}\uhp_F \defeq \uhp^{\Sigma(\R)}\times\uhs^{\Sigma(\C)},\end{equation}
where $\uhp$ is the standard upper half-plane and $\uhs \defeq \{(z,t) \in \C\times\R_{>0}\}$. Then we can descend each $F^i$ to a function $\f^i$ on $\uhp_F$. In the case where $F = \Q$, after renormalising, this recovers the standard theory of modular forms.
\end{longversion}
There is a good theory of Hecke operators on the space of automorphic forms, indexed by ideals of $\roi_F$ and given by double coset operators. We do not go into details here; see \cite{Wei71}, Chapter VI, or \cite{Hid88}, Section 2. Many of the nice properties that Hecke operators satisfy for classical modular forms, such as algebraicity of Hecke eigenvalues, also hold in the general case. By a \emph{Hecke eigenform} we mean an eigenvector of all of the Hecke operators.

%
%

\section{$L$-functions and periods}\label{periods}
In the following section, we attach $L$-functions to automorphic forms, and state some algebraicity results for their critical values.

\begin{longversion}
\subsection{$L$-functions and Fourier expansions}
\end{longversion}
\begin{shortversion}
\subsection{$L$-functions}
\end{shortversion}
Let $\Phi$ be a cuspidal eigenform over $F$ of weight $\lambda = (\mathbf{k},\mathbf{v})$ and level $\Omega_1(\n)$, with $T_{I}$-eigenvalue $\lambda_I$ for each non-zero ideal $I \subset \roi_F$.
\begin{mdef}
Let $\varphi$ be a Hecke character of $F$. The $L$-function of $\Phi$ twisted by $\varphi$ is defined to be
\[L(\Phi,\varphi,s) = \sum_{0\neq I \subset \roi_F}\lambda_I \varphi(I)N(I)^{-s}, \hsp s \in \C.\]
This converges absolutely for Re$(s) >> 0 $ (see \cite{Wei71}, Chapter II). In fact, one can show that it has an analytic continuation to all of $\C$ by writing down an integral formula for $L(\Phi,\varphi,s)$. With this analytic continuation taken as a given, we also define
\[L(\Phi,\varphi) \defeq L(\Phi,\varphi,1).\]
\end{mdef}
\begin{longversion}
As in the rational and Bianchi cases, we can make sense of this in terms of Fourier expansions. To use a nice description of the Fourier expansion of $\Phi$, as in \cite{Hid94}, we pass to the unitisation of $\Phi$.
\begin{mdef}
The \emph{unitisation of $\Phi$} is defined to be
\[\Phi^u(g) = |\det(g)|^{[\kk+2\vv]/2}\Phi(g).\]
\end{mdef}
As before, we obtain individual components $F^{i,u}$ on $\GLt(F_\infty^+)$ by setting $F^{i,u}(g) = \Phi^u(g_ig)$. Hida gives a description of the Fourier expansion of such a function.
\begin{mdef}
Let $K_\alpha$ denote the modified $K$-Bessel function of order $\alpha$ (see \cite{Hid94}, Section 6). 
\end{mdef}
Over $\Q$, a rational modular form has a $q$-expansion in terms of $q = e^{2\pi i z} = e^{2\pi i x}e^{-2\pi y}.$ Recall the definition of the additive adelic character $e_F$ in Definition \ref{ef}; we will use this to write down a general analogue of the $e^{2\pi x}$ term. Finally, we combine the theory of Fourier expansions of rational and Bianchi modular forms to write down the final component we will require, namely, \emph{Whittaker functions}, an analogue of the $e^{-2\pi y}$ term in the rational Fourier expansion.
\begin{mdef}
\begin{itemize}
\item[(i)] For $v \in \Sigma(\C)$, define the Whittaker function at $v$ to be
\begin{align*}
W_v : \C^\times &\longrightarrow V_{k_v^*}(\C),\\
y_v &\longmapsto \sum_{m = 0}^{k_v^*}i^{m-k_{cv}-1}\binomc{k_v^*}{m}y_v^{k_v+1-m}|y_v|^{1+m}\times\\& \hspace{90pt}K_{m-k_{cv}-1}(4\pi|y_v|)S_v^{k_v^*-m}T_v^m.
\end{align*}
\item[(ii)] For $v \in \Sigma(\R)$, define the Whittaker function at $v$ to be
\begin{align*}
W_v: \R^\times &\longrightarrow\R,\\
y_v &\longmapsto |y_v|^{1+k_v/2}e^{-2\pi|y_v|}.
\end{align*}
\item[(iii)] Define the Whittaker function to be
\begin{align*}\mathbf{W}_{\kk,\vv} : F_\infty^+&\longrightarrow V_{\kk^*}(\C),\\
\yy &\longmapsto \prod_{v \in \Sigma(\C)\sqcup\Sigma(\R)}W_v(y_v).
\end{align*}
\end{itemize}
\end{mdef}
Note that the components are what Hida calls $W_\sigma^u$ in equation (6.2a) of \cite{Hid94}.
\begin{mthm}\label{fourexp}
There exists a function $c(\cdot,\Phi)$ on fractional ideals of $F$, supported only on integral ideals, such that
\[\Phi^{u}\left(y_\infty^{-1/2}\matrd{a_i \yy}{\xx}{0}{1}\right) = |D|^{1+[\kk+2\vv]/2}\sum_{\zeta \in F_{+,J}^\times} c(\zeta a_i\Diff,\Phi)\mathbf{W}_{\kk,\vv}(\zeta y_\infty)e_F(\zeta \xx),\]
where $y_\infty$ is the infinite part of $\yy$ and is assumed to be totally positive, $D$ is the discriminant of $F/\Q$, the exponent of $|D|$ is defined by $\kk+2\vv = [\kk+2\vv]\ttt$ and $F_{+,J}^\times \defeq \{\zeta \in F^\times: \text{ for }\sigma \in \Sigma(\R), \hspace{3pt}\sigma(\zeta) > 0 \iff \sigma \in J\}.$
\end{mthm}
\begin{proof}
This is completed in \cite{Hid94}, Chapter 6, using results from \cite{Wei71}. This particular iteration is stated in Chapter 7.
\end{proof}
\begin{mdef}
Let $\varphi$ be a Hecke character as above. Define 
\[\mathcal{L}(\Phi,\varphi,s) \defeq \sum_{0\neq I\subset \roi_F}c(I,\Phi)\varphi(I)N(I)^{-s}.\]
\end{mdef}
\begin{mprop}\label{lfncomp}
Suppose $\Phi$ is an eigenform for all the Hecke operators. Then $\lambda_I = c(I,\Phi)N(I)^{[\kk+\vv]+1}$, and thus we have an identity
\[L(\Phi,\varphi,s) = \mathcal{L}(\Phi,\varphi, s-[\kk/2 + \vv]-1).\]
\end{mprop}
\begin{proof}
See \cite{Hid94}, Section 7, equation (L2), where this is given for $\mathrm{Re}(s)>>0$. We have written $L(\Phi,\varphi,s)$ for what Hida calls $L(s,\lambda\otimes\omega)$, and normalised the function he calls $L(s,f,\omega)$ slightly differently.
\end{proof}
\begin{mrem}
The $L$-function we are truly interested in is $L(\Phi,\varphi,s)$; it is this that comes from the algebraic Hecke eigenvalues and the theory of motives, with the `correct' range of critical values satisfying the algebraicity properties decribed below. The function $\mathcal{L}$ is a more analytic object, normalised so that the functional equation takes a nicer form. We will later link modular symbols to the critical values of $L(\Phi,\varphi,s)$ via the above Fourier expansions and their link to $\mathcal{L}$.
\end{mrem}
\end{longversion}
\begin{shortversion}
As in the case of classical normalised eigenforms, we can make sense of this $L$-function in terms of Fourier coefficients, for a suitable Fourier expansion of $\Phi$. For details of this approach see \cite{Hid94}, Section 6. 
\end{shortversion}
We make one more definition for convenience. The $L$-function has been built using local data at finite primes; here we `complete' it by adding in Deligne's $\Gamma$-factors at infinity.
\begin{mdef}\label{lambda}
Let
\[\Lambda(\Phi,\varphi) \defeq \left[\prod_{v\in\Sigma}\frac{\Gamma(j_v+1)}{(-2\pi i )^{j_v + 1}}\right]L(\Phi,\varphi),\]
where $\varphi$ has infinity type $\jj+\vv$. 
\end{mdef}
%
%
\subsection{Periods and algebraicity}\label{periodssection}
To $p$-adically interpolate $L$-values, we need to renormalise so that they are algebraic. The following is a result proved by Hida in \cite{Hid94}, Theorem 8.1. Earlier, Shimura proved this result over $\Q$ in \cite{Shi77} and later over totally real fields in \cite{Shi78}.

\begin{mthm}\label{deligneperiods}
Let $\Phi$ be a cuspidal eigenform over $F$ of weight $\lambda = (\kk,\vv)$ and level $\Omega_1(\n)$, with associated $L$-function $L(\Phi,\cdot).$ Let $\varphi$ be a Hecke character of infinity type $\jj+\vv$, where $0 \leq \jj \leq \kk$, and let $\varepsilon = \epsphi$ be its associated character on $\pmr$ (as in Section \ref{heckechar}). Let $K$ be a number field containing the normal closure of $F$ and the Hecke eigenvalues of $\Phi$. Then there is a period
\[\Omega_\Phi^\varepsilon \in \C^\times,\]
depending only on $\Phi$ and $\varepsilon$, such that
\[\frac{\Lambda(\Phi,\varphi)}{\Omega_{\Phi}^{\varepsilon}\tau(\varphi)} \in K(\varphi),\]
where $K(\varphi)$ is the number field generated over $K$ by adjoining the values of $\varphi$.
\end{mthm}
\begin{mrems}
\begin{itemize}
\item[(i)]We are assuming that all Hecke characters are arithmetic; if we dropped this assumption, then $K(\varphi)$ need not be finite over $K$ (see \cite{Hid94}, Section 8).
\item[(ii)] There are many choices of such a period, differing by elements of $K^\times$. Throughout the rest of the paper, we shall assume that we fix a period for each character $\varepsilon$.
\item[(iii)] Note that the period depends on the character $\epsphi(\iota) \defeq \varphi|_{\pmr}(\iota)\iota^{\jj+\vv}$ of the Weyl group, and \emph{not} the character $\varphi|_{\pmr}$. \end{itemize}
\end{mrems}
Thus we have a collection of $2^{r_1}$ periods attached to $\Phi$, and each corresponds to a different collection of $L$-values, depending on the parity of the corresponding Hecke characters.

%
%
%
%

\section{Classical modular symbols}\label{seccms}
Modular symbols are algebraic objects attached to automorphic forms that retain Hecke data. As we discard analytic conditions, they are frequently easier to work with than automorphic forms themselves. In this section, we give a brief description of how one associates a $p$-adic modular symbol to an automorphic form. We start with an essential piece of notation.

\begin{mdef}
Let $A$ be a ring. Define $V_\lambda(A)^* \defeq \Hom(V_\lambda(A),A)$ to be the topological dual of the of weight $\lambda$ polynomials over $A$. This inherits a right action of $\GLt(A)^d$ via $(P|\gamma)(f) = P(\gamma\cdot f).$
\end{mdef}

%
%

\subsection{Local systems}\label{localsystems}
We will need to study the interplay between complex and $p$-adic coefficients. We give two ways of defining local systems on $Y_1(\n)$.
\begin{mdef}\label{locsys}For all modules $M$ below, we suppose that the centre of $\GLt(F)\cap\Omega_1(\n)$, which is isomorphic to $\{\varepsilon \in \roi_{F}^\times: \varepsilon \equiv 1 \newmod{\n}\}$, acts trivially on $M$. If this were not the case, the following local systems would not be well-defined.
\begin{itemize}
\item[(i)]Suppose $M$ is a right $\GLt(F)$-module. Then define $\sho{M}$ to be the locally constant sheaf on $Y_1(\n)$ given by the fibres of the projection
\[\GLt(F)\backslash (\GLt(\A_F)\times M)/\Omega_1(\n)K_\infty^+ Z_\infty \longrightarrow Y_1(\n),\]
where the action is given by
\[\gamma(g,m)ukz = (\gamma gukz,m|\gamma^{-1}).\]
\item[(ii)]Suppose $M$ is a right $\Omega_1(\n)$-module. Then define $\sht{M}$ to be the locally constant sheaf on $Y_1(\n)$ given by the fibres of the projection
\[\GLt(F)\backslash (\GLt(\A_F)\times M)/\Omega_1(\n)K_\infty^+ Z_\infty \longrightarrow Y_1(\n),\]
where the action is given by
\[\gamma(g,m)ukz = (\gamma gukz,m|u).\]
\end{itemize}
\end{mdef}

\begin{mremsnum}\label{locsysiso}\begin{itemize}
\item[(i)]Note that if $M$ is a right $\GLt(F\otimes_{\Q}\R)$-module or a right $\GLt(F\otimes_{\Q}\Qp)$-module, then $M$ can be given a $\GLt(F)$-module structure by restriction in the natural way, giving a sheaf $\sho{M}$ as in (i) above.
\item[(ii)] Similarly, for any right $\GLt(F\otimes_{\Q}\Qp)$-module, we have an action of $\Omega_1(\n)$ on $M$ via the projection $\mathrm{Pr}: \GLt(\A_F) \rightarrow \GLt(F\otimes_{\Q}\Qp)$, and we get a sheaf $\sht{M}$ as above. In this case, the sheaves $\sho{M}$ and $\sht{M}$ are naturally isomorphic via the map
\[(g,m) \longmapsto (g,m|g_p)\]
of local systems, where $g_p$ is the image of $g$ under the map $\mathrm{Pr}$ above. 
\item[(iii)] Note that, for a number field $K$ containing the normal closure of $F$, the space $V_\lambda(K)^*$ is naturally a $\GLt(F)$-module via the embedding of $\GLt(F)$ in $\GLt(F\otimes_{\Q} \R)$, whilst if $L/\Qp$ is a finite extension containing $\inc(K)$, then $V_\lambda(L)^*$ is naturally a $\GLt(F\otimes_{\Q}\Qp)$-module. So our above comments apply and we get sheaves attached to $V_\lambda(A)^*$ for suitable $A$.
\end{itemize}
\end{mremsnum}
It will usually be clear which sheaf we must take. However, when the coefficient system is $V_\lambda(L)^*$ (for a sufficiently large finite extension $L/\Qp$) we can associate two \emph{different} (though isomorphic) local systems. As we will later (in Lemma \ref{trackiso}) need to keep track of precisely what this isomorphism does to cohomology elements, throughout the paper we will retain the subscript for clarity.

%
%

\subsection{Operators on cohomology groups}
\subsubsection{Hecke operators}
Recall $q \defeq r_1 + r_2$. We can define actions of the Hecke operators on the cohomology groups $\hcusp(Y_1(\n),V_\lambda(A)^*)$. This is described fully in \cite{Hid88}, Chapter 7, pages 346--347, and \cite{Dim05}, Section 1.14, page 518. We give a very brief description of the definition, following Dimitrov.
$\lb$
For each prime ideal $\pri$ of $\roi_F$, we have a Hecke operator $T_{\pri}$ induced by the double coset $[\Omega_1(\n)a_{\pri}\Omega_1(\n)]$, where $a_{\pri} \in \GLt(\A_F)$ is defined by
\[(a_{\pri})_v = \left\{\begin{array}{ll}\smallmatrd{1}{0}{0}{\pi_{\pri}} &: v = \pri\\
\smallmatrd{1}{0}{0}{1} &: \text{otherwise}.\end{array}\right.\] 
When $\pri|\n$ we write $U_{\pri}$ in place of $T_{\pri}$ in the usual manner.

\subsubsection{Action of the Weyl group}\label{weyl}
We also have an action of the Weyl group $\pmr$ on the cohomology, again described by Dimitrov. To describe this, recall that we took $I_1, ... ,I_h$ to be a complete set of representatives for the class group, with idelic representatives $a_i$, and define
\[g_i = \matrd{a_i}{0}{0}{1} \in \GLt(\A_F).\]
Note that via strong approximation (see \cite{Hid94}, equation (3.4b)), there is a decomposition
\begin{equation}\label{modcurvedec}
Y_1(\n) = \bigsqcup_{i=1}^h Y_1^i(\n),
\end{equation}
where
\begin{align*}
Y_1^i(\n) &= \GLt(F)\backslash\GLt(F)g_i\Omega_1(\n)\GLt^+(F_\infty)/\Omega_1(\n)K_\infty^+Z_\infty^+\\
&= \Gamma_1^i(\n)\backslash \uhp_F.
\end{align*}
Here 
\begin{equation}\label{gammai}\Gamma_1^i(\n) \defeq \SLt(F) \cap g_i\Omega_1(\n)\GLt^+(F_\infty)g_i^{-1}
\end{equation}
and $\uhp_F \defeq \uhp^{\Sigma(\R)}\times\uhs^{\Sigma(\C)},$ where $\uhp$ is the standard upper half-plane and $\uhs \defeq \{(z,t) \in \C\times\R_{>0}\}$ is the upper half-space. 
$\lb$
Now, let $\iota = (\iota_v)_{v\in\Sigma(\R)} \in \pmr.$ Then $\iota$ acts on $\uhp_F$ by $\iota\cdot\mathbf{z} = [(\iota_v\cdot z_v)_{v\in\Sigma(\R)},(z_v)_{v\in\Sigma(\C)}]$, where for $v \in \Sigma(\R)$ we define
\[\iota_v\cdot z_v \defeq \left\{\begin{array}{ll}z_v &: \iota_v = 1\\
-\overline{z_v} &: \iota_v = -1.\end{array}\right.\]
This action induces an action of $\pmr$ on $Y_1^i(\n)$ for each $i$ and hence on $Y_1(\n)$. The action of $\pmr$ on $\hcusp(Y_1(\n),\polyC)$ is then induced by the map of local systems
\[\iota\cdot(g,P) \longmapsto (\iota\cdot g,P).\]
We write this action on the right by $\phi \mapsto \phi|\iota$. The actions of the Hecke operators and the Weyl group commute.

%
%
\subsection{The Eichler--Shimura isomorphism}
The major step in the construction of a modular symbol attached to an automorphic form is the \emph{Eichler--Shimura isomorphism}. \begin{mthm}[Eichler--Shimura]\label{ES}
There is a Hecke-equivariant injection
\[S_{\lambda}(\Omega_1(\n)) \hookrightarrow \hcusp\left(Y_1(\n),\sho{V_\lambda(\C)^*}\right).\]
\end{mthm}
\begin{proof}
An explicit recipe is given in \cite{Hid94}. Note we have composed the classical version of the theorem with the canonical inclusion of cuspidal into compactly supported cohomology.
\end{proof}

\begin{longversion}
\begin{mrems}\begin{itemize}
\item[(i)]In fact, one can define \emph{automorphic forms of type $J$} for a subset $J \subset \Sigma(\R)$, and then if we replace the left-hand side with $\oplus_{J\subset\Sigma(\R)}S_{\lambda,J}(\Omega_1(\n))$, this becomes an isomorphism. In general, an automorphic form of type $J$ satisfies a holomorphicity condition at the places in $J$ and an anti-holomorphicity condition at the remaining real places. In the case where $F =\Q$, the case where $J = \Sigma(\R)$ defines the usual theory of modular forms, whilst if $J = \varnothing$, we get the theory of anti-holomorphic modular forms. We stay exclusively with the case $J = \Sigma(\R)$ for simplicity, but the results should carry over to more general $J$ with only minor modification.
\item[(ii)] In this more general case, there is also a natural action of the Weyl group on the direct sum $\oplus_{J\subset\Sigma(\R)}S_{\lambda,J}(\Omega_1(\n))$, and it permutes the factors in a natural way. The isomorphism is also equivariant with respect to this action.
\end{itemize}
\end{mrems}
\end{longversion}

Under the decomposition of equation (\ref{modcurvedec}), we see that for sufficiently large extensions $A$ of $\Q$ or $\Qp$, there is a (non-canonical) decomposition
\begin{equation}\label{decompco}
\hcusp\left(Y_1(\n),\sho{V_\lambda(A)^*}\right) \cong \bigoplus_{i=1}^h \hcusp\left(Y_1^i(\n),\sho{V_\lambda(A)^*}\right).
\end{equation}

%
%
\subsection{Modular symbols}\label{attach}
Let $L/\Qp$ be a finite extension.
\begin{mdef}
The space of \emph{modular symbols of weight $\lambda$ and level $\Omega_1(\n)$ with values in $L$} is the compactly supported cohomology space $\hc(Y_1(\n),\sht{V_\lambda(L)^*})$.
\end{mdef}
\begin{longversion}
\begin{mrem}In the rational and imaginary quadratic cases, that is for the cases when $q = 1$, there is a simpler description in terms of functions on cusps of $\uhp_F$. See \cite{PS11} and \cite{Wil17} for details respectively.
\end{mrem}
\end{longversion}
Let $\Phi \in S_{\lambda}(\Omega_1(\n))$ be a Hecke eigenform. Then via Theorem \ref{ES} we can attach to $\Phi$ an element 
\[\phi_{\C} \in \hcusp\left(Y_1(\n),\sho{V_{\lambda}(\C)^*}\right).\]
We want to pass from a cohomology class with complex coefficients to one with $p$-adic coefficients. To do this, we use the theory of periods described earlier in Section \ref{periodssection}.
\begin{mdef}
Let $\varepsilon$ be a character of the Weyl group $\{\pm 1\}^{\Sigma(\R)}.$ Then define 
\[\hcusp\left(Y_1(\n),\sho{V_{\lambda}(\C)^*}\right)[\varepsilon] \subset \hcusp\left(Y_1(\n),\sho{V_{\lambda}(\C)^*}\right)\] 
to be the subspace on which $\{\pm 1\}^{\Sigma(\R)}$ acts by $\varepsilon$.
\end{mdef}
\begin{mprop}\label{period}
Let $K$ be a number field containing the normal closure of $F$ and the Hecke eigenvalues of $\Phi$, and let $\varepsilon$ be as above. Let $\Omega_\Phi^{\varepsilon}$ be the period appearing in Theorem \ref{deligneperiods}. Define
\[\phi_{\C}^{\varepsilon} \defeq 2^{-r_1}\sum_{\iota \in \{\pm 1\}^{\Sigma(\R)}}\varepsilon(\iota)\phi_\C|\iota.\]
Then $\phi_{\C}^\varepsilon \in \hcusp(Y_1(\n),\sho{V_\lambda(\C)^*})[\varepsilon]$, and 
\[\phi_K^{\varepsilon} \defeq \phi_{\C}^{\varepsilon}/\Omega_{\Phi}^{\varepsilon} \in \hcusp\left(Y_1(\n),\sho{V_\lambda(K)^*}\right)[\varepsilon].\]
\end{mprop}
\begin{proof}
See \cite{Hid94}, Chapter 8.
\end{proof}
\begin{mdef}\label{thetak}Define
\[\theta_K \defeq \sum_{\varepsilon}\phi_K^\varepsilon \in \hcusp(Y_1(\n),\sho{V_\lambda(K)^*}),\]
where the sum is over all possible characters of the Weyl group $\pmr$.
\end{mdef}
Now let $L/\Qp$ be a finite extension containing $\inc(K)$ (for our fixed embedding $\inc:\overline{\Q} \hookrightarrow \overline{\Qp}$). Then $\inc$ induces an inclusion 
\begin{equation}\label{inclusion1}\hcusp\left(Y_1(\n),\sho{V_\lambda(K)^*}\right) \longhookrightarrow \hcusp\left(Y_1(\n),\sho{V_\lambda(L)^*}\right) \cong \hcusp(Y_1(\n),\sht{V_\lambda(L)^*}).
\end{equation}
Finally, there is a canonical inclusion
\begin{equation}\label{inclusion2}
\hcusp\left(Y_1(\n),\sht{V_\lambda(L)^*}\right) \longhookrightarrow \hc\left(Y_1(\n),\sht{V_\lambda(L)^*}\right).
\end{equation}
\begin{mdef}\label{modsymb}
Let $\Phi$ be an eigenform of weight $\lambda$ and level $\Omega_1(\n)$, and let $L$ be as above. The \emph{modular symbol attached to $\Phi$ with values in $L$} is the image
\[\theta_L \in \hc\left(Y_1(\n),\sht{V_\lambda(L)^*}\right)\]
of the symbol $\theta_K$ under the inclusion of equations (\ref{inclusion1}) and (\ref{inclusion2}).
\end{mdef}

\begin{longversion}
\begin{mrem}
To give some brief motivation for this definition, we will later define an evaluation map $\hc(Y_1(\n),\sht{V_\lambda(L)^*}) \rightarrow L$ corresponding to a critical character $\varphi$ such that if $\varphi$ corresponds to the character $\varepsilon_1$ of $\pmr$, the image of $\phi_{L}^{\varepsilon_2}$ gives the algebraic part of the critical $L$-value at $\varphi$ if $\varepsilon_1 = \varepsilon_2$ and vanishes otherwise. So by taking the sum, we allow ourselves to see \emph{all} critical values.
\end{mrem}
\end{longversion}

%
%
%
%

\section{Automorphic cycles, evaluation maps and $L$-values}\label{autocyclesec}
Let $\Phi$ be a cuspidal automorphic form over $F$. In this section, we give a connection between the cohomology class $\phi_{\C}$ associated to $\Phi$ via the Eichler-Shimura isomorphism and critical values of its $L$-function. We do so via \emph{automorphic cycles}. The cycles we define here are a generalisation of the objects Dimitrov uses in \cite{Dim13} in the totally real case. As a consequence of this section, we also get an integral formula for the $L$-function of $\Phi$, generalising the results of \cite{Hid94}, Section 7, where such a formula is obtained for Hecke characters with trivial conductor.

\subsection{Automorphic cycles}\label{autocycles}
Let $\ff$ be an integral ideal of $F$. We begin with some essential definitions:
\begin{mdef}
Recall $F_\infty^+ \subset (F\otimes_{\Q}\R)^\times$ is the connected component of the identity in the subgroup of infinite ideles, and let $F_\infty^1$ be the subset defined by
\[F_\infty^1 \defeq \{x \in F_\infty^+: |x_v|_v = 1 \text{ for all }v|\infty.\}.\]
\end{mdef}
\begin{mdef}\label{autocyclesdef}
\begin{itemize}
\item[(i)] Recall the definition of $U(\ff) \subset \A_{F,f}^\times$ from Section \ref{notation}, and define a global equivalent
\[E(\ff) \defeq \{x \in \roiplus: x \equiv 1 \newmod{\ff}\} = U(\ff)\cap F^\times.\]
\item[(ii)]We define the \emph{automorphic cycle of level $\ff$} to be
\[X_{\ff} \defeq F^\times\backslash \A_F^\times/U(\ff)F_\infty^1.\]
\end{itemize}
\end{mdef}
\begin{mremnum}\label{autocycledef2}
\begin{longversion}
There is a natural decomposition
\[X_{\ff} = \bigsqcup\limits_{\y \in \clgp{\ff}}X_{\y}.\]
\end{longversion}
\begin{shortversion}
There is a natural decomposition $X_{\ff} = \sqcup_{\y \in \clgp{\ff}}X_{\y},$ where $X_{\y} = \{[x] \in X_{\ff}: x\text{ represents }\y\text{ in }\clgp{\ff}\}.$
\end{shortversion}
\end{mremnum}
There is a natural embedding 
\[\eta_{\ff} : X_{\ff} \longhookrightarrow Y_1(\n)\]
induced by
\begin{align*}\eta : \A_F^\times &\longhookrightarrow \GLt(\A_F)\\
x &\longmapsto \matrd{x}{(x\pi_{\ff}^{-1})_{v|\ff}}{0}{1},
\end{align*}
where $(\pi_{\ff}^{-1})_{v|\ff}$ is the idele defined in Definition \ref{gausssum}. This map is shown to be well-defined in Proposition \ref{welldef} below.
$\lb$
Recall that we have a decomposition $Y_1(\n) = \bigsqcup_{i=1}^h Y_1^i(\n),$ where $Y_1^i(\n)$ is as defined in equation (\ref{modcurvedec}). In particular, $Y_1^i(\n)$ can be described as $\{[g] \in Y_1(\n): \det(g)\text{ represents }i\text{ in }\nclgp\}.$
\begin{mprop}\label{welldef}
The map $\eta_{\ff}$ induces a well-defined map
\[\eta_{\ff}: X_{\ff} \longrightarrow Y_1(\n).\]
Moreover, the restriction of $\eta_{\ff}$ to $X_{\y}$ has image in $Y_1^{i_{\y}}(\n)$, where $i_{\y}$ denotes the element of the narrow class group given by the image of $\y$ under the natural projection $\clgp{\ff} \rightarrow \nclgp$. 
\end{mprop}
\begin{proof}
Suppose $\gamma xur$ is a different representative of $[x]\in X_{\ff}.$ Then
\begin{align}\label{welldefeq} [\eta_{\ff}(\gamma xur)] &= \left[\matrd{\gamma xur}{(\gamma xur\pi_{\ff}^{-1})_{v|\ff}}{0}{1}\right]\notag\\
&= \left[\matrd{\gamma}{0}{0}{1}\matrd{x}{(x\pi_{\ff}^{-1})_{v|\ff}}{0}{1}\matrd{u}{((u-1)\pi_{\ff}^{-1})_{v|\ff}}{0}{1}\matrd{r}{0}{0}{1}\right]\\
&= [\eta_{\ff}(x)] \in Y_1(\n),\notag
\end{align}
showing that the induced map is well-defined. To see that the restriction to $X_{\y}$ lands in $Y_1^{i_{\y}}(\n)$, note that $\det(\eta_{\ff}(x)) = x$, so that if $x$ represents $\y \in \clgp{\ff}$, we see that $\eta_{\ff}(x)$ represents $i_{\y}\in\nclgp$, and in particular, $\eta_{\ff}$ induces a map
\[\{x\in\A_F^\times: [x] = \y \in \clgp{\ff}\} \longrightarrow Y_1^{i_{\y}}(\n),\]
which then descends as claimed.
\end{proof}

\begin{longversion}
\begin{mremnum}\label{multiset}
Later, we will choose $a_{\y}$ as follows. Choose $\{a_i\}$ to be representatives of $\nclgp$ as above, with $(a_i)_{\ff} = (a_i)_\infty = 1$. Now for each $[j] \in (\roi_F/\ff)^\times$, choose $\alpha_j \in \roi_{F,+}^\times$ such that $\alpha_j \equiv j \newmod{\ff}$, and define an idele $b_j$ by
\[(b_j)_v = \left\{\begin{array}{ll}\alpha_j^{-1} &: v\nmid \ff\infty,\\
1 &: v| \ff\infty.\end{array}\right.\]
Then the set 
\[\{a_{ij} \defeq a_ib_j: i \in \nclgp, j \in (\roi_F/\ff)^\times\}\]
is a multiset of representatives of $\clgp{\ff}$, where each class is represented $\#\mathrm{Im}(\roiplus \rightarrow (\roi_F/\ff)^\times)$ times, via the exact sequence
\[\roiplus \longrightarrow (\roi_F/\ff)^\times \longrightarrow \clgp{\ff} \longrightarrow \nclgp \longrightarrow 0.\]
If $a_{ij}$ represents the same element of $\clgp{\ff}$ as $a_{\y}$, then we denote 
\begin{align*}X_{ij} \defeq F^\times\backslash F^\times a_{ij} U(\ff)F_\infty^+/F_\infty^1  = X_{\y}
\end{align*}
and write $\Delta_{ij}$ for the corresponding map. The benefit of this approach is that whilst we do count representatives multiple times, in our later calculations we will be able to write down Gauss sums more effectively, as we can isolate the components coming from $(\roi_F/\ff)^\times$, and means we can use general theory of Gauss sums as devloped by Deligne. The authors apologise for the slightly cumbersome work of carrying around both sets of notation; however, the interplay between them should be apparent, and from now on we will use whichever of the two approaches suits best in particular situations. This will typically be $a_{ij}$ in situations where we develop general theory for individual components (so as to do this in the greatest generality), then using $a_{\y}$ when we want to talk about these objects as a whole (indexed by $\clgp{\ff}$).
\end{mremnum}
\end{longversion}

%
%

\subsection{Evaluation maps}\label{evalumaps}
We now use these automorphic cycles to define \emph{evaluation maps}
\[\mathrm{Ev}:\hc\left(Y_1(\n),\polyC\right) \longrightarrow \C.\]
This will be done in several stages. 

\subsubsection{Pulling back to $X_{\ff}$}
First, we pullback under the inclusion $\eta_{\ff} : X_{\ff} \hookrightarrow Y_1(\n)$. The corresponding sheaf  $\LL_{\ff,1}\polyb \defeq \eta_{\ff}^*\polyC$ can be seen, via equation (\ref{welldefeq}), to be given by the sections of the natural map
\[F^\times\backslash (\A_F^\times \times V_\lambda(\C)^*)/U(\ff)F_\infty^1 \longrightarrow X_{\ff},\]
where the action is given by 
\[f(x,P)ur = \left(fxur,P\bigg|\matrd{f^{-1}}{0}{0}{1}\right).\]

\subsubsection{Passing to individual components}\label{indc}
We can explicitly write
\[X_{\y} \defeq F^\times\backslash F^\times a_{\y} U(\ff)F_\infty^+/U(\ff)F_\infty^1,\]
for $\{a_{\y}: \y \in \clgp{\ff}\}$ a (henceforth fixed) set of class group representatives. Note here that there is an isomorphism
\begin{align}\label{cycleiso}
E(\ff)F_\infty^1\backslash F_\infty^+ &\isorightarrow X_{\y},\\
r &\longmapsto a_{\y}r\notag.
\end{align}
Pulling back under this isomorphism composed with the inclusion $X_{\y} \subset X_{\ff},$ we see that the corresponding sheaf $\LL_{\ff,\y,1} \defeq \tau_{a_{\y}}^*\LL_{\ff,1}\polyb$ is given by the sections of
\[E(\ff)F_\infty^1\backslash (F_\infty^+ \times V_\lambda(\C)^*) \longrightarrow E(\ff)F_\infty^1\backslash F_\infty^+,\]
where now the action is by
\[es(r,P) = \left(esr, P\bigg|\matrd{e^{-1}}{0}{0}{1}\right).\]

\subsubsection{Evaluating}
Let $\jj \in \Z[\Sigma]$ be such that there is a Hecke character $\varphi$ of conductor $\ff$ and infinity type $\jj+\vv$. Note that in this case, for all $e \in E(\ff)$, we have $e^{\jj+\vv} = 1;$ indeed, $e^{\jj+\vv} = \varphi_\infty(e) = \varphi_f(e)^{-1} = 1$, since $e  \equiv 1 \newmod{\ff}.$ Now let $\rho_{\jj}$ denote the map
\[\rho_{\jj} : V_{\lambda}(\C)^* \longrightarrow \C\]
given by evaluating at the polynomial $\XX^{\kk-\jj}\YY^{\jj}$. Then $\rho_{\jj}$ induces a map $(\rho_{\jj})_*$ of local systems on $E(\ff)F_\infty^1\backslash F_\infty^+,$ as
\[P\bigg|\matrd{e^{-1}}{0}{0}{1}(\XX^{\kk-\jj}\YY^{\jj}) = (e^{\jj+\vv})^{-1}P(\XX^{\kk-\jj}\YY^{\jj}) = P(\XX^{\kk-\jj}\YY^{\jj}).\]
We see that the sheaf $(\rho_{\jj})_*\LL_{\ff,\y,1}\polyb$ is the constant sheaf attached to $\C$ over $E(\ff)F_\infty^1\backslash F_\infty^+.$ But note that this space is a connected orientable real manifold of dimension $q$, and hence that there is an isomorphism
\[\hc\left(E(\ff)F_\infty^1\backslash F_\infty^+, \C\right) \cong \C,\]
given by integration over $E(\ff)F_\infty^1\backslash F_\infty^+.$

\begin{mdef}\label{classicalevaluation}
Define 
\[\evyclo : \hc(Y_1(\n),\polyC) \longrightarrow \C\]
to be the composition of the maps 
\begin{align*}\hc(Y_1(\n),\polyC)& \labelrightarrow{\eta_{\ff}^*} \hc(X_{\ff},\LL_{\ff,1}\polyb) \labelrightarrow{\tau_{a_{\y}}^*} \cdots \\
&\hc(\autoc,\LL_{\ff,\y,1}\polyb)\labelrightarrow{(\rho_{\jj})_*} \hc(\autoc,\C) \cong \C.
\end{align*}
\end{mdef}
\begin{mrems}
\begin{itemize}
\item[(i)]
Note that this definition is not restricted to polynomials with coefficients in $\C$. Indeed, the evaluation maps are well-defined for cohomology with coefficients in a number field or an extension of $\Qp$. We will distinguish between the various cases by using a subscript on the cohomology class (for example, $\phi_{\C}$ is a complex modular symbol).
\item[(ii)] The subscript $1$ in $\evyclo$ dictates that this is an evaluation map from the cohomology with coefficients in $\polyC$. Later, we will define an evaluation map $\evyclt$.
\end{itemize}
\end{mrems}

%
%
\begin{shortversion}
\subsection{An integral formula for the $L$-function}
Let $\varphi$ be a Hecke character of conductor $\ff$ and infinity type $\jj+\vv$ for some $0\leq \jj \leq \kk$. The following is a generalisation of a result of Hida:

\begin{mthm}\label{integralformula}
Let $F/\Q$ be a number field, and let $\Phi$ be a cuspidal eigenform over $F$ of weight $\lambda = (\kk,\vv) \in \Z[\Sigma]^2$, where $\kk+2\vv$ is parallel, and let $\varphi$ be a Hecke character of conductor $\ff$ and infinity type $\jj+\vv,$ where $0\leq \jj \leq \kk$. Let $\Lambda(\Phi,\cdot)$ be the normalised $L$-function attached to $\Phi$ defined in Definition \ref{lambda}. Then there is an integral formula
\begin{align*}\sum_{\y\in\clgp{\ff}}\varphi(a_{\y})\evyclo(\phi_{\C}) = (-1)^{R(\jj,\kk)}\left[\frac{|D|\tau(\varphi)}{2^{r_2}}\right]\cdot \Lambda(\Phi,\varphi),
\end{align*}
where:
%
%
\begin{itemize}
\item $\{a_{\y}\}$ is a (fixed) set of adelic representatives for $\clgp{\ff}$ with $(a_{\y})_v = 1$ for $v$ infinite,
\item $R(\jj,\kk) \defeq \sum_{v\in\Sigma(\C)}k_v + \sum_{v\in\Sigma(\R)}k_v + j_v,$
\item $\tau(\varphi)$ is the Gauss sum attached to $\varphi$ defined in Definition \ref{gausssum}, 
\item $D$ is the discriminant of the number field $F$,
\item $\evyclo$ is the classical evaluation map from Definition \ref{classicalevaluation},
\item and $\phi_{\C}$ is the modular symbol attached to $\Phi$ under the Eichler-Shimura isomorphism.
\end{itemize}
\end{mthm}
\begin{proof}(Sketch). The proof is standard but long, messy and technical, and we omit the details. A full and detailed proof can be found in Chapter 12.1.4 of \cite{Wil16}.
$\lb$
The proof relies on explicit computations using the Fourier expansion of the automorphic form. It can be broadly split into several stages, as follows:
\begin{itemize}
\item[(i)] First, we explicitly compute the differential $\delta_{\y}\defeq \tau_{a_{\y}}^*\eta_{\ff}^*\phi_{\C}$. This uses the isomorphism between Betti and de Rham cohomology at the level of complex coefficients, and is done for trivial conductor $\ff$ in \cite{Hid94}, Section 2.5.
\item[(ii)] Write $\delta_{\y} = \sum_{0\leq \jj\leq \kk}\delta_{\y}^{\jj}(z) \mathcal{X}^{\kk-\jj}\mathcal{Y}^{\jj}.$ We then introduce an auxiliary varible $s$ and consider the integral
\[C_{\y}^{\jj}(s) \defeq \int_{E(\ff)F_\infty^1\backslash F_\infty^+} \delta_{\y}^{\jj}(\yy)|\yy|^{s}_{\A_F},\]
where $\yy$ denotes an element of $F_\infty^+$.
\item[(iii)] For $\mathrm{Re}(s)>>0$, we explicitly compute $C_{\y}^{\jj}(s)$, broadly following \cite{Hid94}, Section 7. To do this, we substitute the Fourier expansion of our automorphic form into the expression, and rearrange the result into a product of local integrals at the archimedean places, which are easily computed. We are left with a sum over ideals that are equivalent to $a_{\y}\roi_F$ in $\clgp{\ff}$.
\item[(iv)] By summing over $\y\in\clgp{\ff}$, we get a sum over \emph{all} ideals of $\roi_F$, and this collapses via a Gauss sum to give the value of the $L$-function at $\varphi|\cdot|^s$. We deduce that there is an analytic continuation of $L(\Phi,\varphi,s)$ to the whole complex plane, and that setting $s=0$, we see the (critical) $L$-value at the character $\varphi$.
\item[(v)] We conclude by noting that
\[C_{\y}^{\jj}(0) = \evyclo{\phi_{\C}},\]
from which we deduce the theorem.\qedhere
\end{itemize}
\end{proof}

For later use, it is convenient to record a variant of this theorem here. In particular, in the sequel, we will only be able to consider evaluations at conductors $\ff$ divisible by every prime above $p$. We want to use such evaluations to obtain $L$-values at characters whose conductors do \emph{not} necessarily satisfy this (for example, the trivial character). To do so, we need a compatibility result between evaluation maps for different conductors. By examining the Gauss sum in the proof of the integral formula, we obtain:
\begin{mthm}
Suppose $\phi_{\C}$ is an eigensymbol for all the Hecke operators. Let $\varphi$ be a Hecke character of conductor $\ff$ and infinity type $\jj+\vv$, and let $\pri$ be a prime that divides the level $\n$ but does \emph{not} divide $\ff$. Then
\[\sum_{\mathbf{x}\in\clgp{\ff\pri}}\varphi(a_{\mathbf{x}})\mathrm{Ev}_{\ff\pri,\jj,1}^{a_{\mathbf{x}}}(\phi_{\C}) = (\varphi(\pri)\lambda_{\pri} - 1)\sum_{\y\in\clgp{\ff}}\varphi(a_{\y})\mathrm{Ev}_{\ff,\jj,1}^{a_{\y}}(\phi_{\C}),\]
where $\lambda_{\pri}$ is the Hecke eigenvalue at $\pri$.
\end{mthm}
\begin{mcor}\label{unramified primes}
Suppose $(p)|\n$, and let $\varphi$ be a Hecke character of conductor $\ff|(p^\infty)$ and infinity type $\jj+\vv$. Let $B$ be the set of primes above $p$ for which $\varphi$ is not ramified, and define $\ff' \defeq \ff\prod_{\pri\in B}\pri$. Then $\ff'$ is divisible by every prime above $p$ and we have
\[\sum_{\y\in\clgp{\ff'}}\varphi(a_{\y})\mathrm{Ev}_{\ff',\jj,1}^{a_{\y}}(\phi_{\C}) = \left(\prod_{\pri \in B}(\varphi(\pri)\lambda_{\pri}-1)\right)\sum_{\y\in\clgp{\ff}}\varphi(a_{\y})\mathrm{Ev}_{\ff,\jj,1}^{a_{\y}}(\phi_{\C}).\]
\end{mcor}

\end{shortversion}

\begin{longversion}

%
%
\subsection{An explicit description of $\phi_{\C}$}
We now give an explicit description of the cohomology class $\phi_{\C}$ attached to an automorphic form $\Phi$. We do this by utilising the isomorphism between Betti and de Rham cohomology at the level of complex coefficients (see \cite{Del79}, Section 0.4), which allows us to describe this class as a differential, as in \cite{Hid94}.
$\lb$
Let $\delta_{ij} \defeq \tau_{a_{ij}}^*\eta_{\ff}^*\phi_{\C}.$ Then we can write
\[\delta_{ij}(z) = \sum_{0\leq\jj\leq\kk}\delta_{ij}^{\jj}(z)\mathbf{\mathcal{X}}^{\kk-\jj}\mathbf{\mathcal{Y}}^{\jj}\]
(as elements of the de Rham cohomology), where $\delta^{\jj}_{ij} \in \hc\left(E(\ff)F_\infty^1\backslash F_\infty^+, \C\right).$ Moreover, we see that 
\[\int_{E(\ff)F_\infty^1\backslash F_\infty^+}\delta_{ij}^{\jj} = \evclo(\phi_{\C}).\]

Finally, before giving $\delta_{ij}^{\jj}$ explicitly, we comment on the structure of $E(\ff)F_\infty^1\backslash F_\infty^+.$ We can parametrise this space as the quotient of $\R_{>0}^q$ by units, with one copy of $\R_{>0}$ coming from each real embedding and one from each pair of complex embeddings. This is then isomorphic to $\R_{>0}\times (S^1)^{q-1}.$ The reader should think of this as being an analogue of the path $\{iy \in \uhp: y \in \R_{>0}\}$ (as seen in the rational case when evaluating modular symbols at $\{0\}-\{\infty\}$) in the general setting. 
\begin{mdef}\begin{itemize}
\item[(i)]We parametrise $E(\ff)F_\infty^1\backslash F_\infty^+$ as $E(\ff)\backslash\left\{\y = (y_v)_{v\in \Sigma(\R)\sqcup\Sigma(\C)}: y_v \in \R_{>0}\right\}.$ The use of $\yy$ to mean this rather than a class group representative will be clear from context.
\item[(ii)] If $\jj \in \Z[\Sigma]$, then we define $\jj(\Sigma(\R)) \defeq \sum_{v \in \Sigma(\R)}j_v \in \Z$ (and similarly for $\jj(\Sigma(\C))$). 
\end{itemize}
\end{mdef}
\begin{mrem}
We will use $\yy$ interchangeably to mean an element of $F_\infty^+$ or $F_\infty^1\backslash F_\infty^+$, in the style of Hida. This is for convenience purposes, since the Fourier expansion takes input from the former, whilst the differential has values on the latter. There is, of course, a canonical quotient map between the two, which corresponds to taking norms at complex  places. We also use $\yy$ to mean a representative of $\autoc$, that is, an element of $F_\infty^1\backslash F_\infty^+$ representing a class of this modulo $E(\ff)$. The reader is urged not to get hooked up on the details of this notation!
\end{mrem}
\begin{mprop}\label{c1}
We can explicitly describe $\delta_{ij}^{\jj}$ as follows:
\[\delta_{ij}^{\jj}(\yy) = c_1\yy^{\jj-\kk/2-\ttt}\Phi_{\nn}^u\left(y_\infty^{-1/2}\matrd{a_{ij}\yy}{(\pi_{\ff}^{-1})_{v|\ff}}{0}{1}\right)\bigwedge_{v \in \Sigma(\C)}|y_v|^{-1}d|y_v|\bigwedge_{v \in \Sigma(\R)}dy_v,\]
where:
\begin{itemize}
\item $c_1 \defeq 2^{r_2}N((a_{ij}))^{[\kk+2\vv]/2}\prod_{v\in\Sigma(\C)}(-1)^{k_v+j_{vc}+1}\binomc{k_v^*}{n_v}\prod_{v\in\Sigma(\R)}(-1)^{k_v+j_v}i^{j_v+1},$ for $[\cdot]$ as in Definition \ref{bracket},
\item $\yy \in F_\infty^+$ is considered as an idele by setting $y_v = 1$ for all finite places $v$,
\item $\nn \defeq \sum_{v\in\C}(k_v+j_v-j_{vc}+1)v \in \Z[\Sigma(\C)]$, and $\Phi^u(g) \defeq \sum_{0\leq \rr \leq \kk^*}\Phi_{\rr}^u(g)\XX^{\kk^*-\rr}\YY^{\rr}.$
\end{itemize}
\end{mprop}
\begin{proof}
Most of the work for this is done in \cite{Hid94}, Section 2.5, though with a notable difference stemming from his restriction to characters of trivial conductor. We look at the image of $\y \in \autoc$ under multiplication by $a_{ij}$ composed with the embedding $\eta_{\ff}$. 
$\lb$
Hida's results have no dependence on $j$. To see where this comes in, note that we are restricting $\Phi$ to elements $\smallmatrd{\yy}{\xx}{0}{1}$ where $\xx = (a_{ij}\pi_{\ff}^{-1})_{v|\ff}$. But by definition $(a_{ij})_v = 1$ for all $v|\ff$. Thus $(a_{ij}\pi_{\ff}^{-1})_{v|\ff} = (\pi_{\ff}^{-1})_{v|\ff}$.
$\lb$
The rest of the proof follows exactly as in \cite{Hid94}. Note that he never explicitly uses the dual module $V_\lambda(\C)^*$, instead working implicitly with a particular basis of $V_{\lambda}(\C)$ that corresponds identically to the basis of $V_{\lambda}(\C)^*$ we use under the canonical isomorphism between the two as $\SUt(\C)$-modules (see \cite{Wil17}, Proposition 2.6). We have also passed to the unitisation $\Phi^u(g) \defeq |\det(g)|^{[\kk+2\vv]/2}\Phi(g)$, using the fact that $N((a_{ij})) = |a_{ij}|^{-1}$, as this then has the form we need to use the Fourier expansion developed previously.
\end{proof}

\subsection{An integral formula for the $L$-function}\label{intform}

Let $\varphi$ be a Hecke character of conductor $\ff$ and infinity type $\jj+\vv$ for some $0\leq \jj \leq \kk$. We now look at the image of $\phi_{\C}$ under the evaluation maps, obtaining an integral formula for the $L$-function at $\varphi$. As this calculation is long and messy, for clarity of writing, we have split the work into subsections.

\subsubsection{Notation}
Let $s \in \C$ be an auxiliary variable, and consider the integral
\[C_{ij}^{\jj}(s) \defeq \int_{\autoc}\delta_{ij}^{\jj}(\yy)|\yy|_{\A_F}^{s},\]
where $\y$ denotes an element of $F_\infty^+$ (and can be considered as an idele by setting all finite components equal to 1).
Note here that
\begin{equation}\label{relateCe}
C_{ij}^{\jj}(0) = \evclo(\phi_{\C}).
\end{equation}

\subsubsection{Substituting known expressions}
Now we substitute the explicit value of $\delta_{ij}^{\jj}(\y)$ into $C_{ij}^{\jj}(s)$. This gives
\[C_{ij}^{\jj}(s) = c_1\int_{\autoc}\y^{\jj-\kk/2}\Phi_{n}^u(y_\infty^{-1/2}\matrd{a_{ij}\y}{(\pi_{\ff}^{-1})_{v|\ff}}{0}{1}|\y|_{\A_F}^{s}d^\times|\y|.\]
We can use the Fourier expansion described in Theorem \ref{fourexp}; this yields
\begin{align*}
C_{ij}^{\jj}(s) &= c_1|D|^{1+[\kk+2\vv]/2}\int_{\autoc}\y^{\jj-\kk/2} \sum_{\zeta \in F_+^\times} c(\zeta a_{ij}\Diff,\Phi)e_F(\zeta(\pi_{\ff}^{-1})_{v|\ff})\\
&\times\left[\prod_{v\in\Sigma(\C)} i^{j_v-j_{vc}}\binomc{k_v^*}{n_v}(\zeta y_v)^{j_{vc}-j_v}|\zeta y_v|^{2+k_{vc}+j_v-j_{vc}}K_{j_v-j_{vc}}(4\pi |\zeta y_v|)\right]\\
&\hspace{80pt}\times\left[\prod_{v\in\Sigma(\R)}(|\zeta |y_v)^{1+k_v/2}e^{-2\pi |\zeta |y_v}\right]|\y |_{\A_F}^{s}d^\times|\y|.
\end{align*}
For simplicity, let $\jj^{\#} \defeq (\jj-c\jj)/2$. Grouping together similar terms and rearranging, the above expression simplifies to 
\begin{align*}
C_{ij}^{\jj}(s) &= c_1|D|^{1+[\kk+2\vv]/2}\prod_{v\in\Sigma(\C)}i_{j_v-j_{vc}}\binomc{k_v^*}{n_v}\\
&\times\int_{\autoc}\sum_{\zeta \in F_+^\times} c(\zeta a_{ij}\Diff,\Phi)e_F(\zeta(\pi_{\ff}^{-1})_{v|\ff})\left[\prod_{v\in\Sigma(\C)}|y_v|^{2+j_v+j_{vc}+s}K_{j_v-j_{vc}}(4\pi |\zeta y_v|)\right]\\
&\hspace{80pt}\times\left[\prod_{v\in\Sigma(\R)}(y_v)^{1+j_v+s}e^{-2\pi |\zeta |y_v}\right]\zeta^{\jj^\#}|\zeta|^{1+\kk/2}d^\times|\y|.
\end{align*}

\subsubsection{Rearranging sums and integrals}
For Re$(s) >>0$, we have absolute convergence of the sum; hence we can exchange the order of the sum and integral. Note that in this case, we have
\[\int_{\autoc}\sum_{\zeta\in F_+^{\times}} = \sum_{\zeta\in F_+^\times}\int_{\autoc} = \sum_{\text{ideals }\zeta a_{ij}\Diff}\sum_{\epsilon \in \roi_{F,+}^\times}\int_{\autoc}.\]
Then we have:
\begin{mlem} Let $\epsilon \in E(\ff) \subset \roi_{F,+}^\times$. Then replacing $\zeta$ with $\epsilon \zeta$ in the expression above leaves the integrand unchanged.
\end{mlem}
\begin{proof}
Firstly, it is clear that for any unit, we have $|\epsilon\zeta| = |\zeta|$. There are two other terms involving $\zeta$. One sees from the definition of $e_F$ that if $\epsilon \equiv 1\newmod{\ff}$, then $e_F(\epsilon\zeta(\pi_{\ff}^{-1})_{v|\ff}) = e_F(\zeta(\pi_{\ff}^{-1})_{v|\ff}).$
$\lb$
This just leaves the term $\zeta^{\jj^\#}.$ To deal with this term, recall that we took $\varphi$ to be a Hecke character of conductor $\ff$ and infinity type $\jj+\vv$. Now define 
\[\psi \defeq \varphi|\cdot|^{-[\jj+\vv]}.\]
Then we see that $\psi$ has conductor $\ff$ and infinity type $-\jj^\#$. In particular, we see that 
\[(\epsilon\zeta)^{\jj^\#} = \psi_\infty(\epsilon\zeta)^{-1} = \psi_f(\epsilon\zeta).\]
But $\psi_f(\epsilon) = 1$. Thus it follows that $\zeta^{\jj^\#}$ is invariant under multiplication by $\epsilon$, and we are done.
\end{proof}
In particular, this invariance now allows us to rewrite the integral as
\[[\roi_{F,+}^\times:E(\ff)]\sum_{\text{ideals }\zeta a_{ij}\Diff}\int_{F_\infty^1\backslash F_\infty^+}.\]

\subsubsection{Computing standard integrals}
Using the above, and still assuming that Re$(s) >> 0$, we rearrange further. We can identify $F_\infty^+/F_\infty^1$ with $(\R_{>0})^q$, and hence the integral breaks down into the product of integrals from 0 to $\infty$ at each infinite place. We get:

\begin{align*}
C_{ij}^{\jj}(s) &= c_1[\roi_{F,+}^\times:E(\ff)]|D|^{1+[\kk+2\vv]/2}\prod_{v\in\Sigma(\C)}i_{j_v-j_{vc}}\binomc{k_v^*}{n_v}\\
&\times\sum_{\text{ideals }\zeta a_{ij}\Diff} c(\zeta a_{ij}\Diff,\Phi)e_F(\zeta(\pi_{\ff}^{-1})_{v|\ff})\zeta^{\jj^{\#}}|\jj|^{1+\kk/2}\left[\prod_{v\in\Sigma(\R)}\int_0^\infty(y_v)^{j_v+s}e^{-2\pi |\zeta |y_v}dy_v\right]\\
&\hspace{80pt}\times\left[\prod_{v\in\Sigma(\C)}\int_0^\infty K_{j_v-j_{vc}}(4\pi |\zeta y_v|)|y_v|^{1+j_v+j_{vc}+s}d|y_v|\right].
\end{align*}
These are standard integrals; indeed, we have
%
%
\[\int_0^\infty K_{j_{cv}-j_v}(4\pi|\zeta y_v|)|y_v|^{j_v+j_{vc}+s+1}d|y_v| = (2\pi|\zeta|)^{-j_v-j_{cv}-s-2}2^{-2}\Gamma(j_v+s+1)\Gamma(j_{cv}+s+1),\]
whilst
\[\int_0^\infty e^{-2\pi|\zeta y_v|}|y_v|^{j_v+s}dy_v = (2\pi|\zeta|)^{-j_v-s-1}\int_0^\infty e^{-x}x^{j_v+s}dx = (2\pi|\zeta|)^{-j_v-s-1}\Gamma(j_v+s+1).\]
By substituting these integrals in, we get
%
%
\begin{align*}
C_{ij}^{\jj}(s) &= c_2\sum_{\text{ideals }\zeta a_{ij}\Diff}c(\zeta a_{ij}\Diff)e_F(\zeta(\pi_{\ff}^{-1})_{v|\ff})\zeta^{\jj^{\#}}|\jj|^{1+\kk/2}\\
&\times \left[\prod_{v\in\Sigma(\C)}(2\pi |\zeta|)^{-j_v-j_{cv}-s-2}2^{-2}\Gamma(j_v+s+1)\Gamma(j_{vc} + s + 1)\right]\\
&\hspace{80pt}\times \left[\prod_{v\in\Sigma(\R)}(2\pi |\zeta|^{-j_v-s-1}\Gamma(j_v+s+1)\right],
\end{align*}
which simplifies to
\begin{align*}
C_{ij}^{\jj}(s) &= c_2(2\pi)^{-\jj-(s+1)\ttt}\Gamma(\jj+(s+1)\ttt)\\
&\hspace{80pt}\times\sum_{\text{ideals }\zeta a_{ij}\Diff}c(\zeta a_{ij}\Diff)e_F(\zeta(\pi_{\ff}^{-1})_{v|\ff})\zeta^{\jj^{\#}}|\jj|^{\kk/2 - \jj -s\ttt}.
\end{align*}
Here we have written
\begin{equation}\label{c2}
c_2 = c_1[\roi_{F,+}^\times:E(\ff)]|D|^{1+[\kk+2\vv]/2}\prod_{v\in\Sigma(\C)}i_{j_v-j_{vc}}\binomc{k_v^*}{n_v}\end{equation}
and defined
\[\Gamma(\jj+(s+1)\ttt) \defeq \prod_{v\in\Sigma}\Gamma(j_v+s+1),\]
\[(2\pi)^{\jj+(s+1)\ttt} \defeq \prod_{v\in\Sigma}(2\pi)^{j_v+s+1}\]
for simplicity.

\subsubsection{Simplifying the constant}
We now focus on the term $c_2$. Recall we defined $c_1$ in Proposition \ref{c1}. Subsituting this into equation (\ref{c2}) above, we see that the binomial coefficients cancel, and the signs reduce to give
\begin{align*}
c_2 &= N((a_{ij}))^{[\kk+2\vv]/2}|D|^{1+[\kk+2\vv]/2}[\roi_{F,+}^\times:E(\ff)]2^{r_2}(-1)^{\jj+\ttt}\prod_{v\in\Sigma(\C)}(-1)^{k_v}\prod_{v\in\Sigma(\R)}(-1)^{j_v}\\
&= N(a_{ij}\Diff)^{[\kk+2\vv]/2}|D|[\roi_{F,+}^\times:E(\ff)]2^{r_2}(-1)^{\jj+\ttt}\prod_{v\in\Sigma(\C)}(-1)^{k_v}\prod_{v\in\Sigma(\R)}(-1)^{j_v},
\end{align*}
using the fact that $|D| = N(\Diff)$.

\subsubsection{Rearranging further}
We can massage our formula a bit further; note that 
\begin{align*}|\zeta|^{\kk/2 - \jj - s\ttt} &= N((\zeta))^{[\kk/2 - \jj] - s}\\
&= N(a_{ij}\Diff)^{[\jj-\kk/2]+s}N(\zeta a_{ij}\Diff)^{[\kk/2 - \jj]-s},
\end{align*}
where the first equality follows from the definition of $[\cdot]$ (see Definition \ref{bracket}). When we multiply this by $|D|N(a_{ij}\Diff)^{[\kk+2\vv]/2}$, we obtain $|D|N(a_{ij})^{[\jj+\vv]+s}N(\zeta a_{ij}\Diff)^{\kk/2-\jj]-s}.$ Incorporating all of this, we end up with the formula
\begin{align*}
C_{ij}^{\jj}(s) = c_3 N(a_{ij}\Diff)^{[\jj+\vv]+s}\sum_{\text{ideals }\zeta a_{ij}\Diff}c(\zeta a_{ij}\Diff,\Phi)\zeta^{\jj^{\#}}e_F(\zeta(\pi_{\ff}^{-1})_{v|\ff})N(\zeta a_{ij}\Diff)^{[\kk-\jj/2]-s},
\end{align*}
where
\[c_3 \defeq (-2\pi i)^{-\jj-\ttt}2^{-r_2}(2\pi)^{-qs}|D|\Gamma(\jj+(s+1)\ttt)[\roi_{F,+}^\times:E(\ff)]\prod_{v\in\C}(-1)^{k_v}\prod_{v\in\Sigma(\R)}(-1)^{k_v+j_v}.\]

\subsubsection{Gauss sums}
We now sum over class group representatives and use Gauss sums to obtain the correct twisted $L$-function. Recall that we defined $\varphi$ to be a Hecke character of $F$ of conductor $\ff$ and infinity type $\jj+\vv$, and defined $\psi \defeq \varphi|\cdot|^{-[\jj+\vv]}$, which has infinity type $-\jj^{\#}$. In particular, note that $\zeta^{\jj^{\#}} = \psi_\infty(\zeta)^{-1}$ (since $\zeta$ is totally positive) and we have
\[\varphi(a_{ij})N(a_{ij}\Diff)^{[\jj+\vv]} = N(\Diff)^{[\jj+\vv]}\psi(a_{ij}).\]
Now define
\begin{align*}C_\varphi(s) &\defeq \sum_{i,j} \varphi(a_{ij})N(a_{ij}\Diff)^{-s}C_{ij}^{\jj}(s)\\
&= c_3 \sum_{i,j}\psi(a_{ij})N(\Diff)^{[\jj+\vv]}\sum_{\text{ideals }\zeta a_{ij}\Diff}c(\zeta a_{ij}\Diff,\Phi)\psi_\infty(\zeta)^{-1}\\
&\hspace{80pt}e_F(\zeta(\pi_{\ff}^{-1})_{v|\ff})N(\zeta a_{ij}\Diff)^{[\kk/2-\jj]-s}.
\end{align*}
Here the sum is over $i \in \nclgp$ and $j \in (\roi_F/\ff)^\times$. Ideally, we want the Fourier coefficient to be independent of $j$, so that we can break up the sum and leave a Gauss sum. To achieve this, we scale by $b_j$. Recall $b_j$ is an idele defined to be 1 at places dividing $\ff\infty$ and $\alpha_j^{-1}$ everywhere else, where $\alpha_j \in \roi_{F,+}^\times$ is congruent to $j \newmod{\ff}$. Replace $\zeta$ with $\zeta'\alpha_j$ (noting that we still have absolute convergence). Then $\zeta a_{ij}\Diff = \zeta' a_i\Diff$ (as ideals). This gives
\begin{align*}
C_\varphi(s) = c_3 \sum_{i,j}&\psi(a_{ij})N(\Diff)^{[\jj+\vv]}
\\ &\sum_{\text{ideals }\zeta' a_{i}\Diff}c(\zeta' a_{i}\Diff,\Phi)\psi_\infty(\zeta'\alpha_j)^{-1}e_F(\zeta'\alpha_j(\pi_{\ff}^{-1})_{v|\ff})N(\zeta' a_{i}\Diff)^{[\kk/2-\jj]-s}.
\end{align*}
We need to fix a further piece of notation.
\begin{mnot}Recall: we took $d$ to be a (finite) idele representing the different $\Diff$. We choose a specific $d$. Write $\Diff = I_{i_{\Diff}}(\delta)$, where $\delta \in F_\infty^+$. Then we can taken $d = a_{i_{\Diff}}\delta_f$, where $\delta_f$ is the finite idele with every component equal to $\delta$. It follows that
\begin{equation}\label{expv}
(\zeta'\delta)\alpha_j d^{-1}(\pi_{\ff})_{v|\ff} = \zeta'\alpha_j(\pi_{\ff})_{v|\ff}.
\end{equation}
\end{mnot}
Incorporating equation (\ref{expv}), breaking up $\psi(a_{ij}) = \psi(a_i)\psi(b_j)$, and introducing the term $\psi(d)\psi(d)^{-1}$ (for the Gauss sum), this becomes
\begin{align*}
C_\varphi(s) = c_3 \psi(d)\sum_{i}\psi(a_{i})N(\Diff)^{[\jj+\vv]}
\sum_{\text{ideals }\zeta' a_{i}\Diff}&c(\zeta' a_{i}\Diff,\Phi)\psi_\infty(\zeta')^{-1}N(\zeta' a_{i}\Diff)^{[\kk/2-\jj]-s}\\
&\times\psi(d)^{-1}\sum_{j}\psi(b_j)\psi_\infty(\alpha_j)^{-1}e_F((\zeta'\delta)\alpha_jd^{-1}(\pi_{\ff}^{-1})_{v|\ff}).
\end{align*}
Now, we see that
\begin{align*}\psi(b_j)\psi_\infty(\alpha_j)^{-1} &= \left[\prod_{v\nmid\infty}\psi(\alpha_j^{-1})\prod_{v|\ff}\psi(\alpha_j)\right]\left[\prod_{v\nmid\infty}\psi(\alpha_j)\right]\\
&= \prod_{v|\ff}\psi(\alpha_j) = \psi_{\ff}(\alpha_j),
\end{align*}
so that the second sum becomes
\[\psi(d)^{-1}\sum_{j\in(\roi_F/\ff)^\times} \psi_{\ff}(\alpha_j)e_F((\zeta'\delta)\alpha_jd^{-1}(\pi_{\ff}^{-1})_{v|\ff}) = \left\{\begin{array}{ll}\psi_{\ff}(\zeta'\delta)^{-1}\tau(\psi) &: ((\zeta'\delta),\ff) = 1\\
0 &:\text{otherwise},
\end{array}\right.\]
using the theory of Gauss sums in Section \ref{gsum}. Now note that when $((\zeta'\delta),\ff) = 1$, we have
\[\psi((\zeta'\delta))\psi_{\ff}(\zeta'\delta)\psi_\infty(\zeta'\delta) = 1,\]
by the definition of $\psi_{\ff}$. The sum now becomes
\begin{align*}
C_\varphi(s) = c_3 &\psi(d)\tau(\psi)\sum_{i}\psi(a_{i})N(\Diff)^{[\jj+\vv]}\\
&\times\sum_{\substack{\text{ideals }\zeta' a_{i}\Diff\\ \text{coprime to }\ff}}c(\zeta' a_{i}\Diff,\Phi)\psi((\zeta'\delta))\psi_\infty(\delta)N(\zeta' a_{i}\Diff)^{[\kk/2-\jj]-s}.
\end{align*}
Note that 
\[\psi(d)\psi_\infty(\delta) = \psi(a_{i_{\Diff}})\psi_f(\delta)\psi_\infty(\delta) = \psi(a_{i_{\Diff}}).\]
Rearranging again, and consolidating the terms involving the different and noting that $\psi(\zeta' a_i\Diff) = 0 $ when $((\zeta'\delta),\ff) \neq 1$, this becomes
\begin{align*}
C_\varphi(s) = c_3 \tau(\psi)&N(\Diff)^{[\jj+\vv]}\\
&\times\sum_{i}\sum_{\text{ideals }\zeta' a_{i}\Diff}c(\zeta' a_{i}\Diff,\Phi)\psi(\zeta' a_{i}\Diff)N(\zeta' a_{i}\Diff)^{[\kk/2-\jj]-s}.
\end{align*}
The sum now collapses to one over all ideals of $F$. 

\subsubsection{Obtaining $L$-values}
We have $\psi(d)^{-1}N(\Diff)^{[\jj+\vv]} = \varphi(d)^{-1}$; hence, it is easy to see that $N(\Diff)^{[\jj+\vv]}\tau(\psi) = \tau(\varphi)$. Thus we have
\begin{align*}C_\varphi(s) &= c_3 \tau(\varphi)\mathcal{L}(\Phi,\psi,s-[\kk/2-\jj])\\
&= c_3 \tau(\varphi)\mathcal{L}(\Phi,\varphi,s-[\kk/2 + \vv])\\
&= c_3 \tau(\varphi)L(\Phi,\varphi,s+1).
\end{align*}
With a little extra work, we see that this formula gives an analytic continuation of $L(\Phi,\varphi,s)$ to the complex plane. In particular, setting $s=0$, and recalling that 
\[C_{ij}^{\jj}(0) = \evclo(\phi_{\C}),\]
we see that
\[L(\Phi,\varphi) \defeq L(\Phi,\varphi,1) = \frac{1}{c_3 \tau(\varphi)}\sum_{i,j}\varphi(a_{ij})\evclo(\phi_{\C}).\]
We see that we have proved the following theorem:

\begin{mthm}\label{integralformula}
Let $F/\Q$ be a number field, and let $\Phi$ be a cuspidal eigenform over $F$ of weight $\lambda = (\kk,\vv) \in \Z[\Sigma]^2$, where $\kk+2\vv$ is parallel, and let $\varphi$ be a Hecke character of conductor $\ff$ and infinity type $\jj+\vv,$ where $0\leq \jj \leq \kk$. Let $\Lambda(\Phi,\cdot)$ be the normalised $L$-function attached to $\Phi$ defined in Definition \ref{lambda}. Then there is an integral formula
\begin{align*}\sum_{i,j}\varphi(a_{ij})\evclo(\phi_{\C}) = (-1)^{R(\jj,\kk)}\left[\frac{[\roi_{F,+}^\times:E(\ff)]|D|\tau(\varphi)}{2^{r_2}}\right]\cdot \Lambda(\Phi,\varphi),
\end{align*}
where:
%
%
\begin{itemize}
\item The sum is over $i \in \nclgp$ and $j \in (\roi_F/\ff)^\times$,
\item $R(\jj,\kk) \defeq \sum_{v\in\Sigma(\C)}k_v + \sum_{v\in\Sigma(\R)}k_v + j_v,$
\item $\tau(\varphi)$ is the Gauss sum attached to $\varphi$ defined in Definition \ref{gausssum}, 
\item $D$ is the discriminant of the number field $F$,
\item $\evclo$ is the classical evaluation map from Definition \ref{classicalevaluation},
\item and $\phi_{\C}$ is the modular symbol attached to $\Phi$ under the Eichler-Shimura isomorphism.
\end{itemize}
\end{mthm}

\subsection{Evaluating at ideals other than the conductor}
In the sequel, we will need to look at evaluation maps at ideals other than the conductor of the relevant Hecke character. For example, let $\varphi$ be a Hecke character of conductor $\ff$ and infinity type $\jj+\vv$, and let $\pri$ be a prime not dividing $\ff$; then we will need to consider the expression
\[\sum_{\y\in\clgp{\ff\pri}}\varphi(a_{\y})\mathrm{Ev}_{\ff\pri,\jj,1}^{a_{\y}}(\phi_{\C}).\]
In particular, we need to know how this relates to the evaluation maps at $\ff$ considered above in the case that $\pri$ divides the level $\n$. In this section, we provide a formula for this case.
$\lb$
We start by making the following simple, but crucial, observation about Gauss sums.
\begin{mlem}\label{modified gauss sum}
Let $\varphi$ be a Hecke character of conductor $\ff$, and let $\pri$ be a prime not dividing $\ff$. Let $B$ be a complete set of representatives in $\roi_F$ for the set
\[\{b \newmod{\ff\pri} : b\newmod{\ff} \in (\roi_F/\ff)^\times\}.\]
Then, in the usual set-up for Gauss sums, we have
\[\varphi(d^{-1})\sum_{b\in B}\varphi_{\ff}(b)e_F(\zeta bd^{-1}(\pi_{\ff}^{-1}\pi_{\pri}^{-1})_{v|\ff\pri}) = \left\{\begin{array}{ll}N(\pri)\varphi_{\ff}(\zeta)^{-1}\tau(\varphi) &: ((\zeta),\ff) = 1\text{ and }\pri|(\zeta),\\
0 &: \text{otherwise}.
\end{array}\right.\]
\end{mlem}
\begin{proof}
The sum splits as a product
\[\left[\sum_{\alpha \in \roi_F/\pri}e_F(\zeta\alpha d^{-1}\pi_{\pri}^{-1})\right]\cdot\left[\varphi(d^{-1})\sum_{\beta \in (\roi_F/\ff)^\times}e_F(\zeta\beta d^{-1}(\pi_{\ff}^{-1})_{v|\ff})\right].\]
The first term is non-zero if and only if $\pri|(\zeta)$, in which case $e_F(\zeta\alpha d^{-1}\pi_\pri^{-1}) = 1$ and the sum is $N(\pri)$. The second term is just the usual Gauss sum. The result follows.
\end{proof}

Let $a_{ij} = a_ib_j$, as before, form the usual multiset of representatives for $\clgp{\ff}$. We extend this slightly. To this end, let 
\[\{c_k \in \roi_{F,+}: k\in\roi_F/\pri\}\]
form a complete set of representatives for $\roi_F/\pri$ (noting in particular that we cannot have $c_0 = 0$). Define
\[a_{ijk} = a_{ij}c_k,\]
and note that 
\[\{a_{ijk}: i \in \nclgp, j\in(\roi_F/\ff)^\times, k\in(\roi_F/\pri)^\times\}\]
forms a multiset of representatives for $\clgp{\ff\pri}$, with each representative counted $\#\mathrm{Im}(\roi_{F,+}^\times \longrightarrow (\roi_F/\ff\pri)^\times)$ times. In particular, we see that 
\begin{align*}\sum_{i,j}\sum_{k\in(\roi_K/\pri)^\times}\varphi(a_{ijk})&\mathrm{Ev}_{\ff\pri,\jj,1}^{a_{ijk}}(\phi_{\C}) = \\
&\sum_{i,j}\sum_{k\in\roi_K/\pri}\varphi(a_{ijk})\mathrm{Ev}_{\ff\pri,\jj,1}^{a_{ijk}}(\phi_{\C}) - \sum_{i,j}\varphi(a_{ij}c_0)\mathrm{Ev}_{\ff\pri,\jj,1}^{a_{ij}c_0}(\phi_{\C}),
\end{align*}
where in all three expressions we sum over $i \in \nclgp$ and $j\in(\roi_F/\ff)^\times$. We study each of these terms separately.
\begin{mlem}We have
\[\sum_{i,j}\varphi(a_{ij}c_0)\mathrm{Ev}_{\ff\pri,\jj,1}^{a_{ij}c_0}(\phi_{\C}) =
\frac{[\roi_{F,+}^\times:E(\ff\pri)]}{[\roi_{F,+}^\times:E(\ff)]}\sum_{i,j}\varphi(a_{ij})\mathrm{Ev}_{\ff,\jj,1}^{a_{ij}}(\phi_{\C}),\]
where the sum is over $i \in \nclgp$ and $j\in(\roi_F/\ff)^\times$ in both expressions.
\end{mlem}
\begin{proof}We compute both sides in an almost identical manner to the proof of the integral formula. Firstly, we get the ratio in the unit indices since we are now integrating over $E(\ff\pri)F_\infty^1\backslash F_\infty^+$ rather than $\autoc$. We see that the only other step that changes is the one involving the Gauss sum. To see that this does not affect the final result, note that we have
\[e_F((\zeta'\delta)\alpha_jc_0 d^{-1}(\pi_{\ff\pri}^{-1})_{v|\ff\pri}) = e_F((\zeta'\delta)\alpha_jc_0 d^{-1}(\pi_{\ff}^{-1})_{v|\ff})e_F((\zeta'\delta)\alpha_jc_0 d^{-1}\pi_{\pri}^{-1}),\]
and that the second term of this product is equal to 1 as $c_0 \cong 0 \newmod{\pri}$. The result then follows, since $\{\alpha_jc_0\}$ is again a full set of representatives for $(\roi_F/\ff)^\times$.
\end{proof}
The study of the first term is a little more involved. We give a sketch; the calculations are almost identical to those in the integral formula proved previously, and for the reader's sanity, we do not wish to repeat them.
\begin{mlem}Suppose that $\phi_{\C}$ is an eigenform for all the Hecke operators. We have 
\[\sum_{i,j}\sum_{k\in\roi_K/\pri}\varphi(a_{ijk})\mathrm{Ev}_{\ff\pri,\jj,1}^{a_{ijk}}(\phi_{\C}) = \lambda_{\pri}\varphi(\pri)\frac{[\roi_{F,+}^\times:E(\ff\pri)]}{[\roi_{F,+}^\times:E(\ff)]}\sum_{i,j}\varphi(a_{ij})\mathrm{Ev}_{\ff,\jj,1}^{a_{ij}}(\phi_{\C}),\]
where $\lambda_{\pri}$ is the eigenvalue of $\phi_{\C}$ at the Hecke operator $U_{\pri}$ (recalling that $\pri|\n$).
\end{mlem}
\begin{proof}Again, we examine the proof of the integral formula; the term involving ratios of unit indices is introduced exactly as in the previous lemma. By following the remaining steps in deriving the integral formula, we see again that the only major change is in the Gauss sum, and indeed that we end up with the `modified Gauss sum' of Lemma \ref{modified gauss sum}. In particular, in the calculation of the integral formula, we are left with a sum over ideals that are divisible by $\pri$, in addition to introducing a factor of $N(\pri)$. Since $\phi_{\C}$ is an eigenform, and as $\pri|\n$, the Fourier coefficients satisfy
\[c(I\pri,\Phi) = c(I,\Phi)c(\pri,\Phi),\]
so that the summands are multiplicative and we can recover a sum over all ideals by factoring out the expression $c(\pri,\Phi)\varphi(\pri)N(\pri)^{[\kk/2+\vv]-s}$. After setting $s=0$, we incorporate the extra factor of $N(\pri)$ coming from the Gauss sum and recover the result by noting that 
\[\lambda_{\pri} = c(\pri,\Phi)N(\pri)^{[\kk/2+\vv]+1}\]
from previously (see also \cite{Hid94}, equation (6.2b)).
\end{proof}
We now conclude this section by stating the compatibility results we need. We find we have proved the following:

\begin{mthm}Let $\varphi$ be a Hecke character of conductor $\ff$ and infinity type $\jj+\vv$, and let $\pri$ be a prime dividing $\n$ but \emph{not} dividing $\ff$. Then we have
\[\sum_{i,j,k}\varphi(a_{ijk})\mathrm{Ev}_{\ff\pri,\jj,1}^{a_{ijk}}(\phi_{\C}) = (\varphi(\pri)\lambda_{\pri}-1)\frac{[\roi_{F,+}^\times:E(\ff\pri)]}{[\roi_{F,+}^\times:E(\ff)]}\sum_{i,j}\varphi(a_{ij})\mathrm{Ev}_{\ff,\jj,1}^{a_{ij}}(\phi_{\C}),\]
where the sums are over $i\in\nclgp$, $j\in(\roi_F/\ff)^\times$ and $k\in(\roi_F/\pri)^\times$.
\end{mthm}
In the next section, we will remove the terms of form $[\roi_{F,+}^\times:E(\ff\pri)]$ by summing only over the respective narrow ray class groups (and, in the process, eliminating the double counting of class group representatives; see Theorem \ref{evphithm}). With this in mind, we record the following corollary, which is proved by a simple induction:
\begin{mcor}\label{unramified primes}
Suppose $(p)|\n$, and let $\varphi$ be a Hecke character of conductor $\ff|(p^\infty)$ and infinity type $\jj+\vv$. Let $B$ be the set of primes above $p$ for which $\varphi$ is not ramified, and define $\ff' \defeq \ff\prod_{\pri\in B}\pri$. Then $\ff'$ is divisible by every prime above $p$ and we have
\[\sum_{\y\in\clgp{\ff'}}\varphi(a_{\y})\mathrm{Ev}_{\ff',\jj,1}^{a_{\y}}(\phi_{\C}) = \left(\prod_{\pri \in B}(\varphi(\pri)\lambda_{\pri}-1)\right)\sum_{\y\in\clgp{\ff}}\varphi(a_{\y})\mathrm{Ev}_{\ff,\jj,1}^{a_{\y}}(\phi_{\C}).\]
\end{mcor}

\end{longversion}

%
%

\section{Algebraicity results}\label{algebraic}
So far, all of our work has been done over $\C$. We will now refine these results to connect the \emph{algebraic} modular symbol to the critical $L$-values above.
\begin{longversion}
$\lb$
Recall Theorem \ref{deligneperiods}, which said that the normalised $L$-value $\Lambda(\Phi,\varphi)$ is an algebraic multiple of a period $\Omega_{\Phi}^{\epsphi}$, where $\epsphi$ is the character of $\pmr$ attached to $\varphi$ (see Section \ref{heckechar}). In Section \ref{attach}, for a character $\varepsilon$ of $\pmr$ we defined a modular symbol $\phi_{\C}^\varepsilon$ and stated a result that $\phi_K^\varepsilon \defeq \phi_{\C}^\varepsilon/\Omega_\Phi^\varepsilon$ lived in an algebraic subspace. We now relate $\phi_{\C}^\varepsilon$ (and, by scaling, $\phi_K^\varepsilon$) to the $L$-function using our above formula.
\end{longversion}
\begin{mdef}Let $A_{\ff} = \{a_{\y}: \y \in \clgp{\ff}\}$ denote a fixed set of representatives for $\clgp{\ff}$, with components at infinity that are not necessarily trivial. For a Hecke character $\varphi$ of conductor $\ff$ and infinity type $\jj+\vv$, where $0 \leq \jj \leq \kk$, define a function
\[\ev_\varphi^{A_{\ff}} : \hc\left(Y_1(\n),\polyC\right) \longrightarrow \C\]
by
\[\ev_\varphi^{A_{\ff}}(\phi) = \sum_{\y \in\clgp{\ff}}\epsphi\varphi_{f}(a_{\y})\evyclo(\phi),\]
where as previously we write $\epsphi$ as a function on the ideles by composing it with the natural sign map $\A_F^\times \rightarrow \pmr$.
\end{mdef}
\begin{longversion}
This definition is intimately related to the locally analytic function $\varphi_{p-\mathrm{fin}}$ we defined in Section \ref{rcg}. In particular, note that $\epsphi\varphi_f = \varphi/\varphi_\infty^{\mathrm{alg}}$, where $\varphi_\infty^{\mathrm{alg}}(x) = x_\infty^{\jj+\vv}$ is the unique algebraic function on $F_\infty$ that agrees with $\varphi_\infty$ on $F_\infty^+$.
\end{longversion}
\begin{mlem}The function $\ev_{\varphi}^{A_{\ff}}$ is independent of class group representatives.
\end{mlem}
\begin{proof}
Let $a_{\y}'$ be an alternative representative corresponding to $\y \in \clgp{\ff}$. Then $a_{\yy} = fa_{\y}ur$, where $f \in F^\times, u \in U(\ff)$ and $r \in F_\infty^+$. Looking at the description of the evaluation maps, we see that 
\[\mathrm{Ev}_{\ff,\jj,1}^{a_{\y}'}(\phi) = f^{\jj+\vv}\evyclo(\phi).\]
But 
\begin{align*}\epsphi\varphi_f(a_{\yy}') = \epsphi\varphi_f(fa_{\y}ur) &= \epsphi\varphi_f(f)\epsphi\varphi_f(a_{\y}) \\
&= f^{-\jj-\vv}\epsphi\varphi_f(a_{\y}),
\end{align*}
since $\epsphi\varphi_f$ is trivial on $U(\ff)F_\infty^+$ and by our earlier comment, we have
\[\epsphi\varphi_f(f) = \varphi(f)/\varphi_\infty^{\mathrm{alg}}(f) = f^{-\jj-\vv}.\]
Putting this together, we find that
\[\epsphi\varphi_f(a_{\y}')\ev_{\ff,\jj,1}^{a_{\y}'}(\phi) = \epsphi\varphi_f(a_{\y})\evyclo(\phi),\]
which is the required result.
\end{proof}
\begin{mdef}\label{evphi}Define $\ev_\varphi$ to be the map $\ev_\varphi^{A_{\ff}}$ for any choice of class group representatives $A_{\ff}$. This is well-defined by the above lemma.
\end{mdef}
We will combine this with the following to deduce the result we desire.
\begin{mprop}
Let $\iota \in \pmr$. Then for any idele $a$, we have 
\[\mathrm{Ev}_{\ff,\jj,1}^{\iota a}\bigg(\phi\bigg|\iota\bigg) = \mathrm{Ev}_{\ff,\jj,1}^{a}(\phi).\]
\end{mprop}
\begin{proof}
Recall that the definition of the action of $\iota \in \pmr$ on the cohomology of $Y_1(\n)$ was described in Section \ref{weyl}. There is a well-defined action of $\pmr$ on the local system corresponding to $\LL_{\ff,1}\polyb$ given by
\[\iota\cdot(x,P) = (\iota x,P),\]
where here we have considered $\iota$ to be an idele by setting $\iota_v = 1$ for all complex and finite places $v$. A simple check shows that if $\phi \in \hc(Y_1(\n),\polyC)$ then we have
\[\eta_{\ff}^*(\phi|\iota) = \eta^*_{\ff}(\phi)|\iota\]
coming from the commutative diagram
\[
\CommDia{(g,P)}{\eta_{\ff}^*}{(x,P)}{|\iota}{|\iota}{(\iota\cdot g,P)}{\eta_{\ff}^*}{(\iota x,P)}
\]
of local systems. Continuing to work at the level of local systems, suppose $x$ is an idele that, under the natural quotient map, lies in the component of $X_{\ff}$ corresponding to $a_{\y}$. Then the image of $\iota x$ lies in the component corresponding to $\iota a_{\y}$ (where here we note that if $\{a_{\y}:\y \in \clgp{\ff}\}$ is a complete set of representatives for $\clgp{\ff}$, then so is the set $\{\iota a_{\y}:\y \in \clgp{\ff}\}$). Thus we see that there is a commutative diagram of maps of local systems 
\begin{diagram}
 (x,P)&&\rTo^{\tau_{a_{\y}}^*}&&(r,P)&&\rTo^{\text{ev. at }\XX^{\kk-\jj}\YY^{\jj}}&&(r,c) \\
\dTo_{|\iota}&&&&\dTo{}&&&&\dTo{}\\
\\
\left(\iota x,P\right)&&\rTo^{\tau_{\iota a_{\y}}^*}&&\left(r,P\right)&&\rTo^{\text{ev. at }\XX^{\kk-\jj}\YY^{\jj}}&&(r,c),
\end{diagram}
where the local system on the far right hand side defines the constant sheaf given by sections of $(E(\ff)F_\infty^1\backslash F_\infty^+) \times \C \rightarrow E(\ff)F_\infty^1\backslash F_\infty^+$. The result follows.
\end{proof}
\begin{mcor}We have the relation
\[\mathrm{Ev}_\varphi(\phi|\iota) = \epsphi(\iota)\mathrm{Ev}_\varphi(\phi).\]
\end{mcor}\label{evcor}
\begin{proof}
Considering $\iota$ as an idele in the usual way, we have
\begin{align*}
\ev_{\varphi}(\phi|\iota) &= \sum_{\y \in \clgp{\ff}}\epsphi\varphi_f(\iota a_{\y})\ev_{\jj,\ff,1}^{\iota a_{\y}}(\phi|\iota)\\
&= \epsphi(\iota)\sum_{\y\in\clgp{\ff}}\epsphi\varphi_f(a_{\y})\evyclo(\phi)\\
&=\epsphi(\iota)\ev_{\varphi}(\phi),
\end{align*}
as required.
\end{proof}
\begin{mcor}
We have
\[\mathrm{Ev}_\varphi(\phi_{\C}^\varepsilon) = \left\{\begin{array}{ll}\mathrm{Ev}_\varphi(\phi_{\C}) &: \varepsilon = \epsphi\\
0 &: otherwise.
\end{array}\right.\]
\end{mcor}
\begin{proof}
By definition,
\begin{align*}\ev_\varphi(\phi_{\C}^\varepsilon) &= \ev_\varphi\left(2^{-r_1}\sum_{\iota \in \pmr}\varepsilon(\iota)\phi_{\C}|\iota\right)\\
&= \left[2^{-r_1}\sum_{\iota \in \pmr}\varepsilon(\iota)\epsphi(\iota)\right]\ev_\varphi(\phi_{\C}),
\end{align*}
using linearity of the evaluation maps and Corollary \ref{evcor}. The result then follows from orthogonality of characters, since $\epsphi^2 = 1$.
\end{proof}
Recall that in Definition \ref{thetak}, we set $\theta_K \defeq \sum_{\varepsilon}\phi_K^\varepsilon$. Note that here $\theta_K$ is an element of the cohomology with \emph{algebraic} coefficients in the number field $K$.
\begin{mthm}\label{evphithm}
Let $\varphi$ be a Hecke character of conductor $\ff$ and infinity type $\jj + \vv$, where $0 \leq \jj \leq \kk$, and write $\epsphi$ for the associated character of $\pmr$ defined in Section \ref{heckechar}. Let $\ev_\varphi$ be as in Definition \ref{evphi}. We have
\[\mathrm{Ev}_\varphi(\theta_K) = (-1)^{R(\jj,\kk)}\left[\frac{|D|\tau(\varphi)}{2^{r_2}\Omega_\Phi^{\epsphi}}\right]\cdot \Lambda(\Phi,\varphi),\]
where $R(\jj,\kk) = \sum_{v\in\Sigma(\R)}j_v+k_v + \sum_{v\in\Sigma(\C)}k_v$.
\end{mthm}
\begin{longversion}
\begin{proof}We use Theorem \ref{integralformula}. In particular, note that in the statement of the theorem, we use the multiset $\{a_{ij}\}$ of class group representatives in which each element of the class group is represented $[\roi_{F,+}^\times:E(\ff)] = \#\mathrm{Im}(\roi_{F,+}^\times \rightarrow (\roi_F/\ff)^\times)$ times, so we can cancel this term from the result. Then we see that for these representatives,
\[\epsphi\varphi_f(a_{ij}) = \varphi(a_{ij}),\]
since we chose $(a_{ij})_\infty = 1$, so that the sum we obtained in the statement of Theorem \ref{integralformula} is nothing but $[\roi_{F,+}^\times:E(\ff)]\ev_{\varphi}(\phi_{\C})$. The result follows.
\end{proof}
\end{longversion}
\begin{shortversion}
\begin{proof}
We use Theorem \ref{integralformula}. In particular, note that we choose $(a_{\y})_\infty = 1$, so that
\[\epsphi\varphi_f(a_{\y}) = \varphi(a_{\y}).\]
Thus the sum we obtained in the statement of this theorem is exactly $\ev_{\varphi}(\phi_{\C})$. The result follows.
\end{proof}
\end{shortversion}
To summarise: we have now defined an \emph{algebraic} cohomology class that sees the algebraic parts of all of the critical $L$-values that we hope to interpolate. In particular, by embedding $K$ into a sufficiently large finite extension $L/\Qp$, we get a $p$-adic modular symbol $\theta_L$ that sees all of these critical values.

%
%

\section{Distributions and overconvergent cohomology}\label{distributions}
In this section, we define the distribution modules that we will use as coefficient modules for the spaces of overconvergent modular symbols. This closely follows the analogous section of \cite{Bar13}.
$\lb$
Throughout this section, $L$ is a finite extension of $\Qp$ containing the image of $\inc\circ\sigma: F \hookrightarrow \overline{\Q}_p$ for each embedding $\sigma$ of $F$ into $\overline{\Q}$. First, we give some motivation by reformulating the definition of the space $V_\lambda(L)$. We previously defined this to be the $d$-fold tensor product of the polynomial spaces $V_{k_v}(L)$, with an action of $\GLt(L)$ depending on $\lambda$. Note that $\roifp$ embeds naturally in $\overline{\Q}_p^d$, and in particular, we can see an element of $V_\lambda(L)$ as a function on $\roifp$ in a natural way. We see that the following definition agrees with the definition we gave in Section \ref{setup}.
\begin{mdef}
Let $L/\Qp$ be a finite extension and let $\lambda = (\kk,\vv) \in \Z[\Sigma]$ be admissible (so that, in particular, $\kk\geq 0$). Define $V_\lambda(L)$ to be the space of functions on $\roifp$ that are polynomial of degree at most $\kk$ with coefficients in $L$, with a left action of $\GLt(\roifp)$ given by
\[\matr\cdot P(x) = (ad-bc)^{\vv}(a+cx)^{\kk}P\left(\frac{b+dx}{a+cx}\right).\]
\end{mdef}
We have passed to a non-homogeneous version here. This definition is more easily seen to be compatible with the rest of this section. In particular, it is compatible with the following:
\begin{mdef}
Let $\mathcal{A}(L)$ be the space of locally analytic functions on $\roifp$ that are defined over $L$.
\end{mdef}
We would like to define an action of $\GLt(\roifp)$ on this space, analogously to above. Unfortunately, the action above does not extend to the full space $\mathcal{A}(L)$. We can, however, define an action of a different semigroup.
\begin{mdef}
\begin{itemize}
\item[(i)] Let $\Sigma_0(p)$ be the semigroup
\[\Sigma_0(p) \defeq \left\{\matr \in M_2(\roifp): c \in p\roifp, a\in (\roifp)^\times, ad-bc \neq 0\right\}.\]
\item[(ii)] Define $\mathcal{A}_\lambda(L)$ to be the space $\mathcal{A}(L)$ equipped with a left `weight $\lambda$ action' of $\Sigma_0(p)$ given by
\[\matr\cdot f(z) = (ad-bc)^{\vv}(a+cz)^{\kk}P\left(\frac{b+dz}{a+cz}\right).\]
\end{itemize}
\end{mdef}
Note in particular that this semigroup contains the image of $\Gamma_1(\n)$ under the natural embedding $M_2(\roi_F)\subset M_2(\roifp)$ as well as the matrices that we will need to define a Hecke action at primes above $p$. It is \emph{not} a subset of $\GLt(\roifp)$, but the action of this different semigroup also extends naturally to $V_\lambda(L)$, since both live inside $\GLt(F\otimes_{\Q}\Qp)$.
$\lb$
We are now in a position to define the distribution spaces. 
\begin{mdef}
Define $\dist \defeq \mathrm{Hom}_{\mathrm{cts}}(\mathcal{A}_\lambda(L),L)$ to be the topological dual of $\mathcal{A}_\lambda$, with a right action of $\Sigma_0(p)$ defined by 
\[(\mu|\gamma)(f) \defeq \mu(\gamma\cdot f).\]
\end{mdef}
Note that $\Omega_1(\n)$ acts on $\dist$ via its projection to $\GLt(\Qp)$, giving rise to a local system $\sht{\dist}$ on $Y_1(\n)$.
\begin{mdef}
The space of \emph{overconvergent modular symbols} is the compactly supported cohomology group $\hc(Y_1(\n),\sht{\dist})$.
\end{mdef}
By dualising the inclusion $V_\lambda(L)\subset \mathcal{A}_\lambda(L)$, we get a $\Sigma_0(p)$-equivariant surjection
\[\mathcal{D}_\lambda(L) \longrightarrow V_\lambda(L)^*.\]
This gives rise to a $\Sigma_0(p)$-equivariant \emph{specialisation map}, a map
\[\rho: \hc(Y_1(\n),\sht{\dist}) \longrightarrow \hc(Y_1(\n),\sht{V_\lambda(L)^*}).\]
The space of overconvergent modular symbols is, in a sense, a $p$-adic deformation of the space of classical modular symbols. It was introduced by Glenn Stevens in \cite{Ste94}.
$\lb$
We conclude this section with a result that will be crucial in the following section, where we prove that the space of overconvergent modular symbols admits a slope decomposition with respect to the Hecke operators. For the relevant definitions, see \cite{Urb11}, Section 2.3.12. The space $\dist$ is naturally a nuclear Fr\'{e}chet space\footnote{That is, an inverse limit of Banach spaces in which the projection maps are compact. In \cite{Urb11}, Urban calls this a \emph{compact Fr\'{e}chet space}. We instead follow the terminology utilised in \cite{Sch02}.}; indeed, let $\mathcal{A}_{n,\lambda}(L)$ be the space of functions that are \emph{locally analytic of order $n$}, that is, functions that are analytic on each open set of the form $a+p^n\roifp$. Each $\mathcal{A}_{n,\lambda}(L)$ is a Banach space, and the inclusions $\mathcal{A}_{n,\lambda}(L)\hookrightarrow \mathcal{A}_{n+1,\lambda}(L)$ are compact (\cite{Urb11}, Lemma 3.2.2). We write $\mathcal{D}_{n,\lambda}(L)$ for the topological dual of $\mathcal{A}_{n,\lambda}(L)$. Then $\dist \cong \varprojlim\mathcal{D}_{n,\lambda}(L)$ is equipped with a family of norms coming from the Banach spaces $\mathcal{D}_{n,\lambda}(L)$.
\begin{mdef}Let $M \cong \varprojlim M_n$ be a nuclear Fr\'{e}chet space. We say that an endomorphism $U$ of $M$ is \emph{compact} if it is continuous and there are continuous maps $U_n'$ making the following commute:
\[\CommDia{M}{}{M_{n-1}}{U}{U_{n}'}{M}{}{M_{n}},\]
where the horizontal maps are the natural projections. 
\end{mdef}
\begin{longversion}
In this situation, we obtain compact\footnote{Urban uses `compact' and `completely continuous' interchangeably to describe endomorphisms of Banach spaces that map bounded subsets into relatively compact subsets.} endomorphisms $U_n$ on $M_n$ by composing $U_n'$ with the natural map $M_{n}\rightarrow M_{n-1}$. In \cite{Ser62}, it is proved that if $M_n$ is a Banach space equipped with a compact endomorphism $U_n$, then $M_n$ admits a slope decomposition with respect to $U_n$, and in \cite{Urb11}, Section 2.3.10, Urban uses this -- and compactness of $U$ -- to deduce the existence of a slope decomposition for $M$ with respect to $U$. In particular, the following lemma will be crucial in the next section.
\end{longversion}
\begin{mlem}
Let $\eta \in \GLt(F)\cap\Sigma_0(p)$, which acts naturally on $\dist$. This action is compact. In particular, the action of $\smallmatrd{1}{0}{0}{p}$ is compact on $\dist$.
\end{mlem}
\begin{proof}
See \cite{Urb11}, Lemma 3.2.8.
\end{proof}

%
%

\section{Slope decompositions}\label{slopedecomp}
We start by recalling the relevant definitions about slope decompositions.
\begin{mdef}
Let $L$ be a finite extension of $\Q_{p}$, and let $h \in \Q$. We say a polynomial $Q(X) \in L[X]$ has \emph{slope $\leq h$} if $Q(0)\in \roi_L^\times$ and if $\alpha\in \overline{L}$ is a root of $Q^\ast(X)\defeq X^{\deg(Q)}Q(1/X)$, then $v_p(\alpha) \leq h$.
\end{mdef}
\begin{mdef}\label{slope decomposition}
Let $M$ be an $L$-vector space equipped with the action of an $L$-linear endomorphism $U$. We say that $M$ has a \emph{slope $\leq h$ decomposition with respect to $U$} if there is a decomposition $M\cong M_1\oplus M_2$ such that:
\begin{itemize}
\item[(i)] $M_1$ is finite-dimensional,
\item[(ii)] The polynomial $\det(1-UX)|_{M_1}$ has slope $\leq h$, and
\item[(iii)] For all polynomials $P\in L[X]$ with slope $\leq h$, the polynomial $P^\ast(U)$ acts invertibly on $M_2$.
\end{itemize}
We write $M^{\leq h,U} \defeq M_1$ for the elements of \emph{slope $\leq h$} in $M$. Where the operator $U$ is clear, we drop it from the notation and just write $M^{\leq h}$.
\end{mdef}

\begin{longversion}
In this section, we will prove the following theorem:
\end{longversion}
\begin{shortversion}
The crucial theorem we require is the following:
\end{shortversion}
\begin{mthm} \label{existence of slope decomposition} Let $\lambda= (\mathbf{k}, \mathbf{v})$ be an admissible weight. Then for each $i \in \N$ and any $h \in \Q$, the $L$-vector space $\mathrm{H}_{\mathrm{c}}^{i}(Y_{1}(\mathfrak{n}), \mathcal{L}_{2}(\mathcal{D}_{\lambda}(L)))$ admits a slope $\leq h$ decomposition with respect to the Hecke operator $U_{p}$.
\end{mthm}

\begin{shortversion}
\begin{proof}\emph{(Sketch).} To prove this theorem we follow the arguments given in \cite{Urb11} and \cite{Bar15}, where the same statement is proved in the cases of the cohomology without compact support and $\mathrm{GL}_{2}$ over a totally real field respectively. Both of these rely on general results from earlier in \cite{Urb11}, where Urban proves that any nuclear Fr\'{e}chet space $M$ equipped with a compact endomorphism $U$ admits a slope decomposition with respect to $U$. Given this, the key step is to construct a complex whose cohomology is $\mathrm{H}_{\mathrm{c}}^{\bullet}(Y_{1}(\mathfrak{n}), \mathcal{L}_{2}(\mathcal{D}_{\lambda}(L)))$ and such that each term of the complex is isomorphic to finitely many copies of $\mathcal{D}_{\lambda}(L)$. We can find a lift of the Hecke operators on the cohomology to this complex, and then we use the fact that the action of $(\begin{smallmatrix} 1 & 0 \\ 0 & p \end{smallmatrix})$ on $\mathcal{D}_{\lambda}(L)$ is compact to deduce that this lift acts compactly on the complex. Using Urban's results, we deduce the theorem.
\end{proof}
\end{shortversion}

\begin{longversion}
To prove this theorem we follow the arguments given in \cite{Urb11} and \cite{Bar15}, where the same statement is proved in the cases of the cohomology without compact support and $\mathrm{GL}_{2}$ over a totally real field respectively. Both of these rely on general results from earlier in \cite{Urb11}, where Urban proves that any nuclear Fr\'{e}chet space $M$ equipped with a compact endomorphism $U$ admits a slope decomposition with respect to $U$. Given this, the key step is to construct a complex whose cohomology is $\mathrm{H}_{\mathrm{c}}^{\bullet}(Y_{1}(\mathfrak{n}), \mathcal{L}_{2}(\mathcal{D}_{\lambda}(L)))$ and such that each term of the complex is isomorphic to finitely many copies of $\mathcal{D}_{\lambda}(L)$. We can find a lift of the Hecke operators on the cohomology to this complex, and then we use the fact that the action of $(\begin{smallmatrix} 1 & 0 \\ 0 & p \end{smallmatrix})$ on $\mathcal{D}_{\lambda}(L)$ is compact to deduce that this lift acts compactly on the complex. Using Urban's results, we deduce the theorem.

\subsection{Compactly supported cohomology}

\subsubsection{Complexes}

Let $\Gamma$ be a torsion-free arithmetic subgroup of $\mathrm{SL}_{2}(F)$ and consider the manifold $\Gamma \setminus \mathcal{H}_{F}$. We denote by $\Gamma \setminus \overline{\mathcal{H}}_{F}$ the Borel-Serre compactification of $\Gamma\setminus \mathcal{H}_{F}$ (See \cite{BS74}). Let $\mathfrak{B}$ be the set of proper parabolic $\Q$-subgroups of $\mathrm{SL}_{2}(F\otimes \R)$. To construct the Borel-Serre compactification, recall that we first enlarge $\mathcal{H}_{F}$ to a space $\overline{\mathcal{H}}_{F}$ by adding a euclidean space $e(P)$ of dimension $d$ to each $P \in \mathfrak{B}$. The boundary of $\overline{\mathcal{H}}_{F}$ is given by
\begin{equation}\label{boundary decomposition}
\partial\overline{\mathcal{H}}_{F}= \bigsqcup_{P \in \mathfrak{B}}e(P). 
\end{equation}

The group $\Gamma$ acts on $\overline{\mathcal{H}}_{F}$. The quotient $\Gamma\setminus \overline{\mathcal{H}}_{F}$ is a compact, smooth and $C^{\infty}$-variety with boundary, so there exists a finite triangulation (see \cite{Mun67}) that induces a triangulation on the boundary $\partial (\Gamma\setminus \overline{\mathcal{H}}_{F})$. From this, we obtain a triangulation of $\overline{\mathcal{H}}_{F}$ which contains a triangulation of $\partial(\overline{\mathcal{H}}_{F})$. We consider the complexes of simplicial chains attached to those triangulations, denoted by
$$C_{\bullet}(\Gamma) \ \ \mathrm{and} \ \    C_{\bullet}^{\partial}(\Gamma).$$

These complexes satisfy the following properties:
\begin{itemize}
\item If $i \in \N$, then $C_{i}(\Gamma)$ and $C_{i}^{\partial}(\Gamma)$ are $\Z[\Gamma]$-free modules of finite rank. Since the group $\Gamma$ acts properly on $\overline{\mathcal{H}}_{F}$ and its boundary, this is a consequence of the fact that $\Gamma$ is torsion-free and the (fixed) triangulation of $\Gamma\setminus \overline{\mathcal{H}}_{F}$ is finite.

\item The complex $C_{\bullet}(\Gamma)$ is a resolution of the trivial $\Z[\Gamma]$-module $\Z$, since this complex gives the homology of $\overline{\mathcal{H}}_{F}$, which is contractible.
 
\item The complex $C_{\bullet}^{\partial}(\Gamma)$ is a resolution of the $\Z[\Gamma]$-module $\Z^{\mathfrak{B}}$, where the action of $\Gamma$ on $\Z^{\mathfrak{B}}$ is given by the action on $\mathfrak{B}$. In fact, $C_{\bullet}^{\partial}(\Gamma)$ gives the homology of $\partial\overline{\mathcal{H}}_{F}$, and in decomposition (\ref{boundary decomposition}), each $e(P)$ is contractible.

\end{itemize} 
 
Suppose now that the image of $\Gamma$ in $G(\Q_{p})$ is contained in the Iwahori group (that is, the matrices that are upper-triangular modulo $p$), as is the case for the groups $\Gamma_1^i(\mathfrak{n})$. Then any right $\Omega_1(\n)$-module $M$, as in Definition \ref{locsys} (ii), has an action of $\Gamma$ (the reader should keep the case $M = \dist$ in mind). We define the complexes $C^{\bullet}(\Gamma, M)$ and $C^{\bullet}_{\partial}(\Gamma, M)$ by:
$$C^{i}(\Gamma, M):= \mathrm{Hom}_{\Z[\Gamma]}(C_{i}(\Gamma), M),$$ 
$$C^{i}_{\partial}(\Gamma, M):= \mathrm{Hom}_{\Z[\Gamma]}(C_{i}^{\partial}(\Gamma), M).$$ 
The properties given above for $C_{\bullet}(\Gamma)$ and $C_{\bullet}^{\partial}(\Gamma)$ have the following consequences:
\begin{itemize}
\item[•] $C^{i}(\Gamma, M)$ and $C^{i}_{\partial}(\Gamma, M)$ are isomorphic to finitely many copies of $M$. In particular they are nuclear Fr\'{e}chet $L$-vector spaces.
\item[•] The cohomology of $C^{\bullet}(\Gamma, M)$ is isomorphic to $\mathrm{H}^{\bullet}(\Gamma \setminus \mathcal{H}_{F}, \mathcal{L}_2(M))$.
\item[•]  The cohomology of $C^{\bullet}_{\partial}(\Gamma, M)$ is isomorphic to $\mathrm{H}^{\bullet}(\partial(\Gamma \setminus \mathcal{H}_{F}), \mathcal{L}_2(M))$.

\end{itemize}

Via the natural map $C_{\bullet}^{\partial}(\Gamma) \hookrightarrow C_{\bullet}(\Gamma)$, we obtain a map of complexes $C^{\bullet}(\Gamma, M) \rightarrow C^{\bullet}_{\partial}(\Gamma, M)$. We define:
$$C^{\bullet}_{c}(\Gamma, M)\defeq \mathrm{Cone}\left[C^{\bullet}(\Gamma, M) \rightarrow C^{\bullet}_{\partial}(\Gamma, M)\right].$$

\begin{mprop} \label{complex with compact supports} For each $i \in \N$, the $L$-vector space $C^{i}_{c}(\Gamma, M)$ is a Fr\'{e}chet space. The cohomology of the complex $C^{\bullet}_{c}(\Gamma, M)$ is $\mathrm{H}^{\bullet}_{\mathrm{c}}(\Gamma \setminus \mathcal{H}_{F}, \mathcal{L}(M))$.
\end{mprop}
\begin{proof} By construction we have $C^{i}_{c}(\Gamma, M)= C^{i}(\Gamma, M)\oplus C^{i-1}_{\partial}(\Gamma, M)$. Thus this is a Fr\'{e}chet space as $C^{i}(\Gamma, M)$  and $C^{i-1}_{\partial}(\Gamma, M)$ both are.
$\lb$
For the remainder of the proposition, note that there are isomorphisms 
\[\mathrm{H}^{i}(C^{\bullet}(\Gamma, M))\cong \mathrm{H}^{i}(\Gamma \setminus \mathcal{H}_{F}, \mathcal{L}_2(M)) \hspace{12pt} \mathrm{and}\]
 \[\mathrm{H}^{i}(C^{\bullet}_{\partial}(\Gamma, M))\cong \mathrm{H}^{i}(\partial(\Gamma \setminus \mathcal{H}_{F}), \mathcal{L}_2(M)).\]
 Moreover we have two long exact sequences
 $$... \rightarrow \mathrm{H}^{i}(C^{\bullet}_{c}(\Gamma, M)) \rightarrow \mathrm{H}^{i}(C^{\bullet}(\Gamma, M)) \rightarrow \mathrm{H}^{i}(C^{\bullet}_{\partial}(\Gamma, M)) \rightarrow ...  \ \ \mathrm{and}$$ 
$$... \rightarrow \mathrm{H}^{i}_{c}(\Gamma \setminus \mathcal{H}_{F}, \mathcal{L}(M)) \rightarrow \mathrm{H}^{i}(\Gamma \setminus \mathcal{H}_{F}, \mathcal{L}(M)) \rightarrow \mathrm{H}^{i}( \partial(\Gamma\setminus\mathcal{H}_{F}), \mathcal{L}(M)) \rightarrow ... \hspace{6pt}.$$
Applying the five-lemma to this gives the result.
\end{proof}

\subsubsection{Hecke operators}

Let $\Gamma$ and $\Gamma'$ be as before, let $h: \Gamma\rightarrow \Gamma'$ be a group homomorphism and let $f: M \rightarrow M$ be a linear transformation such that $f(h(\gamma)\mu)= \gamma f(\mu)$. Using $h$ we can consider the complex $C_{\bullet}(\Gamma')$ as a resolution of $\Z$ by $\Z[\Gamma]$-modules. Since $C_{\bullet}(\Gamma)$ is a projective resolution of $\Z$ by $\Z[\Gamma]$-modules, we obtain a map $h_{\bullet}: C_{\bullet}(\Gamma) \rightarrow C_{\bullet}(\Gamma')$ compatible with $h$. Using this last morphism and $f$ we obtain a map:
\begin{align*}
C^{\bullet}(\Gamma', M) &\longrightarrow C^{\bullet}(\Gamma, M),\\
\varphi &\longmapsto  f \circ \varphi \circ h_{\bullet}.
\end{align*}
In the same way as before, we consider the complex $C_{\bullet}^{\partial}(\Gamma')$ as a resolution of $\Z^{\mathfrak{B}}$ by $\Z[\Gamma]$-modules, giving a map of complexes $C^{\partial}_{\bullet}(\Gamma) \rightarrow C_{i}^{\bullet}(\Gamma')$. Then we obtain maps
\[C^{\bullet}_{\partial}(\Gamma', M) \longrightarrow C^{\bullet}_{\partial}(\Gamma, M).\]

From these maps we obtain and the natural compatibility with the maps $C^{\bullet}(\Gamma, M) \rightarrow C^{\bullet}_{\partial}(\Gamma, M)$ and $C^{\bullet}(\Gamma', M) \rightarrow C^{\bullet}_{\partial}(\Gamma', M)$, we obtain maps
$$C^{\bullet}_{c}(\Gamma', M) \longrightarrow C^{\bullet}_{c}(\Gamma, M').$$

Let $\eta \in \GLt(F)\cap\Sigma_0(p)$ (again, the reader should keep in mind $\eta = \smallmatrd{1}{0}{0}{p}$), and define maps
\begin{align*}
h_{1}: \Gamma\cap \eta\Gamma \eta^{-1} &\longrightarrow \Gamma\cap \eta^{-1}\Gamma \eta,\\
\gamma &\longmapsto \eta^{-1}\gamma\eta
\end{align*} and $f_1: M\rightarrow M, f_{1}(\mu)= \mu|\eta$. Let $h_{2}: \Gamma\cap \eta^{-1}\Gamma \eta \hookrightarrow \Gamma$ be the inclusion and $f_{2}$ the identity map. By considering the pairs $(h_{1}, f_{1})$ and $(h_{2}, f_{2})$ in the situation above, we obtain maps
$$[\eta]: C^{\bullet}_{c}(\Gamma\cap \eta^{-1}\Gamma \eta, M) \rightarrow C^{\bullet}_{c}(\Gamma\cap \eta\Gamma \eta^{-1}, M),$$
$$\mathrm{res}^{\Gamma}_{\Gamma\cap \eta^{-1}\Gamma \eta}: C^{\bullet}_{c}(\Gamma, M) \rightarrow C^{\bullet}_{c}(\Gamma\cap \eta^{-1}\Gamma \eta, M).$$
We define corestriction maps $\mathrm{cor}^{\Gamma}_{\Gamma\cap \eta\Gamma \eta^{-1}}$ for the complexes $C^{\bullet}$ and $C^{\bullet}_{\partial}$  in the same way as in \cite[\S 4.2.5]{Urb11} and \cite[\S 2.2.2]{Bar15}. Then we obtain maps
$$\mathrm{cor}^{\Gamma}_{\Gamma\cap \eta\Gamma \eta^{-1}}: C^{i}_{c}(\Gamma\cap \eta\Gamma\eta^{-1}, M) \rightarrow C^{i}_{c}(\Gamma, M).$$

Denote by $[\Gamma\eta\Gamma]= \mathrm{cor}^{\Gamma}_{\Gamma\cap \eta\Gamma \eta^{-1}} \circ [\eta] \circ \mathrm{res}^{\Gamma}_{\Gamma\cap \eta^{-1}\Gamma \eta}$ the composition
$$[\Gamma\eta\Gamma]: C^{i}_{c}(\Gamma, M) \rightarrow C^{i}_{c}(\Gamma, M).$$

\begin{mprop} \label{operator is compact on the complex} The operator $[\Gamma\eta\Gamma]$ is compact. Moreover, if $\eta = \smallmatrd{1}{0}{0}{p},$ it is a lift of the $U_p$ operator on the cohomology to the level of complexes.
\end{mprop}
\begin{proof}
The action of $\eta$ on $\mathcal{D}_{\lambda}(L)$ is compact. Using this property, Proposition \ref{complex with compact supports}, and the fact that composition of a compact map with a continuous map is again compact, we deduce that $[\Gamma\eta\Gamma]$ is compact on $C^{i}_{c}(\Gamma, M)$ (see \cite{Urb11}, Section 4.2.9).
$\lb$
Now fix $\eta = \smallmatrd{1}{0}{0}{p}$. The operators $[\Gamma\eta\Gamma]$ on $C^{\bullet}(\Gamma, M)$ and on $C^{i}_{\bullet}(\Gamma, M)$ lift the corresponding $U_p$ operators on $\mathrm{H}^{\bullet}(\Gamma \setminus \mathcal{H}_{F}, \mathcal{L}_2(M))$ and $\mathrm{H}^{\bullet}_{\partial}(\Gamma \setminus \mathcal{H}_{F}, \mathcal{L}_2(M))$. Hence we deduce the same result for the compact support situation.\end{proof}

\subsection{Proof of Theorem \ref{existence of slope decomposition}} \begin{proof} (Theorem \ref{existence of slope decomposition}). 
Recall that we have the decomposition
\begin{equation}\label{decomposition cohomology}
\mathrm{H}_{\mathrm{c}}^{i}(Y_{1}(\mathfrak{n}), \mathcal{L}_{2}(\mathcal{D}_{\lambda}(L)))= \bigoplus_{j= 1}^{h}\mathrm{H}_{\mathrm{c}}^{i}(Y_{1}^{j}(\mathfrak{n}), \mathcal{L}_2(\mathcal{D}_{\lambda}(L))).
\end{equation} 
Moreover, we can describe the action of the $U_p$ operator on $\mathrm{H}_{\mathrm{c}}^{i}(Y_{1}(\mathfrak{n}), \mathcal{L}_{2}(\mathcal{D}_{\lambda}(L)))$ with respect to this decomposition; indeed, we have
\begin{equation}\label{decomposition of hecke operator}
U_{p}= \bigoplus_{j= 1}^{h}\left[\Gamma_{1}^{j}(\mathfrak{n})\matrd{1}{0}{0}{p}\Gamma_{1}^{j}(\mathfrak{n})\right],
\end{equation}
where $\Gamma_1^j(\n)$ is as defined in equation (\ref{gammai}). For each $i \in \N$, let 
\[C^{i}_{c}(\mathfrak{n}, \mathcal{D}_{\lambda}(L))\defeq \bigoplus_{j= 1}^{h}C^{i}_{c}(\Gamma_{1}^{j}(\mathfrak{n}), \mathcal{D}_{\lambda}(L)).\]
Using Proposition \ref{complex with compact supports} and equation (\ref{decomposition cohomology}), we deduce that each term $C^{i}_{c}(\mathfrak{n}, \mathcal{D}_{\lambda}(L))$ is a Fréchet $L$-vector space and that the cohomology of the complex $C^{\bullet}_{c}(\mathfrak{n}, \mathcal{D}_{\lambda}(L))$ is $\mathrm{H}_{\mathrm{c}}^{\bullet}(Y_{1}(\mathfrak{n}), \mathcal{L}_{2}(\mathcal{D}_{\lambda}(L)))$.
$\lb$
For each $j$ we have an operator
$$\left[\Gamma_{1}^{j}(\mathfrak{n})\matrd{1}{0}{0}{p}\Gamma_{1}^{j}(\mathfrak{n})\right]: C^{i}_{c}\left(\Gamma_{1}^{j}(\mathfrak{n}), \mathcal{D}_{\lambda}(L)\right) \rightarrow C^{i}_{c}\left(\Gamma_{1}^{j}(\mathfrak{n}), \mathcal{D}_{\lambda}(L)\right),$$
 lifting the corresponding operator on the cohomology. We define the operator $$U_{p}: C^{i}_{c}(\mathfrak{n}, \mathcal{D}_{\lambda}(L)) \rightarrow C^{i}_{c}(\mathfrak{n}, \mathcal{D}_{\lambda}(L))$$ by $U_{p}= \bigoplus_{j= 1}^{h}[\Gamma_{1}^{j}(\mathfrak{n})(\begin{smallmatrix} 1 & 0 \\ 0 & p \end{smallmatrix})\Gamma_{1}^{j}(\mathfrak{n})]$. The decomposition of equation (\ref{decomposition of hecke operator}) implies that $U_{p}$ is a lift of the corresponding Hecke operator on the cohomology to the level of complexes. Moreover, from Proposition \ref{operator is compact on the complex} we deduce that $U_{p}$ is a compact operator on $C^{i}_{c}(\mathfrak{n}, \mathcal{D}_{\lambda}(L))$.
 $\lb$
Finally, we complete the proof of Theorem \ref{existence of slope decomposition} by applying \cite[Lemma 2.3.13]{Urb11} to $U_{p}$.
\end{proof}

\end{longversion}

%
%

\section{A control theorem}\label{controltheorem}
In this section, we prove a \emph{control theorem}, showing that the restriction of the specialisation map from overconvergent to classical modular symbols to the `small slope' subspaces is an isomorphism. We actually need a slightly finer definition of slope decomposition; namely, we define the slope decomposition with respect to a finite set of operators rather than just one.
$\lb$
To this end, let $I$ be a finite set, and suppose that for each $i\in I$, we have an endomorphism $U_i$ on the $L$-vector space $M$. Write $A\defeq L[U_{i}, i \in I]$ for the algebra of polynomials in the variables $U_{i}$. Then $A$ acts on $M$, and for $\mathbf{h}= (h_{i}) \in \Q^{I}$ we define the \emph{slope $\leq \mathbf{h}$ subspace with respect to $A$} to be
\[M^{\leq \mathbf{h},A}\defeq \bigcap_{i \in I}M^{\leq h_i, U_{i}}.\]
Where the choice of operators is clear, we will drop the $A$ from the notation and just write $M^{\leq \mathbf{h}}.$
\subsection{Preliminary results}
We start by stating some properties of slope decompositions that will be required in the proof.
\begin{mlem}  \label{lemma finite}  \begin{itemize}
\item[(i)] Let $M,N$ and $P$ be $L$-vector spaces equipped with an action of $A$, and suppose that $M,N$ and $P$ each admit a slope $\leq \mathbf{h}$ decomposition with respect to $A$. If $ 0 \rightarrow M \rightarrow N \rightarrow P \rightarrow 0$ is an exact sequence of $A$-modules, then we have an exact sequence
\[ 0 \rightarrow M^{\leq \mathbf{h}} \rightarrow N^{\leq \mathbf{h}} \rightarrow P^{\leq \mathbf{h}} \rightarrow 0.\]
\item[(ii)] Let $M\cong\varprojlim M_n$ be a nuclear Fr\'{e}chet space equipped with a compact endomorphism $U$ that induces compact operators $U_n$ on $M_n$ for each $n$. Then for each $n$ there is an isomorphism
\[M^{\leq \mathbf{h},U} \cong M_n^{\leq \mathbf{h},U_n}.\]
This fact holds as well for compact maps between complexes of nuclear Fr\'{e}chet spaces and the induced slope decomposition on their cohomology.
\item[(iii)] Let $(M,||\cdot||)$ be an $L$-Banach space equipped with an action of $A$, where $||\cdot||$ denotes the norm on $M$, and suppose that there is a $\roi_L$-submodule 
\[\mathcal{M} \subset \{m\in M: ||m|| \geq 0\}\]
that is stable under the action of $A$. Let $\mathbf{h}= (h_{i})_{i \in I}$ with $h_{i_{0}}< 0$ for some $i_{0} \in I$. Then $M^{\leq \mathbf{h}}= {0}$. 
\end{itemize}
\end{mlem}
\begin{proof} Part (i) is simple (see Corollary 2.3.5 of \cite{Urb11}). Part (ii) is proved in \cite{Urb11}, Lemma 2.3.13. For part (iii), suppose that $M^{\leq \mathbf{h}}\neq {0}$. Then, after possibly replacing $L$ with a finite extension, we can find $\alpha \in L$ and $x \in \mathcal{M}$ such that $v_{p}(\alpha) < 0$ and $U_{i_{0}}x= \alpha x$. Then there exists $n \in \Z$ such that $\alpha^{n} x \notin \mathcal{M}$. This is a contradiction because $\alpha^{n} x= U_{i_{0}}^{n}x \in \mathcal{M}$ by $A$-stability of $\mathcal{M}$.
\end{proof}
In particular, we have the following corollary.
\begin{mdef}
For each $\sigma \in \Sigma$, denote by $\pri(\sigma)$ the unique prime $\pri|p$ such that the embedding $\sigma : F\hookrightarrow \overline{\Q} \subset \C$ extends to an embedding $F_{\pri}\hookrightarrow \overline{\Qp} \subset \Cp$ that is compatible with the fixed embedding $\inc: \overline{\Q}\hookrightarrow \overline{\Qp}$. If $\sigma$ corresponds to $\pri$ under this identification, we write $\sigma \sim \pri$.
\end{mdef}
\begin{mdef}\label{vpri}Let $\nu = (\kk,\vv) \in \Z[\Sigma]^2$ be an admissible weight. Define
\[v_{\pri}(\nu) \defeq \sum_{\sigma \sim \pri}v_\sigma.\]
\end{mdef}
\begin{mcor}\label{vanishing}\begin{itemize}
\item[(i)] Let $\nu = (\kk,\vv) \in \Z[\Sigma]^2$ be a weight with $\kk+2\vv$ parallel (but allowing for negative values of $k_\sigma)$. Let $\mathbf{h} \in \Q^{\{\pri|p\}}$ be such that 
\[h_{\pri} < \frac{v_{\pri}(\nu)}{e_{\pri}}\]
for some prime $\pri$ above $p$. Then for all $r$ we have $\h_{\mathrm{c}}^r(Y_1(\n), \sht{\mathcal{D}_{\nu}(L)})^{\leq \mathbf{h}} = \{0\}$.
\item[(ii)] Under the same hypotheses, the same result holds if we replace $\mathcal{D}_{\nu}(L)$ with any $\Sigma_0(p)$-stable submodule or by quotients by such submodules. 
\end{itemize}
\end{mcor}
\begin{proof}
From Section \ref{distributions}, we know that 
\[\mathcal{D}_\lambda(L) \cong \varprojlim \mathcal{D}_{\lambda,n}(L),\]
where $\mathcal{D}_{\lambda,n}(L)$ is the ($L$-Banach space) of distributions that are locally analytic of order $n$. We also know (from results in the previous section) that the cohomology group $\h_{\mathrm{c}}^r(Y_1(\n),\sht{\mathcal{D}_{\nu,0}(L)})$ is an $L$-Banach space, and we see that $\h_{\mathrm{c}}^r(Y_1,\sht{\mathcal{D}_{\nu,0}(\roi_L)})$ is a $\roi_L$-submodule of the elements of non-negative norm. This space is \emph{not} necessarily preserved by the Hecke operators at $p$, but it \emph{is} preserved by the modified operators
\[U_{\pri}' \defeq \pi_{\pri}^{-v_{\pri}(\nu)}U_{\pri},\]
where we scale by $\pi_{\pri}^{-v_{\pri}(\nu)}$ to ensure integrality in the case $v_{\pri}(\nu)$ is large and negative. Write $A' \defeq L[U_{\pri}']$ for the algebra generated by these modified operators. Applying parts (ii) and (iii) of the above lemma, we see that if $\mathbf{h}' \in \Q^{\{\pri|p\}}$ is chosen such that $h_{\pri}' < 0$ for some prime $\pri$ above $p$, we have
\[\h_{\mathrm{c}}^r(Y_1(\n),\sht{\mathcal{D}_{\nu,0}(L)})^{\leq \mathbf{h}',A'} = \{0\}.\]
By part (ii) of the above lemma, the finite slope cohomologies of $\mathcal{D}_\nu(L)$ and $\mathcal{D}_{\nu,0}(L)$ are isomorphic; hence we conclude that
\[\h_{\mathrm{c}}^r(Y_1(\n),\sht{\mathcal{D}_{\nu}(L)})^{\leq \mathbf{h}',A'} = \{0\}.\] Now note that for any operator $U$ on a nuclear Fr\'{e}chet space $M$, we have a relation
\[M^{\leq h,p^kU} \cong M^{\leq h-k,U}.\]
In particular, define $\mathbf{h} \in \Q^{\{\pri|p\}}$ by
\[h_{\pri} \defeq h_{\pri}' + \frac{v_{\pri}(\nu)}{e_{\pri}}.\]
Note that $h_{\pri}'<0$ for some $\pri$ above $p$ if and only if $h_{\pri} < \frac{v_\pri(\nu)}{e_{\pri}}$ for some $\pri$ above $p$, and that the space on which the Hecke operators at $p$ act with slope $\leq \mathbf{h}$ is isomorphic to the space on which the operators $U_{\pri}'$ act with slope $\leq \mathbf{h}'$. Part (i) follows.
$\lb$
The proof for submodules is identical. The case of quotients then follows by taking a long exact sequence, applying Lemma \ref{lemma finite}(i), and using the result for submodules.
\end{proof}
%
%
%
%
%

\subsection{Theta maps and partially overconvergent coefficients}
We now introduce modules of partially overconvergent coefficients that will play a key role in the proof. 
$\lb$
For any $\sigma \in \Sigma$, let $\lambda_{\sigma}= (\mathbf{k}',\mathbf{v}')$ be the weight defined by 
\begin{equation}\label{lambdasigma}
k'_{\tau}= \left\{\begin{array}{ll} k_\tau &: \tau\neq \sigma,\\
-2-k_{\sigma} &: \tau = \sigma.\end{array}\right.,
\hspace{20pt} v'_{\sigma}= \left\{\begin{array}{ll} v_\tau &: \tau\neq \sigma,\\
v_{\sigma} + k_{\sigma} + 1 &: \tau = \sigma.\end{array}\right.
\end{equation}
Let $f$ be a locally analytic function on $\roifp$, and let $\{V\}$ be an open cover of $\roifp$ such that $f|_V$ is analytic for each $V$. Then we can consider $f|_V$ as a power series in the $d$ variables $\{z_\sigma:\sigma\in\Sigma\}$. We can consider the operator $(d/dz_{\sigma})^{k_\sigma + 1}$ on such power series in the natural way, and note that this induces a map
\[\Theta_{\sigma}: \mathcal{A}_\lambda(L) \longrightarrow \mathcal{A}_{\lambda_{\sigma}}(L).\]
For more details about this map, see \cite[Prop. 3.2.11]{Urb11}. Taking the continuous dual of this map, we obtain a map 
\[\Theta_{\sigma}^{\ast}: \mathcal{D}_{\lambda_{\sigma}}(L) \longrightarrow \mathcal{D}_{\lambda}(L).\]
\begin{mrem}
This map is equivariant with respect to the action of $\Sigma_0(p)$. Note, however, that the action of the $U_{\pri}$ operator is \emph{different} on $\mathcal{D}_{\lambda_{\sigma}}(L)$ and $\mathcal{D}_\lambda(L)$, due to the scaling of $\vv$ at $\sigma$. Indeed, we introduce a factor of the determinant of the component at $\sigma$ to the power of $k_\sigma + 1$.
\end{mrem}
Now label the elements of $\Sigma$ as $\sigma_{1},\sigma_{2},..., \sigma_{d}$, where we can choose any ordering of the elements. We write $\Theta_{0}^{\ast}: \lbrace 0\rbrace \rightarrow \mathcal{D}_{\lambda}$, and for each $s = 1,...,d$, we denote by $\Theta_{s}^{\ast}$ the map
\[\Theta_s^{\ast} \defeq \sum_{i= 1}^{s}\Theta_{\sigma_{i}}^{\ast}: \bigoplus_{i= 1}^{s}\mathcal{D}_{\lambda_{\sigma_{i}}}(L) \longrightarrow \mathcal{D}_{\lambda}(L).\]
The cokernels of the maps $\Theta_s^*$ play a crucial role in the sequel. In particular, from the definition it is clear that $\coker(\Theta_0^*) = \dist$. Consider now the map $\Theta_1^*$. If $\mu \in \mathcal{D}_{\lambda_{\sigma_1}}(L),$ then $\Theta_1^*(\mu)$ is $0$ on elements of $\A_\lambda(L)$ that are locally polynomial in $z_{\sigma_1}$ of degree at most $k_{\sigma_1}$. Hence, for $\mu \in \dist$, we have $\mu \notin \mathrm{Im}(\Theta_1^*)$ if and only if there exists a monomial $\mathbf{z}^{\mathbf{r}} \defeq \prod_{\sigma\in\Sigma}z_{\sigma}^{r_{\sigma}}$ with $r_{\sigma_1} \leq k_{\sigma_1}+1$ such that $\mu(\mathbf{z}^{\mathbf{r}}) \neq 0$. From this one can see that $\coker(\Theta_1^*)$ can be seen as the module of coefficients that are classical at $\sigma_1$ and overconvergent at $\sigma_2, ..., \sigma_d.$ This motivates the following:
\begin{mdef}
Let $J \subset \Sigma$. For $\nu  = (\kk,\vv)\in \Z[\Sigma]^2$, define $\mathcal{A}^J_\nu(L)$ to be the space of functions on $\roifp$ defined over $L$ that are locally analytic in the variables $z_\sigma$ for $\sigma \notin J$ and locally algebraic of degree at most $\max(k_\sigma,0)$ in the variables $z_\sigma$ for $\sigma\in J$. Define $\mathcal{D}^J_\nu(L)$ be the topological dual of $\mathcal{A}^J_\nu(L)$.
\end{mdef}
Thus we see that $\coker(\Theta_1^*) = \mathcal{D}^{\{\sigma_1\}}_\lambda(L)$. Continuing in the same vein, we see that $\coker(\Theta_s^*) = \mathcal{D}^{J_s}_\lambda(L)$, where $J_s \defeq \{\sigma_1,...,\sigma_s\}.$ In particular, if we write $V_{\lambda,\mathrm{loc}}(L)$ for the space of locally algebraic polynomials on $\roifp$ of degree at most $\kk$, with the natural action of $\Sigma_0(p)$ depending on $\lambda$, then we get:
\begin{mprop}\label{bgg}
There is an exact sequence
\[\bigoplus_{\sigma \in \Sigma}\mathcal{D}_{\lambda_{\sigma}}(L) \labelrightarrow{\Theta_d^*} \dist \longrightarrow V_{\lambda,\mathrm{loc}}(L) \longrightarrow 0.\]
In particular, we have
\[\coker(\Theta_d^*) = \mathcal{D}^\Sigma_\lambda(L) \cong V_{\lambda,\mathrm{loc}}(L)^*.\]
\end{mprop}
These are the last terms of the locally analytic BGG resolution introduced in \cite{Urb11}, Section 3.3. See Proposition 3.2.12 of Urban's paper for further details of this exact sequence.

\begin{longversion}
\begin{mrem}
This bears comparison to the results in \cite{Wil17}, Section 6, where an analagous theorem is proved in the case that $F$ is imaginary quadratic and $p$ splits as $\pri\pribar$ in $F$. In particular, it is shown that in this setting, we can lift small slope classical symbols to small slope `half-overconvergent' symbols using the $U_{\pri}$-operator, and then to small slope fully overconvergent symbols using the $U_{\pribar}$-operator. In this way, we get independent notions of small slope at each of the primes above $p$, allowing us to extend the notion of `non-critical' slope and prove a control theorem for a larger range of classical eigensymbols. 
\end{mrem}
\end{longversion}

\subsection{The control theorem}
The following theorem is the main result of this part of the paper, and allows us to canonically lift small-slope classical modular symbols to overconvergent modular symbols. 
\begin{mthm}\label{comparaison} 
Let $\lambda = (\kk,\vv)$ be an admissible weight, and let $\mathbf{h}= (h_{\pri})_{\pri | p} \in \Q^{\lbrace\mathfrak{p}\mid p\rbrace}$. Let $k_{\pri}^0\defeq \min\{k_\sigma: \sigma\sim\pri\}$ and recall the definition of $v_{\pri}(\lambda)$ from Definition \ref{vpri}. If for each prime $\pri$ above $p$ we have 
\begin{equation}\label{smallslopecond}
h_{\pri} < \frac{k_{\pri}^0+v_{\pri}(\lambda)+ 1}{e_{\mathfrak{p}}},\end{equation}
then, for each $r$, the restriction
\[
 \xymatrix{
\rho: \hct^{r}(Y_{1}(\mathfrak{n}), \mathcal{L}_{2}(\mathcal{D}_{\lambda}(L)))^{\leq \mathbf{h}} \ar^{\sim}[r]&
\hct^{r}(Y_{1}(\mathfrak{n}), \mathcal{L}_{2}(V_{\lambda}(L)^{\ast}))^{\leq \mathbf{h}}
}
\]
of the specialisation map to the slope $\leq \mathbf{h}$ subspaces with respect to the $U_{\pri}$-operators is an isomorphism.
\end{mthm}

To prove this, we make use of:

\begin{mlem}In the set-up of Theorem \ref{comparaison}, if $\mathbf{h}$ satisfies equation (\ref{smallslopecond}), then for any $s$ there is an isomorphism
\[
 \xymatrix{
\hct^{r}(Y_{1}(\mathfrak{n}), \mathcal{L}_{2}(\mathcal{D}^{J_{s-1}}_\lambda(L)))^{\leq \mathbf{h}} \ar^{\sim}[r]&
\hct^{r}(Y_{1}(\mathfrak{n}), \mathcal{L}_{2}(\mathcal{D}^{J_s}_\lambda(L)))^{\leq \mathbf{h}}
}
\]
induced from the natural specialisation maps.
\end{mlem}

\begin{proof} We follow \cite{Urb11}. For any $\sigma \in \Sigma$, let $\lambda_{\sigma}= (\mathbf{k}',\mathbf{v}')$ be the weight defined in equation (\ref{lambdasigma}), and recall the theta maps
\[\Theta_{s}^{\ast}: \bigoplus_{i = 1}^s \mathcal{D}_{\lambda_{\sigma_i}}(L) \rightarrow \mathcal{D}_{\lambda}(L).\]
Recall that $\coker(\Theta_s^*) = \mathcal{D}^{J_s}_\lambda(L)$ can be viewed as a module of distributions that are classical at $\sigma_1,...,\sigma_s$ and overconvergent at $\sigma_{s+1},...,\sigma_d$. In particular, there are natural projection maps $\mathcal{D}^{J_{s-1}}_\lambda(L)\rightarrow \mathcal{D}^{J_{s}}_\lambda(L)$ given by specialising from overconvergent to classical coefficients at $\sigma_{s}$. Moreover, from the definition of $\Theta_{\sigma_s}^*$ there is an exact sequence
\[\mathcal{D}_{\lambda_s}(L) \labelrightarrow{\Theta_{\sigma_s}^*} \mathcal{D}^{J_{s-1}}_\lambda(L) \longrightarrow \mathcal{D}^{J_{s}}_\lambda(L) \longrightarrow 0,\] 
and a closer inspection shows that the sequence
 \begin{equation}\label{exact quotient}
   0 \longrightarrow  \mathcal{D}_{\lambda_{\sigma_s}}^{J_{s-1}}(L) \longrightarrow \mathcal{D}^{J_{s-1}}_\lambda(L)\longrightarrow \mathcal{D}^{J_{s}}_\lambda(L) \longrightarrow 0
\end{equation}
is exact for the quotient $\mathcal{D}_{\lambda_{\sigma_s}}^{J_{s-1}}(L)$ of $\mathcal{D}_{\lambda_{\sigma_s}}(L)$.
$\lb$
Using Lemma \ref{lemma finite} on the exact sequence of equation (\ref{exact quotient}), we obtain the exact sequence
\begin{align*}\cdots \longrightarrow \hct^{i}(Y_{1}(\mathfrak{n}), &\mathcal{L}_{2}(\mathcal{D}_{\lambda_{\sigma_s}}^{J_{s-1}}(L)))^{\leq \mathbf{h}} \longrightarrow \hct^{i}(Y_{1}(\mathfrak{n}), \mathcal{L}_{2}(\mathcal{D}^{J_{s-1}}_\lambda(L)))^{\leq \mathbf{h}}\\
&\longrightarrow \hct^{i}(Y_{1}(\mathfrak{n}), \mathcal{L}_{2}(\mathcal{D}^{J_{s}}_\lambda(L)))^{\leq \mathbf{h}} \longrightarrow \hct^{i+ 1}(Y_{1}(\mathfrak{n}), \mathcal{L}_{2}(\mathcal{D}_{\lambda_{\sigma_s}}^{J_{s-1}}(L)))^{\leq \mathbf{h}}\longrightarrow \cdots,\end{align*}
where here we are taking slope decompositions with respect to the Hecke operators at $p$. 
$\lb$
If $h_{\pri} < (k_{\pri}^0+v_{\pri}(\lambda)+1)/e_{\pri}$ for all primes above $p$, it follows that 
\[h_{\pri(\sigma_s)} < \frac{k_{\sigma_s} + v_{\pri(\sigma_s)}(\lambda) + 1}{e_{\pri(\sigma_s)}} = \frac{v_{\pri(\sigma_s)}(\lambda_{\sigma_s})}{e_{\pri(\sigma_s)}}.\]
Now, by Corollary \ref{vanishing} (ii), as $\mathcal{D}_{\lambda_{\sigma_s}}^{J_{s-1}}$ is a quotient of $\mathcal{D}_{\lambda_{\sigma_s}}$, we must have 
\[\h_{\mathrm{c}}^r(Y_1(\n),\sht{\mathcal{D}_{\lambda_{\sigma_s}}^{J_{s-1}}(L)})^{\leq \mathbf{h}} = \{0\}\]
for all $r$. Then, using the long exact sequence, for all $r$ we have
\[\hct^{r}(Y_{1}(\mathfrak{n}), \mathcal{L}_{2}(\mathcal{D}^{J_{s-1}}_\lambda(L)))^{\leq \mathbf{h}}\cong  \hct^{r}(Y_{1}(\mathfrak{n}), \mathcal{L}_{2}(\mathcal{D}^{J_{s}}_\lambda(L)))^{\leq \mathbf{h}},\]
as required.
\end{proof}

\begin{proof}(Theorem \ref{comparaison}). Recall that we defined $V_{\lambda,\mathrm{loc}}(L) \subset \mathcal{A}(L)$ to be the subspace of functions which are locally polynomial of degree at most $\kk$. We see that $V_{\lambda,\mathrm{loc}}(L) \cong \varprojlim V_{\lambda,n}(L),$ where $V_{\lambda,n}(L) \defeq \mathcal{A}_{\lambda,n}(L)\cap V_{\lambda,\mathrm{loc}}(L).$ Note that $V_\lambda(L) = V_{\lambda,0}(L)$. In particular, using part (ii) of Lemma \ref{lemma finite}, we have
\[\hc(Y_1(\n),\sht{V_{\lambda,\mathrm{loc}}(L)^*})^{\leq\mathbf{h}} \cong \hc(Y_1(\n),\sht{V_{\lambda}(L)^*})^{\leq\mathbf{h}}.\]
Hence it suffices to prove the theorem by considering the coefficients of the target space to be in $V_{\lambda, \mathrm{loc}}(L)^{\ast}$ instead of $V_{\lambda}(L)^{\ast}$.
$\lb$
We use the lemma. For this, note that $\mathcal{D}^\Sigma_\lambda(L) = V_{\lambda,\mathrm{loc}}(L)^*$ and $\mathcal{D}^{\varnothing}_\lambda(L) = \dist$. A simple induction on $s$ then shows that we have the required isomorphism.
\end{proof}

%
%

\section{Construction of the distribution}\label{constructdist}
\label{consdist}

Let $\Phi$ be a cuspidal eigenform over $F$ that has small slope (in the sense of the previous section). Then via Eichler--Shimura, we can attach to $\Phi$ a small slope $p$-adic classical modular eigensymbol, and using the results of previous sections, we can lift this to a unique small slope overconvergent eigensymbol. In the work of Pollack and Stevens in \cite{PS11} and \cite{PS12}, and the work of the second author in \cite{Wil17}, once one has such a symbol, one can evaluate it at the cycle $\{0\}-\{\infty\}$ to obtain the $p$-adic $L$-function we desire. This, however, relies on the identification of $\oms$ with the space $\mathrm{Hom}_\Gamma(\mathrm{Div}^0(\Proj(F)),\dist)$, an identification that exists only for $q = 1$, that is, for $F = \Q$ or an imaginary quadratic field. To generalise this to the totally real case, in \cite{Bar13} the first author used automorphic cycles, as introduced in Section \ref{autocycles}, writing down overconvergent analogues of the evaluation maps we used with classical coefficients. Here, we generalise his results to the case of general number fields. The notation we use here was fixed in Section \ref{autocycles}.

\subsection{Evaluating overconvergent classes}\label{evalumapsoc}
\begin{longversion}Recall that in Section \ref{evalumaps}, we used automorphic cycles to define evaluation maps on the space of classical modular symbols with complex coefficients. Here, we adapt these evaluation maps to the case of \emph{overconvergent} modular symbols with \emph{$p$-adic} coefficients. As we are considering a different local system on $Y_1(\n)$, this will be slightly different. In the sequel, we will link all of the various evaluation maps together by explicitly examining the interplay between the various local systems.
$\lb$
\end{longversion}
Suppose $\Psi \in \oms.$ Here recall that we consider the local system given by fibres of
\[\GLt(F)\backslash(\GLt(\A_F)\times \dist)/\Omega_1(\n)K_\infty^+Z_\infty \longrightarrow Y_1(\n),\]
where the action is by 
\[\gamma(x,\mu)uk = (\gamma xuk, \mu * u).\]
In this setting, slightly different versions of the evaluation maps will allow us to associate a distribution to such a class.
\subsubsection{Step 1: Pulling back to $X_{\ff}$}
First we pullback along the map $\eta_{\ff}:X_{\ff} \rightarrow Y_1(\n)$. We have
\[\eta_{\ff}^*\Psi \in \hc(X_{\ff},\eta_{\ff}^*\sht{\dist}).\]
We can see (by examining equation (\ref{welldefeq})) that here the local system corresponding to $\LL_{\ff,2}'(\dist) \defeq \eta_{\ff}^*\sht{\dist}$ is given by the fibres of
\[F^\times\backslash (\A_F^\times \times \dist)/U(\ff)F_\infty^1 \longrightarrow X_{\ff},\]
with action
\[\gamma(x,\mu)ur = \left(\gamma xur, \mu * \matrd{u}{((u-1)\pi_{\ff}^{-1})_{v|p}}{0}{1}\right).\]

\subsubsection{Step 2: Twisting the action}
Unlike in the complex case described earlier, the action describing the local system above  is \emph{not} a nice action, so we twist to get a nicer action of units. To this end, the matrix
\begin{align*}
\matrd{1}{-1}{0}{(\pi_{\ff})_{v|p}} &\in \GLt\bigg(\prod_{\pri|p}F_{\pri}\bigg)\\
&= \GLt(F\otimes_{\Q}\Qp)
\end{align*}
lies in $\Sigma_0(p)$. So we twist our local system by this; denote this twist on distributions by
\begin{align*}
\zeta: \dist &\longrightarrow \dist,\\
\mu &\longmapsto \mu * \matrd{1}{-1}{0}{(\pi_{\ff})_{v|p}},
\end{align*}
and consider
\[\zeta_*\eta_{\ff}^*\Psi \in \hc(X_{\ff},\LL_{\ff,2}(\dist)),\]
where now the local system $\LL_{\ff,2}(\dist)$ is given by
\[F^\times\backslash(\A_F^\times \times \dist)/U(\ff)F_\infty^1 \longrightarrow X_{\ff},\]
\[\gamma(x,\mu)ur = \left(\gamma xur, \mu * \matrd{u}{0}{0}{1}\right).\]

\subsubsection{Step 3: Passing to individual components}
In identical fashion to Section \ref{indc}, we pull back under the isomorphism $\tau_{a_{\y}} : E(\ff)F_\infty^1\backslash F_\infty^+ \isorightarrow X_{\y} \longhookrightarrow X_{\ff}$ given by multiplication by $a_{\y}$. Then we have
\[\tau_{a_{\y}}^*\zeta_*\eta_{\ff}^*\Psi \in \hc(E(\ff)F_\infty^1\backslash F_\infty^+, \LL_{\ff,\y,2}(\dist)),\]
where the local system $\LL_{\ff,\y,2}(\dist)$ is given by
\[E(\ff)F_\infty^1\backslash(F_\infty^+ \times \dist) \longrightarrow E(\ff)F_\infty^1\backslash F_\infty^+,\]
\[er(z,\mu) = \left(erz, \mu*\matrd{e^{-1}}{0}{0}{1}\right).\]
(Note here that whilst $u \in U(\ff)$ acts as $\smallmatrd{u}{0}{0}{1}$, in this step we now have an inverse. This because $u$ is considered as an element of the finite ideles whilst we instead see $e$ as a diagonal infinite idele, which is equivalent under multiplication by $F^\times$ to $e^{-1}$ as a diagonal finite idele and thus an element of $U(\ff)$).

\subsubsection{Step 4: Restricting the coefficient system}
We would like a constant local system. This would allow us to evaluate the cohomology class easily. We see that if we restrict to a quotient of $\dist$ such that, for all $e \in E(\ff)$, the matrix $\smallmatrd{e}{0}{0}{1}$ acts trivially, then we have precisely this. With this in mind, we make the following definitions:
\begin{mdef}
\begin{itemize}
\item[(i)] Define $\locplusf$ to be the subspace of $\loc$ given by
\[\locplusf \defeq \left\{f \in \loc: \matrd{e}{0}{0}{1}*f = f\hspace{3pt} \forall e \in E(\ff)\right\}.\]
Note that equivalently this is the set of all $f \in \loc$ such that $f(ez) = e^{\kk+\vv}f(z)$. 
\item[(ii)] Define $\distplusf$ to be the topological dual of $\locplusf$. Note that $\distplusf$ is a quotient of $\dist$.
\end{itemize}
(Henceforth, we'll drop $\ff$ from the notation, as the level will be clear from context).
\end{mdef}
Now, if we pushforward via the map
\begin{align*}\nu: \dist &\longrightarrow \distplus,\\
\mu &\longmapsto \mu|_{\locplus},
\end{align*}
then the resulting local system is constant. We see that
\begin{align*}
\nu_*\tau_{a_{\y}}^*\zeta_*\eta_{\ff}^*\Psi \in &\hc(E(\ff)F_\infty^1\backslash F_\infty^+, \distplus)\\
&\cong \distplus,
\end{align*}
where the isomorphism is given by integrating over $E(\ff)F_\infty^1\backslash F_\infty^+.$
\subsubsection{Definition of the evaluation map}
\begin{mdef}
We write $\evy$ for the composition
\[\evy : \oms \longrightarrow \distplus\]
of the maps
\begin{align*}
\hc(Y_1(\n),&\sht{\dist}) \labelrightarrow{\zeta_*\eta_{\ff}^*}\hc(X_{\ff},\LL_{\ff,2}(\dist)) \labelrightarrow{\tau_{a_{\y}}^*} \cdots \\
&\hc(\autoc,\LL_{\ff,\y,2}(\dist) \labelrightarrow{\nu_*} \hc(\autoc,\distplus) \cong \distplus.
\end{align*}
\end{mdef}
In particular, we have maps $\evy$ for each $\y \in \clgp{\ff}$. Note that these maps are dependent on the choice of representatives. In any case, for a fixed choice of representatives $\{a_{\y} \in \A_F^\times: \y \in \clgp{\ff}\}$, we have now defined a map
\[\oms \labelrightarrow{\oplus_{\y}\evy} \bigoplus_{\clgp{\ff}}\distplus.\]
\begin{longversion}
To summarise the above construction: to an overconvergent modular symbol, for a fixed ideal $\ff$, we attach a collection of distributions (indexed by $\clgp{\ff}$) with a compatible action of the totally positive units of $F$. This construction depends on a choice of idelic representatives for $\clgp{\ff}$.
\end{longversion}

%
%
\subsection{Locally analytic functions on $\clgp{p^\infty}$}
\begin{longversion}
The $p$-adic $L$-function should be a locally analytic distribution on the ray class group $\clgp{p^\infty}$. Before constructing such a distribution, we must take a digression to describe what locally analytic functions on this space actually look like.
\end{longversion}
\begin{shortversion}
Let $L$ be a (not necessarily finite) extension of $\Qp$ contained in $\Cp$, the completion of an algebraic closure of $\Qp$. Denote by $\mathcal{A}(\clgp{p^\infty},L)$ the space of locally analytic functions on $\clgp{p^\infty}$ defined over $L$, and denote by $\mathcal{D}(\clgp{p^\infty},L)$ its topological dual over $L$. The $p$-adic $L$-function should be an element of this space of distributions; we now give some properties of locally analytic functions that will be required in the sequel.
\end{shortversion}
\subsubsection{The geometry of $\clgp{p^\infty}$}
We first recall the geometry of $\clgp{p^\infty}.$ It is defined as follows:
\[\clgp{p^\infty} \defeq F^\times\backslash\A_F^\times/U(p^\infty)F_\infty^+.\]
Letting $\ff$ range over all ideals dividing $(p)^\infty$ and taking the inverse limit of the series of exact sequences
\[\roiplus \longrightarrow (\roi_F/\ff)^\times \longrightarrow \clgp{\ff} \longrightarrow \nclgp \longrightarrow 0,\]
we see that we have an exact sequence
\[\overline{\roiplus} \longrightarrow (\roifp)^\times \longrightarrow \clgp{p^\infty} \longrightarrow \nclgp \longrightarrow 0,\]
so that -- after picking a choice of representatives for $\nclgp$ -- we have
\[\clgp{p^\infty} \cong \bigsqcup_{\nclgp}(\roifp)^{\times}/\overline{E(1)}.\]
(Here note that $E(1) = \roiplus$, and we have taken $\overline{E(1)}$ to be the image of $E(1)$ in $(\roifp)^\times$). Indeed, for any $\ff$, we can go further, and write
\[\clgp{p^\infty} \cong \bigsqcup_{\y \in \clgp{\ff}}G_{\y},\]
where
\begin{equation}\label{gy}
G_{\y} \defeq \{z \in \clgp{p^\infty}: z \mapsto \y \text{ under the map }\clgp{p^\infty} \rightarrow \clgp{\ff}\}.
\end{equation}
Note that multiplication by $a_{\y}^{-1}$ gives an isomorphism
\[G_{\y} \cong G \defeq \{z \in (\roifp)^{\times}: z \equiv 1 \newmod{\ff}\}/\overline{E(\ff)}.\]

\begin{longversion}
\subsubsection{Analytic functions}
Let $L$ be a (not necessarily finite) extension of $\Qp$ contained in $\Cp$, the completion of the algebraic closure of $\Qp$. Suppose the $L$ is large enough to contain the image of all completions of $F$ at primes above $p$ under their embeddings into $\Cp$.
\begin{mdef}\label{locan}A locally analytic function on $\roifp$ defined over $L$ is a function $f: \roifp \rightarrow L$, such that for each $z_0 \in \roifp$, there is a neighbourhood $U \subset \roifp$ of $z_0$ such that
\[f\bigg|_U(z) = \sum_{\mathbf{r}\geq 0}a_{\mathbf{r}}(z-z_0)^{\mathbf{r}}, \hsp a_{\mathbf{r}} \in L,\]
where $\mathbf{r} \in \Z[\Sigma]$.
\end{mdef}
For a choice of idelic representatives $\{a_i\} \subset \A_F^\times$ of $\nclgp$, we can consider any function
\[\varphi : \clgp{p^\infty} \longrightarrow L\]
as a collection
\[\varphi_{a_i} : (\roifp)^\times/\overline{E(1)} \longrightarrow L,\]
where
\[\varphi_{a_i}(z) \defeq \varphi (a_i^{-1}z).\]
Then, in a slight abuse of notation, $\varphi_{a_i}$ can be thought of as a function
\[\varphi_{a_i} : \roifp \longrightarrow L\]
with support on $(\roifp)^\times$ and with $\varphi_{a_i}(ez) = \varphi_{a_i}(z)$ for all $e \in \overline{E(1)}$. 
\begin{mdef}
We say that $\varphi$ is \emph{locally analytic on $\clgp{p^\infty}$} if each $\varphi_{a_i}$ is locally analytic on $\roifp$.
\end{mdef}
This is independent of the choice of class group representatives. Before we prove this, we rephrase the condition slightly. If $U \subset \roifp$ is an open set, then we say a function $f : U \rightarrow L$ is \emph{analytic} if it can be written as a single convergent power series. By definition, as $\roifp$ is compact, $\varphi_{a_i}$ is locally analytic as a function on $\roifp$ if and only if there exists an ideal $\ff|p^\infty$ such that $\varphi_{a_i}$ is analytic on each of the sets $\{a + \ff(\roifp): a \in \roi_F\otimes\Zp\}$. To make this more precise: given a function $\varphi: \clgp{p^\infty} \rightarrow L$, and an ideal $\ff|p^\infty$, pick some choice of idelic representatives $a_{\y}$ for $\clgp{\ff},$ and define
\begin{align*}\varphi_{a_{\y}} : G &\longrightarrow L,\\
z &\longmapsto \varphi(a_{\y}^{-1}z),
\end{align*}
where $G \cong G_{\y}$ is as above. Also as above, $\varphi_{a_{\y}}$ can be viewed as a function on $\roi_F\otimes_{\Z}\Zp$. We say that $\varphi$ is locally analytic if there exists some $\ff$ such that $\varphi_{a_{\y}}$ is analytic for each $\y$.
\begin{mprop}
Let $\varphi: \clgp{p^\infty} \rightarrow L$ be a function. Then in the above construction, if $\varphi$ is locally analytic for some choice of idelic class group representatives for $\clgp{\ff}$, then it is locally analytic for \emph{any} choice of representatives.
\end{mprop}
\begin{proof}
Suppose there exists an $\ff$ and a set $\{a_{\y}\}$ of representatives such that each $\varphi_{a_{\y}}$ is analytic. Consider a different choice of representatives $\{a_{\y}'\}$. Then, for a fixed $\y$, we have $a_{\y}' = a_{\y}\gamma ur$, where $\gamma \in F^\times, u \in U(\ff),$ and $r \in F_\infty^+.$ Then
\begin{align*}\varphi_{a_{\y}'}(z) = \varphi((a_{\y}')^{-1}z) &= \varphi(a_{\y}^{-1}u^{-1}\gamma^{-1}r^{-1}z)\\
&= \varphi(a_{\y}^{-1}u^{-1}z)\\
&= \varphi_{a_{\y}}(\widetilde{u}^{-1}z),
\end{align*}
where $\widetilde{u}$ is the image of $u \in U(\ff)$ in $U(\ff)/U(p^\infty) \widetilde{\subset} (\roifp)^\times.$ Then, if $\varphi_{a_{\y}}(z) = \sum a_{\mathbf{r}}z^{\mathbf{r}}$, we have
\[\varphi_{a_{\y}'}(z) = \sum(a_{\mathbf{r}}\widetilde{u}^{-r})z^{\mathbf{r}}.\]
As we picked $L$ large enough to contain the image of $\widetilde{u}$ under any embedding into $\Cp$, and $\widetilde{u}$ is a unit, this is a convergent power series over $L$.
$\lb$
Thus $\varphi_{a_{\y}}$ is analytic if and only if $\varphi_{a_{\y}'}$ is analytic, and as $\y$ was arbitrary, this proves the proposition.
\end{proof}
\begin{mremnum}
Note that we have the dictionary 
\[\varphi_{a_{\y}}(z) = \varphi_{a_{\y}'}(\widetilde{u}z)\]
as functions $\roifp \rightarrow L$, where $a_{\y}' = ua_{\y}\gamma r$. We will use this later to prove that the distribution we obtain is independent of choices.
\end{mremnum}
\begin{mdef}We write $\mathcal{A}(\clgp{p^\infty},L)$ for the space of locally analytic functions on $\clgp{p^\infty}$ defined over $L$. We also define the space of \emph{locally analytic distributions} on the ray class group to be
\[\rayclgp \defeq \Hom_{\mathrm{cts}}(\mathcal{A}(\clgp{p^\infty},L),L).\]
\end{mdef}
\end{longversion}

\begin{shortversion}
\subsubsection{Properties of locally analytic functions}
For a choice of idelic representatives $\{a_{\y}\} \subset \A_F^\times$ of $\clgp{\ff}$, we can consider any function
\[\varphi : \clgp{p^\infty} \longrightarrow L\]
as a collection $\{\varphi_{a_{\y}}: \y \in \clgp{\ff}\}$, for
\begin{align*}
\varphi_{a_{\y}} : G &\longrightarrow L,\\
z &\longmapsto \varphi(a_{\y}^{-1}z).
\end{align*}
Then, in a slight abuse of notation, $\varphi_{a_{\y}}$ can be thought of as a function $\varphi_{a_{\y}} : \roifp \longrightarrow L$ with support on a subset of $(\roifp)^\times$ and with $\varphi_{a_{\y}}(ez) = \varphi_{a_{\y}}(z)$ for all $e \in \overline{E(1)}$. A simple calculation then shows:
\begin{mprop}
Suppose $a_{\y}' = a_{\y}\gamma ur$ is a different representative of the class $\y\in\clgp{\ff}$, where $\gamma \in F^\times, u\in U(\ff)$ and $r\in F^+_\infty$. Then $\varphi_{a_{\y}}(z) = \varphi_{a_{\y}'}(\widetilde{u}z)$ as functions on $G$, where $\widetilde{u}$ is the image of $u\in U(\ff)$ in $U(\ff)/U(p^\infty)\subset(\roi_F\otimes_{\Z}\Zp)^\times.$
\end{mprop}
\end{shortversion}

\subsection{Constructing $\mu_{\Psi}$ in $\rayclgp$}
\begin{longversion}With our family of maps $\evy$ from $\oms$ to $\distplus$, we can construct a candidate distribution for the $p$-adic $L$-function. Fix an ideal $\ff|p^\infty$.
\end{longversion}
\begin{mnot}We write $A_{\ff} = \{a_{\y}\}$ to denote our system of class group representatives for $\clgp{\ff}$. 
\end{mnot}
We now construct a distribution $\mu_{\Psi}^{A_{\ff}}$ associated to \emph{this} choice of representatives. Let $\varphi$ be a locally analytic function on $\clgp{p^\infty}$. Via the above construction, we obtain functions $\varphi_{a_{\y}}: G_{\y} \rightarrow L$, each of which we can view as a function
\[\varphi_{a_{\y}}: \roifp \longrightarrow L\]
with support on the open subset $\{z \in (\roifp)^\times: z \equiv 1 \newmod{\ff}\}$ and satisfying
\[\varphi_{a_{\y}}(ez) = \varphi_{a_{\y}}(z) \hsp \forall e \in E(\ff).\]
Now, $\evy(\Psi) \in \distplus$. This is a distribution that takes as input functions $\psi: (\roifp)^\times \rightarrow L$ with $\psi(ez) = e^{\kk+\vv}\psi(z).$ To force $\varphi_{a_{\y}}$ to satisfy this condition, we twist it.
\begin{mdef}
If $\psi : \roifp \rightarrow L$ is a function with support on elements congruent to $1 \newmod {\ff}$ and that satisfies $\psi(ez) = \psi(z)$ for all $e \in E(\ff)$, then we define $\psi^* \in \locplus$ by 
\[\psi^*(z) = \left\{\begin{array}{ll}z^{\kk+\vv}\psi(z^{-1}) &: z \in (\roi_F\otimes_{\Z}\Zp)^\times, \\
0 &: \text{ otherwise}.\end{array}\right.\]
Since $\psi$ has support inside the units, this remains continuous. It is simple to see that this now satisfies the condition required. We use $z^{-1}$ rather than $z$ for reasons of compatibility in later calculations.
\end{mdef}
Now we can evaluate $\evy(\Psi)$ at $\varphi_{a_{\y}}^*$. This motivates:
\begin{mdef}\label{defndist}
Define $\mu_{\Psi}^{A_{\ff}} \in \rayclgp$ by 
\[\mu_{\Psi}^{A_{\ff}}(\varphi) = \sum_{\y \in \clgp{\ff}}\evy(\Psi)(\varphi_{a_{\y}}^*) \in L.\]
\end{mdef}
\begin{mprop}For fixed $\ff$, this is independent of the choice of class group representatives.
\end{mprop}
\begin{proof}
There are two layers to this. Choosing representatives fixes:
\begin{itemize}
\item[(a)] The collection of maps $\{\evy(\Psi):a_{\y} \in A_{\ff}\}$, and
\item[(b)] The identification of $\varphi$ with $(\varphi_{a_{\y}})_{\y \in \clgp{\ff}}.$
\end{itemize}
We prove that these choices cancel each other out. To do so, we examine the local systems; see Section \ref{autocycles} for descriptions of each local system.
$\lb$
Recall that we have $\zeta_*\eta_{\ff}^*\Psi \in \hc(X_{\ff},\LL_{\ff,2}(\dist))$ (canonically), and then that we can pull back to $X_{\y}$ under the canonical inclusion. At the first stage where our representatives come into play, the map of local systems induced by 
\[\tau_{a_{\y}}: E(\ff)F_\infty^1\backslash F_\infty^+ \isorightarrow X_{\y}\]
can be described by the map
\begin{align}\label{localsystem}
F^\times\backslash( F^\times a_{\y} U(\ff)F_\infty^+ \times \dist)/U(\ff)F_\infty^1 &\longrightarrow E(\ff)F_\infty^1\backslash(F_\infty^+ \times \dist)\\
(\gamma a_{\y}ur, \mu) &\longmapsto \left(r,\mu*\matrd{\widetilde{u}^{-1}}{0}{0}{1}\right),\notag
\end{align}
recalling that $\tau_{a_{\y}}$ is given by $z \mapsto a_{\y}z$ and that $\widetilde{u}$ is the image of $u$ in $(\roifp)^\times$. This map is well-defined; indeed, consider
\begin{align*}\gamma' [(\gamma a_{\y}ur,\mu)]vs &= \left[\left(\gamma'\gamma a_{\y}uvrs, \mu * \smallmatrd{\widetilde{v}}{0}{0}{1}\right)\right]\\
& \longmapsto \left[\left(rs, \left(\mu *\smallmatrd{\widetilde{v}}{0}{0}{1}\right)* \smallmatrd{(\widetilde{uv})^{-1}}{0}{0}{1}\right)\right]\\
& = \left[\left(r,\mu*\smallmatrd{\widetilde{u}^{-1}}{0}{0}{1}\right)\right] = \mathrm{Im}([\gamma a_{\y}ur, \mu]).
\end{align*}
Now suppose we choose a different set of representatives $\{a_{\y}'\}$, with, as before, 
\[a_{\y}' = a_{\y}\gamma ur, \hsp \gamma \in F^\times, u \in U(\ff), r \in F_\infty^1.\]
Then under the map of equation (\ref{localsystem}), we have
\[[(a_{\y}',\mu)] = [(a_{\y}\gamma ur,\mu)] \longmapsto \left[\left(r,\mu*\smallmatrd{\widetilde{u}^{-1}}{0}{0}{1}\right)\right].\]
Thus, when we restrict, we find that
\[\mathrm{Ev}_{\ff,\dagger}^{a_{\y}'}(\Psi) = \evy(\Psi)*\matrd{\widetilde{u}^{-1}}{0}{0}{1}.\]
We have already shown that, for $\varphi \in \mathcal{A}(\clgp{p^\infty},L),$ we have 
\[\varphi_{a_{\y}'}(\widetilde{u}z) = \varphi_{a_{\y}}(z).\]
\begin{shortversion}
Then an easy calculation shows that
\[\varphi_{a_{\y}'}^*(z) = \matrd{\widetilde{u}}{0}{0}{1}*\varphi_{a_{\y}}^*(z).\]
\end{shortversion}
\begin{longversion}
Then
\begin{align*}\varphi_{a_{\y}}^*(z) &= z^{\kk+\vv}\varphi_{a_{y}}(z)^{-1}\\
&= z^{\kk+\vv}\varphi_{a_{\y}'}(\widetilde{u}z)^{-1}\\
&= \widetilde{u}^{-\kk-\vv}(\widetilde{u}z)^{\kk+\vv}\varphi_{a_{\y}'}(\widetilde{u}z)^{-1}\\
& = \widetilde{u}^{-\kk-\vv}\varphi_{a_{\y}'}^*(\widetilde{u}z)\\
&= \matrd{\widetilde{u}^{-1}}{0}{0}{1}*\varphi_{a_{\y}'}^*(z).
\end{align*}
Thus
\[\varphi_{a_{\y}'}^*(z) = \matrd{\widetilde{u}}{0}{0}{1}*\varphi_{a_{\y}}^*(z).\]
\end{longversion}
Accordingly,
\begin{align*}\mathrm{Ev}_{\ff,\dagger}^{a_{\y}'}(\Psi)(\varphi_{a_{\y}'}^*) &= \evy(\Psi)*\matrd{\widetilde{u}^{-1}}{0}{0}{1}\left(\matrd{\widetilde{u}}{0}{0}{1}* \varphi_{a_{\y}}^*\right)\\
&= \evy(\Psi)(\varphi_{a_{\y}}^*).
\end{align*}
Thus this is independent of the choice of representatives, as desired.
\end{proof}

\begin{mdef}For some choice of representatives $A_{\ff} = \{a_{\y}\}$ of $\clgp{\ff}$, define
\[\mu_{\Psi}^{\ff} \defeq \mu_{\Psi}^{A_{\ff}}.\]
(Note that, by the proposition, this is well-defined for each $\ff$).
\end{mdef}

%
%
\subsection{Compatibility over choice of $\ff$}
We have defined, for each $\ff|p^\infty$, a distribution $\mu_{\Psi}^{\ff} \in \rayclgp.$ We now investigate how this distribution varies with the choice of $\ff$. Since we have independence of choice, we now choose class group representatives that are compatible in the following sense.
\begin{mnot}Throughout this section, take $\ff|p^\infty$ and let $\pri|p$ be a prime. We will make the following important assumption throughout this section:
\[\textit{The ideal }\ff\textit{ is divisible by all of the primes above }p.\]
Let $A_{\ff} =\{a_{\y}\}$ be a full set of representatives for $\clgp{\ff}$, and let $\{u_r \in U(\ff): r\in R\}$, for $R = U(\ff)/E(\ff)U(\ff\pri)$, be elements of $U(\ff)$ such that the set
\[A_{\ff\pri} \defeq \{a_{\y}u_r: \y \in \clgp{\ff}, r\in R\}\] 
is a full set of representatives for $\clgp{\ff\pri}$.
\end{mnot}

\begin{longversion}
\begin{mlem}\label{compatf}\begin{itemize}
\item[(i)] There is a commutative diagram 
\[\CommDia{\hc(Y_1(\n),\sht{\dist})}{U_{\pri}}{\hc(Y_1(\n),\sht{\dist})}
{\zeta_*\eta_{\ff\pri}^*}{\zeta_*\eta_{\ff}^*}
{\hc(X_{\ff\pri},\LL_{\ff\pri,2}(\dist))}{\mathrm{Tr}}{\hc(X_{\ff},\LL_{\ff,2}(\dist))},\]
where the bottom map is the natural trace map on cohomology (see, for example, \cite{Hid93}, Section 7).

\item[(ii)] Write $\mathrm{Ev}_{\ff,\dagger}^{A_{\ff}} = \oplus_{\y \in \clgp{\ff}}\evy$, and similarly for $\ff\pri$ with respect to $A_{\ff\pri}$. Then we have the commutative diagram
\[\TallCommDia{\hc(Y_1(\n),\sht{\dist})}{\mathrm{Ev}_{\ff\pri,\dagger}^{A_{\ff\pri}}}{\bigoplus_{\clgp{\ff\pri}}\distplus}{U_{\pri}}{\mathrm{tr}_{\ff}}{\hc(Y_1(\n),\sht{\dist})}{\mathrm{Ev}_{\ff,\dagger}^{A_{\ff}}}{\bigoplus_{\clgp{\ff}}\distplus},\]
where
\[\mathrm{tr}_{\ff}\bigg((\mu_{a_{\y}u_r})_{\y,r}\bigg) = \left(\sum_{r\in R}\mu_{a_{\y}u_r}\bigg|\matrd{\widetilde{u_r}^{-1}}{0}{0}{1}\right)_{\mathbf{y}\in\clgp{\ff}}.\]
\end{itemize}
\end{mlem}
\begin{proof}
We construct commutative diagrams at each step in the definition of the evaluation maps. For convenience, we drop $\dist$ from the notation and instead write $\LL = \LL(\dist)$, for appropriate subscripts on $\LL$. 
$\lb$
Note that there is a natural projection map $\mathrm{pr}:X_{\ff\pri} \rightarrow X_{\ff}$ that induces a trace map on the cohomology. We wish to construct a map $\hc(X_{\ff\pri},\LL_{\ff\pri,2}') \rightarrow \hc(X_{\ff},\LL_{\ff,2}')$, and for this it hence suffices to construct a map of sheaves
\[\mathrm{pr}_*\LL_{\ff\pri,2}' \longrightarrow \LL_{\ff,2}'.\]
We do this as follows.  Note that there is a natural map $\alpha: \LL_{\ff\pri,2}' \rightarrow \LL_{\ff,2}'$ given by 
\[(x,\mu) \longmapsto \left(x,\mu * \matrd{1}{0}{0}{\pi_{\pri}}\right).\]
Let $G\defeq \mathrm{Gal}(X_{\ff\pri}/X_{\ff}) \cong \roi_F/\pri.$ Now let $\mathcal{U} \subset X_{\ff}$ be open and sufficiently small that 
\[\mathrm{pr}^{-1}(\mathcal{U}) = \bigsqcup_{g\in G}\mathcal{U}_g \subset X_{\ff\pri},\]
where $\mathrm{pr}$ induces a homeomorphism $i_g: \mathcal{U} \rightarrow \mathcal{U}_g.$ Then 
\[\mathrm{pr}_*\LL_{\ff,2}'(\mathcal{U}) \defeq \LL_{\ff\pri,2}'(\mathrm{pr}^{-1}(\mathcal{U})) = \bigoplus_{g\in G}\LL_{\ff\pri,2}'(\mathcal{U}_g).\]
Then define a map
\begin{align*}
\mathrm{pr}_*\LL_{\ff\pri,2}'(\mathcal{U}) &\longrightarrow \LL_{\ff,2}'(\mathcal{U}),\\
s = (s_g)_{g\in G} &\longmapsto \sum \alpha\circ s_g \circ i_g.
\end{align*}
We see that we have defined a map
\[\hc(X_{\ff\pri},\LL_{\ff\pri,2}') \labelrightarrow{*\smallmatrd{1}{0}{0}{\pi_{\pri}}} \hc(X_{\ff},\LL_{\ff,2}'),\]
and moreover we see that this makes the diagram
\[\CommDia{\hc(Y_1(\n),\sht{\dist})}{U_{\pri}}{\hc(Y_1(\n),\sht{\dist})}
{\eta_{\ff\pri}^*}{\eta_{\ff}^*}
{\hc(X_{\ff\pri},\LL_{\ff\pri,2}')}{*\smallmatrd{1}{0}{0}{\pi_{\pri}}}{\hc(X_{\ff},\LL_{\ff,2}')}
\]
commute.
$\lb$
Via a similar construction, and by replacing the map $\alpha$ with the map $\alpha': \LL_{\ff\pri,2} \rightarrow \LL_{\ff,2}$ defined by $(x,\mu) \mapsto (x,\mu)$, we see that we construct a map $\hc(X_{\ff\pri},\LL_{\ff\pri,2}) \rightarrow \hc(X_{\ff},\LL_{\ff,2})$ that is nothing but the trace map on cohomology, and in particular that we have the following commutative diagram:
\[\CommDia{\hc(X_{\ff\pri},\LL_{\ff\pri,2}')}{*\smallmatrd{1}{0}{0}{\pi_{\pri}}}{\hc(X_{\ff},\LL_{\ff,2}')}
{\zeta_*}{\zeta_*}
{\hc(X_{\ff\pri},\LL_{\ff\pri,2})}{\mathrm{tr}}{\hc(X_{\ff},\LL_{\ff,2})}.
\]
Finally, we bring in the dependence on our choice of class group representatives. We see that there is commutative diagram
\[\BigCommDia{\hc(X_{\ff\pri},\LL_{\ff\pri,2})}{\mathrm{tr}}{\hc(X_{\ff},\LL_{\ff,2})}
{\tau_{a_{\y}u_r}^*}{\tau_{a_{\y}}^*}
{\hc(E(\ff\pri)F_\infty^1\backslash F_\infty^+,\LL_{\ff\pri,\mathbf{x},2})}{\left[* \smallmatrd{\widetilde{u_r}^{-1}}{0}{0}{1}\right]_*\circ \mathrm{pr}_*}{\hc(\autoc,\LL_{\ff,\yy,2})},
\]
where $\mathbf{x}$ denotes the class in $\clgp{\ff\pri}$ represented by $a_{\y}u_r$, and $\mathrm{pr}$ is the natural projection map. Here we have used the results of the previous section. This shows that for $\Psi \in \oms$, we have
\[\mathrm{Ev}_{\ff\pri,\dagger}^{a_{\y}u_r}(\Psi)* \matrd{\widetilde{u_r}^{-1}}{0}{0}{1} = \mathrm{Ev}_{\ff,\dagger}^{a_{\y}}(\Psi|U_{\pri})\bigg|_{G_{\mathbf{x}}},\]
where $G_{r} = \{z \in \roifp: z \cong u_re \newmod{\ff\pri}\text{ for some }e\in E(\ff)\}.$ Summing over the relevant narrow ray class groups gives the diagram as stated.
\end{proof}

\end{longversion}

\begin{shortversion}
\begin{mlem}\label{compatf}\begin{itemize}
\item[(i)] There is a commutative diagram 
\[\CommDia{\hc(Y_1(\n),\sht{\dist})}{U_{\pri}}{\hc(Y_1(\n),\sht{\dist})}
{\zeta_*\eta_{\ff\pri}^*}{\zeta_*\eta_{\ff}^*}
{\hc(X_{\ff\pri},\LL_{\ff\pri,2}(\dist))}{\mathrm{Tr}}{\hc(X_{\ff},\LL_{\ff,2}(\dist))},\]
where the bottom map is the natural trace map on cohomology (see, for example, \cite{Hid93}, Section 7).

\item[(ii)]We have, for $\Psi \in \oms$, the relation
\[\mathrm{Ev}_{\ff\pri,\dagger}^{a_{\y}u_r}(\Psi)*\matrd{\widetilde{u_r}}{0}{0}{1} = \evy(\Psi|U_{\pri})\bigg|_{G_r},\]
where \[G_r \defeq \{z \in \roifp: \text{ there exists }e \in E(\ff)\text{ such that }ez \equiv u_r \newmod{\ff\pri}\}.\]
\end{itemize}
\end{mlem}
\begin{proof}
For part (i), see \cite{Bar13}, Lemme 5.2; the proof generalises immediately to the general number field setting. For part (ii), we bring in our explicit dependence on class group representatives. In particular, note that there is a commutative diagram
\[\BigCommDia{\hc(X_{\ff\pri},\LL_{\ff\pri,2}(\dist))}{\mathrm{Tr}}{\hc(X_{\ff},\LL_{\ff,2}(\dist))}
{\nu_*(\tau_{a_{\y}u_r}^{\ff\pri})^*}{\nu_*(\tau_{a_{\y}u_r}^{\ff})^*}
{\mathcal{D}_{\lambda}^{\ff\pri,+}(L)}{\text{restriction to }G_r}{\distplusf},\]
where we have written $(\tau_{a_{\y}u_r}^{\ff})^*$ to emphasise the dependence of this map on the ideal. Hence 
\[\mathrm{Ev}_{\ff\pri,\dagger}^{a_{\y}u_r}(\Psi) = \mathrm{Ev}_{\ff,\dagger}^{a_{\y}u_r}(\Psi|U_{\pri})\bigg|_{G_r}.\]
Using the results of the previous section, we have the equality
\[\evy(\Psi|U_{\pri}) = \mathrm{Ev}_{\ff,\dagger}^{a_{\y}u_r}(\Psi|U_{\pri})*\matrd{\widetilde{u_r}}{0}{0}{1},\]
 hence the result.

\end{proof}

\end{shortversion}
\begin{mprop}
Let $\ff|p^\infty$ be divisible by all of the primes above $p$, and let $\pri$ be a prime above $p$. Let 
\[\Psi \in \hc(Y_1(\n),\sht{\dist})\]
be an eigensymbol for all the Hecke operators at $p$, with $U_{\pri}$-eigenvalue $\lambda_{\pri}.$ Then
\[\mu_{\Psi}^{\ff\pri} = \lambda_{\pri}\mu_{\Psi}^{\ff}.\]
\end{mprop}
\begin{proof}
Let $\varphi \in \mathcal{A}(\clgp{p^\infty})$. We evaluate $\mu_{\Psi}^{\ff\pri}$ at $\varphi$ by using the class group representatives $A_{\ff\pri}$, and then evaluate $\mu_{\Psi|U_{\pri}}^{\ff}$ at $\varphi$ using the representatives $A_{\ff}$, and use the previous lemma to show that they are equal.
$\lb$
Fix $\y\in\clgp{\ff}$ and $r \in R$. Then we see that
\[\varphi_{a_{\y}u_r}(z) = \varphi_{a_{\y}}(u_r^{-1}z)\]
for $z\in G_r$. In particular, we have
\[\varphi_{a_{\y}}^*\bigg|_{G_r} = \matrd{\widetilde{u_r}^{-1}}{0}{0}{1} *\varphi_{a_{\y}u_r}^*\left(z\right)\]
Observe now that by the previous lemma, we have
\begin{align*}
\evy(\Psi|U_{\pri})(\varphi_{a_{\y}}^*) &= \sum_{r\in R}\mathrm{Ev}_{\ff\pri}^{a_{\y}u_r}(\Psi)\bigg|\matrd{\widetilde{u_r}}{0}{0}{1} \left(\varphi_{a_{\y}}^*\bigg|_{G_r}\right)\\
&= \sum_{r\in R}\mathrm{Ev}_{\ff\pri}^{a_{\y}u_r}(\Psi)\left(\varphi_{a_{\y}u_r}^*\right).
\end{align*}
Summing over $\y \in \clgp{\ff}$ on both sides, and replacing $\Psi|U_{\pri}$ with $\lambda_{\pri}\Psi$ on the left hand side, now shows the result.
\end{proof}
We have now proved the following:
\begin{mthm}\label{construction}
Let $\Psi \in \hc(Y_1(\n),\sht{\dist})$ be an eigenclass for the $U_{\pri}$ operators for all $\pri|p$, and let $\ff|(p^\infty)$ be some choice of ideal with $\ff$ divisible by all the primes above $p$. Define $U_{\ff} \defeq \prod_{\pri^r || \ff}U_{\pri}^r$, write $\lambda_{\ff}$ for the eigenvalue of $U_{\ff}$, and define
\[\mu_{\Psi} \defeq \lambda_{\ff}^{-1}\mu_{\Psi}^\ff.\]
This is well-defined and independent of choices up to a fixed choice of uniformisers at primes above $p$.
$\lb$
Thus for such $\Psi$ there is way of attaching an element $\mu_{\Psi}$ of $\rayclgp$ to $\Psi$ that is independent of choices.
\end{mthm}

\begin{mdef}
In the set-up of above, we call $\mu_{\Psi}$ the \emph{$p$-adic $L$-function of $\Phi$}.
\end{mdef}

\subsection{Evaluating at Hecke characters}\label{evathc}
Let $\varphi$ be a Hecke character of infinity type $\rr \in \Z[\Sigma]$ and conductor $\ff|(p^\infty)$, where $\ff$ is divisible by every prime above $p$. In this section we describe the evaluation of the distribution $\mu_{\Psi}$ at $\varphi_{p-\mathrm{fin}}$ (as defined in Section \ref{rcg}).
$\lb$
Choosing representatives $\{a_{\y}\}$ for $\clgp{\ff}$, we see that
\[(\varphi_{p-\mathrm{fin}})_{a_{\y}} = \mathbf{1}_{G_{\y}}\epsphi\varphi_f(a_{\y})\mathbf{z}^{\rr},\]
where $\mathbf{1}_{G_{\y}}$ is the indicator function of the open subset of $\clgp{p^\infty}$ corresponding to $\y \in \clgp{\ff}$ (see equation (\ref{gy})), and $\mathbf{z}$ is a variable on $\roi_F\otimes_{\Z}\Zp$. We see that, for $\Psi$ as above,
\begin{equation}\label{expev}
\mu_{\Psi}(\varphi_{p-\mathrm{fin}}) = \lambda_{\ff}^{-1}\sum_{\y}\epsphi\varphi_f(a_{\y})\evy(\Psi)(\mathbf{z}^{\kk+\vv-\rr}).
\end{equation}

%
%

\section{Interpolation of $L$-values}\label{interpolation}
In previous sections, we have defined the maps denoted by solid arrows in the following diagram:
\begin{equation}\label{compat}
\begin{diagram}
\hc(Y_1(\n),\sho{V_\lambda(L)^*}) &&&&\rTo^{\evyclo}&&&& L\\
\rotatebox{270}{\scalebox{3}[1]{$\hspace{-4pt}\cong\hspace{-2pt}$}} &&&&&&&&\dDotsto_{\beta}\\
\hc(Y_1(\n),\sht{V_\lambda(L)^*}) &&&&\rDotsto^{\evyclt}&&&& L\\
\uTo_{\rho} &&&&&&&&\uDotsto_{\delta}\\
\hc(Y_1(\n),\sht{\dist}) &&\rTo^{\evy} && \distplus &&\rTo^{\text{ev.\ at }z^{\kk-\jj}} && L
\end{diagram}
\end{equation}
In particular, the isomorphism is induced by the isomorphism of local systems given in Remark \ref{locsysiso}, the top (classical) evaluation map was defined in Section \ref{evalumaps}, the map $\rho$ is induced from the specialisation $\dist \rightarrow V_\lambda(L)^*,$ and the bottom (overconvergent) evaluation map was defined in Section \ref{evalumapsoc}. In this section, we define the maps above denoted by dotted arrows in a manner such that the diagram commutes. By doing so, we will be able to use our previous results to relate the evaluation of the distribution $\mu_\Phi$ at Hecke characters with critical $L$-values of $\Phi$.

\subsection{Classical evaluations, II}\label{classical 2}
We start by defining the`missing' evaluation map. We have already touched on all of the key points of this construction; it is essentially a blend of our previous two evaluation maps. Taking notation from Section \ref{autocyclesec}, we pullback along $\eta_{\ff}$, giving a local system $\eta_{\ff}^*\sht{V_\lambda(L)^*}$ on $X_{\ff}$ that can be described by sections of the projection
\[F^\times\backslash (\A_F^\times \times V_\lambda(L)^*)/U(\ff)F_\infty^1,\]
with action
\[f(x,P)ur = \left(fxur, P*\matrd{u}{((u-1)\pi_{\ff}^{-1})_{v|\ff}}{0}{1}\right).\]
This bears relation with the overconvergent case, in that we have an action of units that is not particularly nice. As in that case, we `untwist' this action using the map $(\zeta_{\ff})_*$ from Section \ref{evalumapsoc}, so that units act via the matrix $\smallmatrd{u}{0}{0}{1}$. We can then pull-back under the injection 
\[\tau_{a_{\y}} : E(\ff)F_\infty^1\backslash F_\infty^+ \longhookrightarrow X_{\ff}\]
of previous sections. Finally, as in the classical case, we pushforward under evaluation at the polynomial $\XX^{\kk-\jj}\YY^{\jj}$, which lands us in a cohomology group with coefficients in a constant sheaf (see Section \ref{evalumaps}). Combining all of these maps, we get a map
\[\evyclt : \hc(Y_1(\n),\sht{V_\lambda(L)^*}) \longrightarrow L,\]
which gives the definition of the dotted horizontal arrow in the diagram.
$\lb$
The following lemma determines the definition of the map $\beta$ in the diagram. For ease of notation, write $\mathrm{Ev}_k$ for the map $\mathrm{Ev}_{\ff,\jj,k}^{a_{\y}}$.
\begin{mlem}
\label{trackiso}
Let $\alpha$ denote the isomorphism
\[\alpha : \hc(Y_1(\n),\sho{V_\lambda(L)^*}) \isorightarrow \hc(Y_1(\n),\sht{V_\lambda(L)^*})\]
induced by the isomorphism $\sho{V_\lambda(L)^*} \isorightarrow \sht{V_\lambda(L)^*}$ of local systems given by
\[(g,P) \longmapsto (g,P|g_p)\]
(see Remark \ref{locsysiso}). Then 
\[\mathrm{Ev}_2(\alpha(\phi)) = \pi_{\ff}^{\jj+\vv}\mathrm{Ev}_1(\phi).\]
\end{mlem}
\begin{mrem}
Here, in an abuse of notation, we write $\pi_{\ff}$ for the natural element of $L$ corresponding to $(\pi_{\ff})_{v|p} \in \roi_F\otimes_{\Z}\Zp$ under our fixed choice of uniformisers at primes above $p$. Note that under this map, a uniformiser $\pi_{\pri}$ is mapped to $N_{K_{\pri}/\Qp}(\pi_{\pri})$, so that as elements of $L$, we have $\pi_{\ff}^{\jj+\vv}= N(\ff)^{[\jj+\vv]}$ up to multiplication by a $p$-adic unit. In particular, `multiplication by $\pi_{\ff}^{\jj+\vv}$' is a well-defined concept.
\end{mrem}
\begin{proof}
We look at the local systems in each case. A simple check shows that there is a commutative diagram
\[\BigCommDia{\hc(Y_1(\n),\sho{V_\lambda(L)^*})}{\eta_{\ff}^*}{\hc(X_{\ff},\LL_{\ff,1}(V_\lambda(L)^*))}{\alpha}{\alpha'}{\hc(Y_1(\n),\sht{V_\lambda(L)^*})}{(\zeta_{\ff})_*\eta_{\ff}^*}{\hc(X_{\ff},\LL_{\ff,2}(V_{\lambda}(L)^*))},
\]
where $\alpha'$ is the map induced by the map
\[(x,P) \longmapsto \left(x,P\left|\matrd{x_p}{0}{0}{(\pi_{\ff})_{v|\ff}}\right.\right)\]
of local systems. Then continuing, we see that there is a commutative diagram
\[
\CommDia
{\hc(X_{\ff},\LL_{\ff,1}(V_\lambda(L)^*))}
{\tau_{a_{\y}}^*}
{\hc(\autoc,\LL_{\ff,\y,1}(V_\lambda(L)^*))}
{\alpha'}{\alpha''}
{\hc(X_{\ff},\LL_{\ff,2}(V_\lambda(L)^*))}
{\tau_{a_{\y}}^*}
{\hc(\autoc,\LL_{\ff,\y,2}(V_\lambda(L)^*))}
\]
where $\alpha''$ is the map induced by the map
\[(r,P) \longmapsto \left(r,\left|\matrd{1}{0}{0}{(\pi_{\ff})_{v|\ff})}\right.\right)\]
of local systems. Finally, there is a commutative diagram
\[\BigCommDia{\hc(\autoc,\LL_{\ff,\y,1}(V_\lambda(L)^*))}
{(\text{ev.\ at }\XX^{\kk-\jj}\YY^{\jj})_*}
{L}
{\alpha''}{\times \pi_{\ff}^{\jj+\vv}}
{\hc(\autoc,\LL_{\ff,\y,2}(V_\lambda(L)^*))} 
{(\text{ev.\ at }\XX^{\kk-\jj}\YY^{\jj})_*}
{L}.
\]
Putting these diagrams together gives the required result.
\end{proof}
Recall the definition of $\ev_\varphi$ in Definition \ref{evphi}, and relabel $\ev_{\varphi,1} \defeq \ev_\varphi$. Similarly define
\[\ev_{\varphi,2} \defeq \sum_{\y\in\clgp{\ff}}\epsphi\varphi_f(a_{\y})\evyclt,\]
where this makes sense, and note that by an identical argument to previously this is independent of class group representatives. Using the results above with the results in Section \ref{algebraic}, we obtain:
\begin{mcor}\label{fincor}
Recall the definition of $\theta_K \in \hcusp(Y_1(\n),\sho{V_\lambda(K)^*})$ from Definition \ref{thetak}, and recall that we set $\theta_L$ to be its image in $\hc(Y_1(\n),\sht{V_\lambda(L)^*})$ under the inclusions of equation (\ref{inclusion1}) and (\ref{inclusion2}). Then 
\[\ev_{\varphi,2}(\theta_L) = \pi_{\ff}^{\jj+\vv}\ev_{\varphi,1}(\theta_K) = (-1)^{R(\jj,\kk)}\left[\frac{|D|\tau(\varphi)\pi_{\ff}^{\jj+\vv}}{2^{r_2}\Omega_\Phi^{\epsphi}}\right]\cdot \Lambda(\Phi,\varphi),\]
where $R(\jj,\kk) = \sum_{v\in\Sigma(\R)}j_v+k_v + \sum_{v\in\Sigma(\C)}k_v$.
\end{mcor}
Note here that this holds for \emph{any} conductor $\ff|(p^\infty)$, with no condition on ramification.

\subsection{Relating classical and overcovergent evaluations}\label{classical and ovcgt}
Returning to the commutative diagram in equation (\ref{compat}), we now show that the map $\delta$ is actually nothing but the identity map. For a suitable automorphic form $\Phi$, this will then allow us to prove the required interpolation property for the distribution $\mu_{\Phi}$.
\begin{mprop}
There is a commutative diagram
\[\BigCommDia{\hc(Y_1(\n),\sht{\dist})}{\evy}{\distplus}{\rho}{\textit{ev. at }z^{\kk-\jj}}{\hc(Y_1(\n),\poly)}{\evyclt}{L},\]
where the left vertical arrow is the specialisation map and the right vertical arrow is evaluation at the polynomial $z^{\kk-\jj}$.
\end{mprop}
\begin{proof}
This is easily shown by looking at each step of the construction of the maps $\evy$ and $\evyclt$ in the previous sections. At each of steps 1, 2 and 3 we can write down a specialisation map by restricting the coefficients, and by looking at the level of local systems, we can clearly see that these specialisations commute with the maps $\eta_{\ff}$, $\zeta_{\ff}$ and $\tau_{a_{\y}}$. It remains to show compatibility over step 4, where the construction is slightly different. This amounts to showing that the diagram
\[\BigCommDia{\hc(\autoc,\LL_{\ff,\y,2}(\dist))}{\textit{res}}{\distplus}{}{\textit{ev. at }z^{\kk-\jj}}{\hc(\autoc,\LL_{\ff,\y,2}(V_\lambda(L)^*))}{\text{ev. at }\mathbf{X}^{\kk-\jj}\mathbf{Y}^{\jj}}{L}\]
commutes, where the lefthand map is restriction of the coefficients, the map $\textit{res}$ is the restriction of coefficients to $\distplus$ followed by integration over a fixed de Rham cohomology class, and the bottom map is the composition of $(\rho_{\jj})_*$ with integration over the same de Rham cohomology class. Since $\polyt \hookrightarrow \loc$ via $P(\XX,\YY) \mapsto P(z,1)$, we see that when we look at the corresponding local systems, we are evaluating at the same element in each case; thus the diagram commutes.
\end{proof}
By combining this with equation (\ref{expev}) for $\mu_{\Psi}(\varphi_{p-\mathrm{fin}})$, we get the following corollary:
\begin{mcor}\label{fincor2}
Let $\phi \in \hc(Y_1(\n),\sht{V_\lambda(L)^*})$ be a small slope Hecke eigensymbol with $U_{\ff}$-eigenvalue $\lambda_{\ff}$ and with (unique) overconvergent eigenlift $\Psi$, and let $\mu_{\Psi}$ be the corresponding ray class distribution. Then for a Hecke character $\varphi$ of infinity type $\jj+\vv$ and conductor $\ff|(p^\infty)$, where $0\leq \jj \leq \kk$ and $\ff$ is divisible by every prime above $p$, we have 
\[\mu_{\Psi}(\varphi_{p-\mathrm{fin}}) = \lambda_{\ff}^{-1}\ev_{\varphi,2}(\phi).\]
\end{mcor}
In the case that $\phi$ is the modular symbol attached to an automorphic form, this then gives the desired interpolation property at Hecke characters that ramify at all primes above $p$ as an immediate corollary (see Theorem \ref{summary} below).


\subsection{Interpolating at unramified characters}
We now consider interpolation of $L$-values at Hecke characters that are not necessarily ramified at all primes above $p$. For this, we use Corollary \ref{unramified primes}. Whilst the results of this section up until now have been for arbitrary modular symbols, to use this corollary we need to restrict to the case where the cohomology classes we consider are those attached to automorphic forms via the Eichler-Shimura isomorphism. Let $\Phi$ be such an automorphic form of weight $\lambda$ and level $\Omega_1(\n)$, and suppose that $\Phi$ is a Hecke eigenform that has small slope at the primes above $p$. Let $\phi_L$ be the ($p$-adic) modular symbol attached to $\Phi$, and let $\Psi$ be the associated (unique) overconvergent modular symbol corresponding to $\phi_L$ under the control theorem. Then we have the following lemma:
\begin{mlem}\label{interpolation unramified}
Let $\varphi$ be a Hecke character of conductor $\ff|(p^\infty)$ (with no additional conditions on $\ff$) and infinity type $\jj+\vv$, where $0\leq \jj\leq \kk$. Let $B$ be the set of primes above $p$ that do not divide $\ff$, and define $\ff' \defeq \ff\prod_{\pri\in B}\pri,$ so that $\ff'$ is divisible by all the primes above $p$. Then we have
\begin{align}\label{expev unramified}
\mu_\Psi(\varphi_{p-\mathrm{fin}}) &= \lambda_{\ff'}^{-1}\pi_{\ff'}^{\jj+\vv}\left[\prod_{\pri \in B}(\varphi(\pri)\lambda_{\pri} - 1)\right]\mathrm{Ev}_{\varphi,1}(\phi_L)\notag\\
&= \lambda_{\ff}^{-1}\pi_{\ff}^{\jj+\vv}\left[\prod_{\pri \in B}\varphi_{p-\mathrm{fin}}(\pi_{\pri})(1 - \lambda_{\pri}^{-1}\varphi(\pri)^{-1})\right]\mathrm{Ev}_{\varphi,1}(\phi_L).
\end{align}
\end{mlem}

\begin{proof}
By definition, $\mu_\Psi \defeq \lambda_{\ff'}^{-1}\mu_\Psi^{\ff'}$. Hence we see that
\[\mu_{\Psi}(\varphi_{p-\mathrm{fin}}) = \lambda_{\ff'}^{-1}\sum_{\y\in\clgp{\ff'}}\epsphi\varphi_f(a_{\y})\mathrm{Ev}_{\ff',\dagger}^{a_{\y}}(\Psi)(\mathbf{z}^{\kk-\jj}).\]
Using the results of Section \ref{classical and ovcgt}, we can replace the overconvergent evaluations with classical ones, and then using the results of Section \ref{classical 2}, we get
\[\mu_{\Psi}(\varphi_{p-\mathrm{fin}}) = \lambda_{\ff'}^{-1}\pi_{\ff'}^{\jj+\vv}\sum_{\y\in\clgp{\ff'}}\epsphi\varphi_f(a_{\y})\mathrm{Ev}_{\ff',\jj,1}^{a_{\y}}(\phi_L).\]
We now use Corollary \ref{unramified primes}, which directly gives the first equality. The second equality follows since for $\pri$ not dividing $\ff$, we have $\pi_{\pri}^{\jj+\vv} = \varphi_{p-\mathrm{fin}}(\pi_{\pri})\varphi(\pri)^{-1},$ an identity which follows from the definition of $\varphi_{p-\mathrm{fin}}$.
\end{proof}
\begin{longversion}
\begin{mrem}
As an example, consider the case where $F=\Q$ and $\varphi = |\cdot|^j$ is the $j$th power of the norm map, which has trivial conductor. Then $\pi_p = p$, and $\varphi_{p-\mathrm{fin}}(p) = \varphi_p(p)p^j = p^{-j}p^j = 1$, whilst $1-\lambda_p^{-1}\varphi((p))^{-1} = 1-p^j/\lambda_p$. Thus in this case we get precisely the $p$-adic multiplier explored in Section 14 of \cite{MTT86}.
\end{mrem}
\end{longversion}

%
%

\section{Main results}\label{padiclfns}
\begin{longversion}The results of the previous section, and in particular Corollaries \ref{fincor} and \ref{fincor2} and Lemma \ref{interpolation unramified}, give the desired interpolation property for our distribution. The following is a summary of the main results of this paper. 
$\lb$
Recall the set-up. Let $F/\Q$ be a number field and $p$ a rational prime. Let $\Phi$ be a small slope cuspidal eigenform over $F$ of weight $\lambda = (\kk,\vv)\in \Z[\Sigma]^2$, where $\kk + 2\vv$ is parallel, and level $\Omega_1(\n)$, where $(p)|\n$. Let $\Lambda(\Phi,\cdot)$ be the normalised $L$-function attached to $\Phi$ in Definition \ref{lambda}. Write $\theta_L$ for the $p$-adic modular symbol associated to $\Phi$ in Definition \ref{modsymb}, where $L$ is a sufficiently large extension of $\Qp$. Using the control theorem, we may lift $\theta_L$ to a unique small slope overconvergent eigensymbol $\Psi$, and using Theorem \ref{construction} we may construct a distribution
\[\mu_{\Psi} \in \rayclgp\]
attached to $\Psi$.
\end{longversion}

\begin{shortversion}
The following is a summary of the main results of this paper. Recall the setting; let $\Phi$ be a small slope cuspidal eigenform for $\GLt$ over a number field $F$, of weight $\lambda = (\kk,\vv)\in \Z[\Sigma]^2$, where $\kk + 2\vv$ is parallel, and level $\Omega_1(\n)$, where $(p)|\n$. Let $\Lambda(\Phi,\cdot)$ be the normalised $L$-function attached to $\Phi$ in Definition \ref{lambda}. To $\Phi$, one can attach a unique overconvergent modular symbol $\Psi$ using Theorem \ref{comparaison}. Using Theorem \ref{construction} we may construct a distribution $\mu_{\Psi} \in \rayclgp$ attached to $\Psi$, which we defined to be the $p$-adic $L$-function of $\Phi$.
\end{shortversion}
\begin{mthm}\label{summary}
Let $\varphi$ be a Hecke character of conductor $\ff|(p^\infty)$ and infinity type $\jj+\vv$, where $0 \leq \jj \leq \kk$, and let $\varepsilon_\varphi$ be the character of $\{\pm 1\}^{\Sigma(\R)}$ attached to $\varphi$ in Section \ref{heckechar}. Let $\varphi_{p-\mathrm{fin}} \in \mathcal{A}(\clgp{p^\infty},L)$ be the $p$-adic avatar of $\varphi$. Let $B$ be the set of primes above $p$ that do not divide $\ff$. Then
\[\mu_{\Psi}(\varphi_{p-\mathrm{fin}}) = (-1)^{R(\jj,\kk)}\left[\frac{|D|\tau(\varphi)\pi_{\ff}^{\jj+\vv}}{2^{r_2}\lambda_{\ff}\Omega_{\Phi}^{\varepsilon_\varphi}}\right]\left(\prod_{\pri\in B}Z_{\pri}\right)\Lambda(\Phi,\varphi),\]
where 
\[Z_{\pri} \defeq \varphi_{p-\mathrm{fin}}(\pi_{\pri})(1 - \lambda_{\pri}^{-1}\varphi(\pri)^{-1})\]
(noting here that $\varphi(\pri)$ is well-defined since $\varphi$ is unramified at $\pri$).
\end{mthm}
Here $R(\jj,\kk) = \sum_{v\in\Sigma(\R)}j_v+k_v + \sum_{v\in\Sigma(\C)}k_v$, $D$ is the discriminant of $F$, $\tau(\varphi)$ is the Gauss sum of Definition \ref{gausssum}, $r_2$ is the number of pairs of complex embeddings of $F$, $\lambda_{\ff}$ is the $U_{\ff}$-eigenvalue of $\Phi$, $\Omega_{\Phi}^{\varepsilon_\varphi}$ is the fixed period attached to $\Phi$ and $\varepsilon_\varphi$ in Theorem \ref{deligneperiods}, and $\Lambda(\Phi,\cdot)$ is the normalised $L$-function of $\Phi$ as defined in Definition \ref{lambda}.

\section{Remarks on uniqueness}\label{admissibility}
When $F$ is a totally real or imaginary quadratic field, we can prove a \emph{uniqueness} property of this distribution. In particular, we prove that the distribution constructed above is \emph{admissible} in a certain sense, and any admissible distribution is uniquely determined by its values at functions coming from critical Hecke characters (see \cite{Colm10} and \cite{Loe14}). For further details of admissibility conditions in these cases, see \cite{Bar13} and \cite{Wil17} for the totally real and imaginary quadratic situations respectively. In the general case, things are more subtle. There is a good notion of admissibility for distributions on $\roi_F\otimes_{\Z}\Zp$, but it is not at all clear how this descends to a `useful' admissibility condition on $\clgp{p^\infty}$. 
$\lb$
In particular, recall that $\clgp{p^\infty} = \bigsqcup_{\nclgp} (\roi_F\otimes\Zp)^\times/\overline{E(1)}.$ When $F$ is imaginary quadratic, the unit group is finite, and in particular in passing to the quotient we do not change the rank. In this case, growth properties pass down almost unchanged. When $F$ is totally real, the unit group is in a sense `maximal' if we assume Leopoldt's conjecture. In particular, provided this, the quotient is just one dimensional, and we have a canonical `direction' with which to check growth properties.
$\lb$
Let us illustrate the difficulties of the general case with a conceptual example, for which the authors would like to thank David Loeffler. Let $F = \Q(\sqrt[3]{2})$, and note that $F$ is a cubic field of mixed signature. We see that $(\roi_F\otimes_{\Z}\Zp)^\times$ is a $p$-adic Lie group of rank 3, and that the quotient by $\overline{E(1)}$ has rank 2 (since the unit group has rank 1 by Dirichlet's unit theorem). In particular, a distribution on $\clgp{p^\infty}$ can `grow' in two independent directions. 
$\lb$
As the maximal CM subfield of $F$ is nothing but $\Q$, it follows that the only possible infinity types of Hecke characters of $F$ are parallel. In particular, there is only \emph{one} `dimension' of Hecke characters. In this sense, even though we have constructed a distribution that interpolates all critical Hecke characters, there are simply not enough Hecke characters to hope that we can uniquely determine a ray class distribution by this interpolation property.
$\lb$
One might be able to obtain nice growth properties using the extra structure that we obtain from our overconvergent modular symbol; in particular, one might expect the overconvergent cohomology classes we construct to take values in the smaller space of admissible distributions on $\roi_F\otimes\Zp$, which makes sense \emph{before} we quotient to obtain distributions on $\clgp{p^\infty}$. Without the theory of admissibility at hand in the latter situation, however, we cannot show that the distribution constructed in this paper is (in general) unique. We have tried to rectify this by proving that the distribution we obtain is independent of choices. As seen in the previous sections, we were able to do this up to a (fixed) choice of uniformisers at the primes above $p$. Hence, in the spirit of Pollack and Stevens in \cite{PS12}, we simply \emph{define} the $p$-adic $L$-function to be this distribution.
$\lb$
\begin{longversion}
It remains to comment on the dependence on choices of uniformisers. This dependence seems to be intrinsic to this more explicit approach; indeed, the evaluation maps at the level of $p$-adic coefficients depend on the choice of uniformiser, and accordingly the distribution we have defined to be the $p$-adic $L$-function does as well. However, the interpolation property \emph{also} has an explicit dependence on the uniformisers (coming from the Gauss sum and the term $\pi_{\ff}^{\jj+\vv}$), so by changing the uniformisers we are changing both the distribution and its interpolating property, so do not `break' any potential uniqueness property. Since submission, Bergdall and Hansen have given a similar construction in the Hilbert case which removes this dependency on uniformisers (see \cite{BH17}).
\end{longversion}
\begin{shortversion}
It remains to comment on the dependence on choices of uniformisers. Whilst this dependence seems intrinsic to our more explicit approach, since submission, Bergdall and Hansen have given a similar, but less hands-on, construction in the Hilbert case which removes this dependency on uniformisers (see \cite{BH17}).
\end{shortversion}

\small
\renewcommand{\refname}{\normalsize References} 
\bibliography{references}{}
\bibliographystyle{alpha}
\Addresses
\end{document}